\documentclass[10pt]{smfart}
\usepackage{a4wide}
\usepackage{amsmath,amssymb,smfthm}
\usepackage[francais]{babel}
\usepackage[latin1]{inputenc}
\usepackage[all]{xy}
\usepackage{graphicx}
\usepackage{epic,eepic}
\usepackage{epsfig}
\usepackage{color}

\def\R{\mathbb{R}}

\def\P{\mathbb{P}}
\def\C{\mathbb{C}}
\def\CC{\mathcal{C}}

\def\Fp{\mathbb{F}_p}

\def\L{\mathcal{L}}
\def\La{\Lambda}

\def\O{\mathcal{O}} 
\def\F{\mathcal{F}}

\def\*{^\times }
\def\dpt{\displaystyle}
\def\G{\mathcal{G}}

\def\a{\alpha}
\def\b{\beta}

\def\s{\sigma} 
\def\ph{\varphi}

\def\e{\epsilon}

\def\lssi{\Longleftrightarrow}

\def\limpl{\Longrightarrow}
\def\drt{\rightarrow}
\def\ldrt{\longrightarrow}

\def\Qp{\mathbb{Q}_p}
\def\Qpb{\overline{\mathbb{Q}}_p}
\def\Zp{\mathbb{Z}_p}
\def\Z{\mathbb{Z}}
\def\N{\mathbb{N}}

\def\Hom{\text{Hom}}

\def\Gal{\text{Gal}}
\def\={\! = \!}

\def\spec{\text{Spec}}
\def\spf{\text{Spf}}

\def\E{\mathcal{E}}

\def\limp{\underset{\longleftarrow}{\text{ lim }}\;}
\def\limi{\underset{\longrightarrow}{\text{ lim }}\;}

\def\iso{\xrightarrow{\;\sim\;}}

\def\GL{\hbox{GL}}

\def\xrig{\xrightarrow}

\def\M{\mathcal{M}}

\def\X{\mathfrak{X}}
\def\GG{\Gamma}

\def\bc{\backslash}
\def\AA{\mathcal{A}}
\def\spa{\hbox{Spa}}
\def\<{<\hspace{-1mm}}
\def\>{\hspace{-1mm}>}

\def\dem{{\it Démonstration. }}

\def\unp{[ \frac{1}{p}]}
\def\LT{\mathcal{L}\mathcal{T}}
\def\D{{\mathcal{D}r}}

\def\DD{\mathbb{D}}
\def\unpi{[\frac{1}{\pi}]}
\def\YY{\mathfrak{Y}}
\def\ZZ{\mathfrak{Z}}
\def\Ec{\mathcal{E}\! c}
\def\OEc{\mathcal{O}\!\mathcal{E}\! c}
\def\Ad{\mathcal{A} d}
\def\Spa{\text{Spa}}
\def\et{\text{ét}}
\def\top{^{\widetilde{\;\;}}}
\def\qet{\text{qét}}
\def\qetc{\text{qétc}}
\def\lss{\text{lss}}

\author{  Laurent Fargues} 
\address{CNRS-université Paris-Sud-IHES} 
\email{laurent.fargues@math.u-psud.fr} 
\date{}

\begin{document}

\title[Comparaison de la cohomologie des deux tours]{
 Comparaison de la cohomologie des tours de Lubin-Tate et de Drinfeld et correspondance de
  Jacquet-Langlands géométrique}
\maketitle

\begin{abstract}
Cet article est le dernier consacré  à l'isomorphisme entre les tours
de Lubin-Tate et de Drinfeld. Nous y démontrons l'existence d'un
isomorphisme entre la cohomologie à support compact des tours de Lubin-Tate et de
Drinfeld, et plus généralement de leurs complexes de cohomologie
équivariants. Nous y démontrons également l'existence d'une
correspondance de Jacquet-Langlands locale géométrique entre certains faisceaux
étales rigides équivariants sur l'espace des périodes de Gross-Hopkins $\mathbb{P}^{n-1}$
et sur l'espace de Drinfeld $\Omega$.
\end{abstract}

\begin{altabstract}
This article is the last one about the isomorphism between Lubin-Tate
and Drinfeld towers. We  prove the existence of an isomorphism
between the compactly supported cohomology of the Lubin-Tate and
Drinfeld towers, and more generally their equivariant cohomology
complex. 
We also prove the existence of a geometric local Jacquet-Langlands
correspondence between some equivariant rigid étale sheaves on
Gross-Hopkins period space $\mathbb{P}^{n-1}$ and Drinfeld one $\Omega$. 
\end{altabstract}

\section*{Introduction}

Cet article est le dernier consacré à l'isomorphisme entre les tours
de Lubin-Tate et de Drinfeld. Les  principaux résultats sont :
\begin{itemize}
 \item Soit $(X_i)_{i\in I}$  un système projectif filtrant d'espaces rigides
  quasicompacts. Soit $(\X_i)_{i\in I}$  un modèle entier de ce
  système projectif formé de schémas formels admissibles au sens de
  Raynaud dont les morphismes de transition sont affines. Soit
  $\X_\infty =\underset{i}{\limp}{\X_i}$, par exemple si $\X_i=\spf
  (R_i)$ alors $\X_\infty =\spf ((\underset{i}{\limi}
  R_i)^{\widehat{\;\;\;}})$. Alors le topos limite projective 
$
\underset{i}{\limp} (X_i)_\et
$
ne dépend que de $\X_\infty$ et la cohomologie $\underset{i}{\limi}
H^\bullet ((X_i)_\et, \La)$ également. 
\item L'existence d'une correspondance ``de Jacquet-Langlands locale
  géométrique'' entre faisceaux étales $D^\times$-équivariants sur
  l'espace des périodes de Gross-Hopkins $(\P^{n-1})^{rig}$ pour
  lesquels l'action de $D^\times$ est lisse (au sens de la théorie des
  représentations des groupes $p$-adiques : le stabilisateur d'une
  section est ouvert) et faisceaux $\GL_n (F)$-équivariants sur
  l'espace de Drinfeld $\Omega$ pour lesquels l'action de $\GL_n (F)$
  est lisse. Il s'agit du théorème \ref{QFKSDuqza236gjyza}. 
\item Le fait que les complexes de cohomologie à support compact de la
  tour de Lubin-Tate et de Drinfeld, vus comme éléments de la
  catégorie dérivée $\GL_n (F)\times D^\times \times W_F$-équivariante
  lisse, sont isomorphes. Il s'agit du théorème \ref{ksdsrpytg79324Fgp}
\end{itemize}

Les deux derniers résultats utilisent bien sûr le théorème principal de \cite{iso4} ainsi que le premier résultat.
Pour le premier résultat
les points clef consistent en une utilisation du théorème d'approximation d'Elkik (\cite{Elkik}) et du théorème de platification de Raynaud-Gruson (\cite{Ray2}). 

\subsection*{La correspondance de Jacquet-Langlands locale
  géométrique}

Rappelons que la tour de Lubin-Tate est une tour d'espaces analytiques
rigides 
$
(\M_K^{\LT})_{K\subset GL_n (\O_F)}
$
munie d'une action ``verticale'' de $\GL_n (F)$ par correspondances de
Hecke et ``horizontale'', étage par étage, de $D^\times$. Elle forme un ``pro-$\GL_n
(\O_F)$-torseur étale'' au dessus de la tour sans niveau 
$$
\M^{\LT}_{GL_n (\O_F)} = ``\text{Tour de L.T.}/\GL_n (\O_F) `` \simeq
\coprod_{\Z} \overset{\hspace{-6mm}\circ}{\mathbb{B}^{n-1}}
$$
une union disjointe de boules $p$-adiques ouvertes. De plus il y a une
application des périodes de Gross-Hopkins (\cite{HopkinsGross}, \cite{RZ}, \cite{Cellulaire})
$$
\M^{\LT}_{GL_n (\O_F)} \ldrt \P^{n-1}
$$
surjective, $D^\times$-équivariante, invariante sous les
correspondances de Hecke sphériques et dont les fibres sont exactement
les orbites de Hecke. On a donc envie d'écrire 
$$
\text{Tour de L.T.}/\GL_n (F)\simeq \P^{n-1}
$$
ou plutôt au sens des champs $\left [  \text{Tour de L.T.}/\GL_n (F)
\right ]\simeq \P^{n-1}$ i.e. les faisceaux rigides-étales
Hecke-équivariants sur la tour de Lubin-Tate sont en bijection avec
les faisceaux rigides-étales sur $(\P^{n-1})^{rig}$. 
\\

Quant à la tour de Drinfeld $(\M_K^\D)_{K\subset \O_D^\times}$ elle
est munie d'une action ``verticale'' de $D^\times$ et horizontale de
$\GL_n (F)$. Elle forme un ``pro-$\O_D^\times$-torseur étale'' au
dessus de la tour sans niveau
$$
\M^\D_{\O_D^\times} = ``\text{Tour de Dr.}/\O_D^\times `` \simeq
\coprod_{\Z} \Omega
$$
où $\Omega$ désigne l'espace de Drinfeld, si $n=2\;\; \Omega (\C_p)
=\C_p\setminus \Qp$. Il y a de plus une application des périodes
beaucoup plus simple que la précédente
$$
\M^\D_{\O_D^\times} = \coprod_{\Z} \Omega \ldrt \Omega
$$
qui est $\GL_n (F)$-équivariante, un isomorphisme sur chacune des
composantes de l'union $\coprod_\Z$, $D^\times$-invariante et dont les
fibres sont exactement les $D^\times/\O_{D}^\times$-orbites. On a donc
envie d'écrire 
$$
\left [ \text{Tour de Dr.}/D^\times \right ] \simeq \Omega
$$
c'est à dire que les faisceaux rigides-étales $D^\times$-équivariants
sur la tour de Drinfeld sont en bijection avec les faisceaux
rigides-étales sur $\Omega$.
\\

Rappelons maintenant (\cite{iso4}) qu'on a démontré l'existence d'un 
``isomorphisme en niveau infini'' 
$$
\M_\infty^{\LT}\iso \M_\infty^{\D}
$$
$\GL_n (F)\times D^\times$-équivariant. On a donc envie d'écrire
qu'il y a un isomorphisme
$$
\left [ D^\times \bc \P^{n-1} \right ] \simeq \left [ \GL_n (F)\times
  D^\times \bc \M_\infty^{\LT} \right ] \simeq \left [ \GL_n (F)\times
  D^\times \bc \M_\infty^\D \right ] \simeq \left [ \GL_n (F)\bc
  \Omega \right ]
$$
c'est à dire une équivalence entre faisceaux $D^\times$-équivariants
sur $\P^{n-1}$ et faisceaux $\GL_n (F)$-équivariants sur $\Omega$. 

Cet énoncé est en fait incorrect. Pour expliquer pourquoi prenons une
analogie. Soit $k$ un corps de clôture séparable $\overline{k}$ et $X$
un $k$-schéma de type fini. Soit $p:X_{\overline{k}}\ldrt X$ la
projection. Le foncteur qui à un faisceau étale $\F$ sur $X$ associe
le faisceau $p^*\F$ muni de son action de $\Gal (\overline{k}|k)$
compatible à celle sur $X_{\overline{k}}$
induit une équivalence entre faisceaux étale sur $X$ et ceux sur
$X_{\overline{k}}$ munis d'une action de $\Gal (\overline{k}|k)$
compatible à celle sur $X_{\overline{k}}$ et continue. La condition de
continuité signifie que le stabilisateur d'une section du faisceau sur
un ouvert quasicompact est un sous-groupe ouvert de $\Gal
(\overline{k}|k)$. Dit en d'autres termes, il y a une condition de
continuité sur la donnée de descente pour descendre des objets en
niveau infini sur la tour $(X\otimes_k L)_{L|k\text{ finie}}$ à $X$.

Il en est de même pour les tours de Lubin-Tate et de Drinfeld. Il y a
bien une équivalence entre faisceaux équivariants sur $\P^{n-1}$ et
$\Omega$ mais il faut ajouter une condition de continuité sur l'action
: le stabilisateur d'une section est un sous-groupe compact-ouvert.
\\

Pour expliquer cette condition de continuité revenons à notre
analogie. Si $U\ldrt X_{\overline{k}}$ est un morphisme étale avec $U$
quasicompact il existe une extension de degré fini $L|k$ et un ouvert
étale $U'\ldrt X_L$ tel que $U'=U\otimes_L \overline{k}$. Ainsi en
restriction au sous-groupe ouvert $\Gal (\overline{k}|L)$ sur
$X_{\overline{k}}$ s'étend en une action sur $U$
$$
\forall \s\in \Gal (\overline{k}|L) \;\;
\xymatrix@R=5mm@C=6mm{
U \ar[r]\ar[d] & U \ar[d] \\
X_{\overline{k}} \ar[r]^{\s^*} & X_{\overline{k}} 
}
$$
et si $\F$ est un faisceau étale sur $X_{\overline{k}}$ muni d'une
action de $\Gal (\overline{k}|k)$, $\forall  \s\in \Gal
(\overline{k}|L)\;\; (\s^*\F) (U) \iso \F (U)$ et il y a une action de
$
\Gal (\overline{k}|L)$ sur $\F (U)$. La condition de lissité de
l'action de $\Gal (\overline{k}|L)$ est donc bien définie.

L'analogue en géométrie analytique a été étudié par Berkovich dans
\cite{Berk2}. Expliquons le plutôt dans le contexte des espaces
rigides. Soit $Y=\text{Sp} (\mathcal{B})\ldrt \text{Sp} (\mathcal{A})$
un morphisme étale entre espaces affinoïdes. Soit $G$ un groupe
profini agissant continûment sur $X$ au sens où si $f\in \mathcal{A}$ 
$$
\underset{g\drt e}{\lim}\; \|g^* f-f\|_\infty =0
$$
Il existe alors un moyen de relever ``canoniquement'' l'action de $G$
sur $X$ à $Y$, quitte à se restreindre à un sous-groupe ouvert
suffisamment petit $U$ de $G$
$$
\xymatrix@R=6mm@C=7mm{
Y \ar[rr]  \ar[d] & \ar@(lu,ru)[]|U & Y \ar[d] \\
X  \ar[rr] & \ar@(lu,ru)[]|U & X
}
$$
Ainsi si $\F$ est un $G$-faisceau étale sur $X$ il y a une action de
$U$ sur $\F (Y)$ et on peut parler de lissité de l'action de $U$ sur
$\F (Y)$. 

Ce fait est une généralisation du lemme de Krasner. Si $\YY=\spf (A)\ldrt \spf
(B)=\X$ est un modèle entier du morphisme $Y\ldrt X$, i.e. $A$ et $B$
sont deux algèbre topologiquement de type fini sur $\Zp$ sans $p$-torsion, $A\unp = \mathcal{A}$
et
$B\unp = \mathcal{B}$ il existe alors un entier $N$ tel que pour tout
$\X$-schéma formel topologiquement de type fini $\ZZ$ l'application de
réduction
$$
\Hom_\X ( \ZZ, \YY) \ldrt \Hom_{\X\otimes \Z/p^N} ( \ZZ\otimes \Z/p^N,
\YY\otimes \Z/p^N)
$$
soit une bijection. Il s'agit essentiellement d'une application du théorème
d'approximation d'Elkik (\cite{Elkik}). 
Si le morphisme entier $\YY\ldrt \X$ est étale on peut bien sûr
prendre $N=1$ et en général l'entier $N$  dépend d'une puissance
suffisamment grande de $p$ telle qu'un certain idéal discriminant
divise cette puissance (entier qui existe puisque après inversion de
$p$ le morphisme est étale). 
En particulier deux $\X$-morphismes
proches de $\YY$ dans lui-même coïncident et tout isomorphisme défini
modulo $p^N$ peut se relever en caractéristique zéro.
\\

Finissons par expliquer la terminologie ``correspondance de Jacquet Langlands''. Soit $\F$ un $D^\times$-faisceau lisse sur $\P^{n-1}$ et $\text{JL} (\F)$ le faisceau associé sur $\Omega$.
Soit $\delta\in D^\times$ semi-simple régulier et $\gamma\in \GL_n
(F)$ stablement conjugué à $\delta$. Il devrait y avoir un lien entre
$\dpt{
\sum_{x\in \text{Fix} (\delta ; \P^{n-1})} \text{tr} ( \delta ; \F_x
)}$ et $\dpt{\sum_{y\in \text{Fix} (\gamma ;\Omega)} \text{tr} ( \gamma ; \text{JL} (\F)_y)}
$. 
Cela n'est pas démontré dans cet article, l'auteur prévoit d'en écrire la démonstration ainsi que d'autres propriétés de cette correspondance dans un futur article.

\subsection*{Cohomologie à support compact des deux tours}

On démontre l'existence d'isomorphismes
$$
\underset{K\subset GL_n (\O_F)}{\limi} H^\bullet_c (\M_K^{\LT}\hat{\otimes}
\C_p,\Z/\ell^n \Z) \simeq \underset{K\subset \O_D^\times}{\limi} 
H^\bullet_c (\M_K^\D\hat{\otimes}
\C_p,\Z/\ell^n \Z)
$$
en tant que $\Z/\ell^n\Z [\GL_n (F)\times D^\times \times
W_F]$-modules lisses. En fait on démontre un résultat beaucoup plus
précis au niveau des complexes de cohomologie à support compact dans
la catégorie dérivée équivariante-lisse. Si $\F$ est un
$D^\times$-faisceau lisse sur $\P^{n-1}$ on lui associe canoniquement,
en tirant en arrière ce faisceau par l'application des périodes en
chaque niveau de l'espace de Lubin-Tate, 
un complexe de cohomologie à support compact de la tour de Lubin-Tate 
$$
R\GG_c (\LT, \F) \in \DD^+ ( \Z/\ell^n \Z [\GL_n (F)\times D^\times
\times W_F]-\text{Mod-lisses})
$$
De même si $\G$ est un $\GL_n (F)$-faisceau lisse sur $\Omega$ on lui
associe $R\GG_c (\D,\G)$. Soit alors $\F\longmapsto \text{JL} (\F)$
la correspondance de Jacquet-Langlands décrite précédemment. Il y a
alors un isomorphisme de foncteurs
$$
R\GG_c (\LT,-) \iso R\GG_c (\D,-)\circ \text{JL}
$$

\subsection*{Description des différentes parties}

\subsubsection*{Espaces rigides généralisés}

Dans les chapitres 1 à 3 on définit et étudie les propriétés de base
des espaces rigides ``fibre générique'' de schémas formels $p$-adiques qui ne
sont pas nécessairement topologiquement de type fini, du type de ceux
intervenant dans \cite{Cellulaire} et \cite{iso4}. Le point de vue
choisi est celui de Raynaud consistant à voir la catégorie des espaces rigides 
comme un localisé de la catégorie des schémas formels $p$-adiques
relativement aux éclatements formels admissibles.
\begin{itemize}
\item Dans le premier chapitre on définit et étudie les éclatements
  formels admissibles de tels schémas formels. 
Il s'agit
  essentiellement de vérifier que certaines notions intervenant dans
  \cite{BLI} restent valables dans un contexte plus général. Une des
  difficultés de la théorie est que contrairement au cas classique de
  \cite{BLI} les théorèmes de cohérence des algèbres $p$-adiques
  topologiquement de type fini sans $p$-torsion ne sont plus
  valables. 

On retiendra une des propriétés importantes de ces éclatements formels
admissibles généralisés n'ayant pas d'équivalent dans le cadre
``classique'' : la proposition \ref{vuindyuyszpgyt} assurant la
compatibilité de ces éclatements au passage à la limite projective de
schémas formels.
\item Dans le deuxième chapitre on définit et étudie le topos
  admissible de tels schémas formels. On utilise pour cela le langage
  des limites inductives de sites et des limites projectives de topos
  de \cite{SGA4_exp6}. Le topos admissible est vu comme la limite
  projective des topos Zariskiens des éclatements de notre modèle
  formel.

L'une des propriétés fondamentales est le théorème de
``décompletion'', la proposition \ref{vhkqpuzg} qui dit qu'un ouvert
admissible d'une limite projective de schémas formels
$\underset{i}{\limp}\X_i$ est un germe d'ouvert admissible sur les
$\X_i$, et la proposition \ref{jfyyusouzrpgoie53TU5} son
interprétation en termes de topos.

Plus tard on appliquera le même type de procédures pour le site étale
rigide au lieu du site admissible. Le lecteur peut donc considérer ce
chapitre comme un entraînement à la manipulation des topos limites
projectives dans un cas ``simple''. 
\item Dans le chapitre 3 on définit et étudie l'espace de
  Zariski-Riemann associé à nos espaces rigides généralisés. Pour les
  espaces rigides classiques cet espace coïncide avec l'espace
  topologique adique défini par Huber (\cite{Hu1}) et également étudié
  par Fujiwara (\cite{Fuji2}). On montre en particulier que le topos
  admissible s'identifie au topos des faisceaux sur cet espace
  topologique. 
On interprète également cet espace comme un recollé de spectres
valuatifs. 

On utilisera plus tard l'espace de Zariski-Riemann dans le but de
définir la notion de famille couvrante dans notre site étale rigide. 

On définit également l'espace de Berkovich associé aux générisations
maximales dans l'espace de Zariski-Riemann (les valuations de rang 1
i.e. à valeurs dans $\R$). Cela nous permettra de définir la notion de
topos surconvergent correspondant dans le cas classique au topos étale
de Berkovich.
\end{itemize}

\subsubsection*{Le topos étale rigide}

Dans les chapitres 4,5 et 6 on définit et étudie le topos étale rigide
de nos schémas formels. Il s'agit du coeur de l'article.

Expliquons d'abord les difficultés rencontrées.
\begin{itemize}
\item
Soit $A$ une algèbre $p$-adique sans $p$-torsion. Si $B$ est une
$A$-algèbre $p$-adique topologiquement de présentation finie et si 
$B\simeq A<T_1,\dots,T_N>/(f_1,\dots,f_q)$ est une présentation on
peut définir comme dans le chapitre 4 de \cite{Elkik} 
en utilisant un idéal jacobien explicite 
la notion d'être
rig-étale après inversion de $p$ :
une certaine puissance de $p$ appartient à cet idéal jacobien. 
 Néanmoins il n'y a pas de raison en
général pour que cette définition ne dépende pas du choix d'une telle
présentation comme c'est le cas si $A$ est topologiquement de type
fini sur un anneau de valuation de hauteur $1$ (le cas rigide
classique). 
\item
On veut de plus appliquer le théorème 6 de \cite{Elkik} à nos algèbres
rig-étales. Or, même après avoir fixé une présentation comme
précédemment la démonstration donnée dans \cite{Elkik} ne s'adapte pas
sans des théorèmes de cohérence. Essentiellement le résultat qui pose
problème est que si $R$ est une algèbre $p$-adique et $I=\{x\in R\;|\;
\exists k\; p^k x=0\}$ alors en général $I$ n'est pas un idéal fermé de $R$,
et si l'on veut tuer la $p$-torsion dans $R$ il ne faut pas prendre
$R/I$ mais $R/\overline{I}$, le séparé de $R/I$. En général on n'a aucun contrôle sur $\bar{I}$. 
Bien sûr si $R$ est topologiquement de type fini sur un anneau de valuation de hauteur $1$ d'après \cite{BLI} on a $\overline{I}=I$, mais dans notre situation il n'y a pas de raison pour que ce soit le cas.
\item Soit $(R_i)_{i\in \N}$ un système inductif d'algèbres $p$-adiques topologiquement de type fini sur $\Zp$, $R_\infty = \underset{i}{\limi} R_i$. On peut même supposer que $R_i\hookrightarrow R_{i+1}$, $R_{i+1}|R_i$ est fini, $R_i$ est intégralement fermé dans $R_i\unp$, 
$R_i\unp\hookrightarrow R_{i+1}\unp$ est étale fini et que le système inductif $(R_i \unp)_{i\in \N}$ est un torseur étale sous un groupe profini et les $R_i \unpi$ sont rig-lisses sur $\Zp$ 
(exemple : prendre l'image réciproque d'un ouvert affine dans les tours de Lubin-Tate et de Drinfeld). 

On veut qu'il y  ait une équivalence entre $\widehat{R}_\infty$-algèbres $p$-adiques rig-étales et ``germes'' d'algèbres $p$-adiques rig-étales sur les $(R_i)_{i}$, via l'application 
$$(R_i\ldrt B)\longmapsto (\widehat{R}_\infty \ldrt B\hat{\otimes}_{R_i} \widehat{R}_\infty )$$
Il devrait s'agir typiquement d'une application du théorème d'approximation d'Elkik à l'anneau Hensélien $(R_\infty, pR_\infty)$. Néanmoins il n'y a pas de raison dans le cas des tours de Lubin-Tate et de Drinfeld, pour les modèles construits dans \cite{Cellulaire} et \cite{iso4} par normalisation,  pour que les morphismes $R_i\ldrt R_{i+1}$ soient plats. Il n'y  a pas non plus de raison pour que le système inductif $(R_i)_{i\in \N}$ satisfasse à une hypothèse de presque platitude au sens de Faltings, c'est à dire si $i_0\in \N$, $I\subset R_{i_0}$ est un idéal de type fini alors
$$
\exists \a\in\N\;\forall i\geq i_0\;\; p^\a.\text{Tor}_1^{R_{i_0}} ( R_i, R_{i_0}/I)=0
$$
Par exemple on aurait bien aimé que le morphisme $R_\infty \ldrt \widehat{R}_\infty$ soit fidèlement plat, mais il n'y a pas de raison non plus pour que ce soit le cas. 

Le problème que cela soulève est que si $B$ est une $R_i$-algèbre rig-étale sans $p$-torsion alors la $p^\infty$-torsion des algèbres $(B\hat{\otimes}_{R_i} R_j)_{j\geq i}$ peut ``exploser'' lorsque $j\ldrt +\infty$. 
\end{itemize}
\vspace{4mm}

Voici comment on procède pour palier à ces problèmes.
\begin{itemize}
\item Dans le chapitre 4 on montre que si $A$ est une algèbre $p$-adique sans $p$-torsion et que l'on se restreint aux $A$-algèbres $p$-adiques finies localement libres étales après inversion de $p$ on a une bonne théorie pour de telles algèbres rig-étales. Ainsi si $\X$ est un schéma formel $p$-adiques sans $p$-torsion quelconque on a une bonne notion de $\X$-schéma formel fini localement libre rig-étale. 

On montre de plus que le théorème d'approximation d'Elkik s'applique à
ces morphismes rig-étales : si $(\X_i)_{i\in \N}$ est un système
projectif de schémas formels $p$-adiques sans $p$-torsion
à morphismes de transition affines, $\X_\infty =\underset{i\in \N}{\limp} \X_i$, les $\X_\infty$-schémas formels finis localement libres rig-étales sont équivalents aux ``germes'' de $\X_i$-schémas formels finis localement libres rig-étales lorsque $i$ varie via
$$
(\YY\ldrt \X_i)\longmapsto (\YY\times_{\X_i} \X_\infty \ldrt \X_\infty )
$$
\item Dans le chapitre 5 on définit et étudie la classe générale de morphismes rig-étales utilisée pour définir le topos étale rigide. On les appelles morphismes de type $(\E)$. Il sont construits à partir d'éclatements formels admissibles, de morphismes étales de schémas formels et de morphismes finis localement libres rig-étales (ceux étudiés dans le chapitre 4). On montre qu'ils vérifient le théorème d'approximation comme précédemment : les $\X_\infty$-schémas formels de type $(\E)$ sont équivalents aux germes de $\X_i$-schémas formels de type $(\E)$ lorsque $i$ varie.

Le point clef justifiant leur introduction est que d'après le théorème de platification de Raynaud-Gruson si $\YY\ldrt \X$ est un morphisme de schémas formels admissibles quasicompacts tel que $\YY^{rig}\ldrt \X^{rig}$ soit étale il existe alors un diagramme
$$
\xymatrix@R=5mm@C=5mm{
\ZZ \ar[rr]\ar[rd] & & \YY\ar[ld] \\
 & \X
}
$$
où $\ZZ\ldrt \X$ est de type $(\E)$ et $|\ZZ^{rig}|\ldrt |\YY^{rig}|$ est surjectif et que donc les $\X$-schémas formels de type $(\E)$ engendrent topologiquement le site étale de $\X^{rig}$.
\item Dans le chapitre 6 on définit et étudie un site ainsi que le topos étale rigide des schémas formels ne vérifiant pas d'hypothèse de finitude. Pour cela on utilise les morphisme de type $(\E)$ du chapitre 5. 

On a vu que dans le cas classique ces morphismes engendrent topologiquement le site étale rigide usuel. D'après le théorème 4.1 de l'exposé III de SGA4 (cf. théorème \ref{odpmmlfjhjutrr} de cet article) ils suffisent pour reconstruire le topos étale d'un espace rigide classique. La contrepartie est qu'ils ne forment pas une prétopologie de Grothendieck à cause de l'absence de certains produits fibrés dans cette catégorie.

On utilise alors toute la puissance du formalisme de Grothendieck qui permet de définir une topologie en toutes généralités sans cette hypothèse d'existence de produits fibrés. L'analogie suivante éclairera peut-être le lecteur : on est dans la situation où l'on a un espace topologique $X$ muni d'une famille d'ouverts $\mathcal{C}$ telle que tout ouvert de $X$ puisse s'écrire comme une union d'ouverts de $\mathcal{C}$, mais si $U,V\in\mathcal{C}$ $\; U\cap V$ n'appartient pas nécessairement à $\mathcal{C}$. Néanmoins, grâce au formalisme des cribles, $\mathcal{C}$ est muni d'une topologie de Grothendieck qui n'est pas définie par une prétopologie mais qui permet tout de même de retrouver le topos $X\top$ comme équivalent à $\mathcal{C}\top$; un faisceau sur $X$ est la même chose qu'un foncteur, i.e. un préfaisceau, $\F$ défini sur $\mathcal{C}$ vérifiant $\forall U\in \mathcal{C} \; \forall (V_i)_{i\in I}$ un recouvrement de $U$ par des objets de $\mathcal{C}$
$$
\F (U) = \ker ( \prod_{i\in I} \F (V_i) \xymatrix@C=6mm{\ar@<.6ex>[r] \ar@<-.6ex>[r] &} \prod_{i,j} \Hom_{\text{préfaisceaux}} ( h_{V_i}\times_{h_U} h_{V_j},\F)  )
$$
où $h_V$ désigne le préfaisceau représenté par $V$ et la même formule permet d'étendre $\F$ à tout ouvert de $X$ qui n'est pas dans $\mathcal{C}$. 

Le théorème principal s'énonce alors en disant que le topos rig-étale d'une limite projective de schémas formels $p$-adiques est la limite projective des topos rig-étales de chaque schéma formel de la limite projective. Il s'agit d'une application des résultats d'approximation des chapitres précédents.
\end{itemize}

\subsubsection*{Faisceaux étales munis d'une action lisse d'un groupe
  $p$-adique}

Dans les chapitres 8 à 10 on s'intéresse à la cohomologie à support
compact équivariante des espaces analytiques de Berkovich ainsi que
des espaces rigides généralisés introduits dans les section
précédentes.

On s'inspire des travaux de Berkovich (\cite{Berk2} et
\cite{BerkNotes}) sur le sujet. \'Etant donné un groupe topologique
$G$ agissant continûment sur un espace analytique de Berkovich $X$,
resp. un espace rigide généralisé, étant donné un faisceau étale $\F$
muni d'une action lisse de $G$ compatible à celle sur $X$ on montre
que l'action de $G$ sur $H^\bullet_c (X,\F)$ est lisse.  En fait le
but plus général est d'associer fonctoriellement à $\F$ un complexe de
cohomologie $R\GG_c (X,\F) \in \DD^+ ( \La [G]-\text{Mod-lisses})$.

\begin{itemize}
\item Pour cela on développe dans le chapitre 8 un formalisme général
  des $G$-faisceaux lisses qui peut s'appliquer aussi bien aux espaces
  de Berkovich qu'aux espaces rigides généralisés. Le résultat
  principal est le théorème \ref{OOPfyzY136Uetet} qui dit en
  particulier qu'on peut résoudre un $G$-faisceau lisse par des
  $G$-faisceaux lisses flasques après oubli de l'action de $G$.
\item Le chapitre 9 contient l'application de ce formalisme aux
  espaces analytiques de Berkovich. Le site étale d'un espace de
  Berkovich n'est pas adapté au formalisme des faisceaux lisses, de
  plus plus tard dans l'article nous devrons jongler entre site étale
  d'un espace de Berkovich et site étale de l'espace rigide associé,
  c'est pourquoi nous travaillons avec les sites quasi-étales définis
  dans \cite{Berk2} (le site quasi-étale correspond au site étale de
  l'espace rigide et le site étale au site étale surconvergent). Par
  contre la cohomologie à support compact des espaces de Berkovich
  sans bord (i.e. surconvergents sur leur corps de base) est plutôt
  bien adaptée au site étale. D'où les jonglages permanents entre
  sites étales et quasi-étales.
\item Le chapitre 10 est consacré à l'analogue pour les espaces
  rigides généralisés.
\end{itemize}

\subsubsection*{Les résultats principaux}

Les chapitres 11 et 12 sont consacrés à la démonstration des
principaux résultats de l'article. On y récolte les résultats des
chapitres précédents pour les appliquer aux tours de Lubin-Tate et de
Drinfeld.

\vspace{1cm}

{\it Remerciements : L'auteur tient à remercier Jean François Dat dont
les travaux et les diverses discussions que l'auteur a eues avec lui ont
influencé la forme de cet article, notamment le chapitre 11.
Il remercie également Ofer Gabber pour des discussions sur le sujet.  
 Il apparaîtra comme
clair au lecteur que les travaux de Raynaud en géométrie rigide ont eu
une forte influence sur l'auteur ainsi que l'article \cite{Fuji2} de
Fujiwara dont l'auteur s'est largement inspiré.
}

\section{Schémas formels $\pi$-adiques}

On fixe $\O_K$ un anneau de valuation de hauteur $1$ et $\pi$ un
élément de $\O_K$ de valuation strictement positive. On vérifiera
facilement que les définitions données ne dépendent pas du choix d'un
tel $\pi$. 
Tous les schémas formels considérés seront supposés quasi-séparés.

\subsection{Rappels sur les schémas formels $\pi$-adiques}

\begin{defi}
On appelle schéma formel $\pi$-adique (sous-entendu sur $\spf (\O_K)$) un schéma formel $\mathfrak{Z}$ sur $\spf (\O_K)$ tel que $\pi\O_\mathfrak{Z}$ soit un idéal de définition de $\mathfrak{Z}$. 
\end{defi}

Soit 
la   catégorie dont les objets sont  $(\spec (\O_K/\pi^k \O_K))_{k\geq 1}$
et les flèches sont les morphismes de réduction modulo des puissances de $\pi $, $ \spec (\O_K/\pi^k \O_K)\hookrightarrow \spec (\O_K/\pi^l \O_K)$ pour $k\leq l$ :
$$
\spec (\O_K/\pi\O_K) \drt \dots \drt \spec  (\O_K/\pi^k \O_K) \drt 
\spec (\O_K/\pi^{k+1}\O_K) \drt \dots
$$
On définit une catégorie fibrée au dessus de cette petite catégorie en posant  que la fibre sur $\spec (\O_K/\pi^k \O_K)$ est la catégorie des schémas sur $\spec (\O_K/\pi^k\O_K)$. La catégorie des schémas formels $\pi$-adiques est alors équivalente à la limite projective de cette catégorie fibrée 
 (cf. \cite{SGA4} exposé 6 section 6.10 page 273 pour la notion de limite projective de catégories fibrées). 

 Cela signifie que se donner un schéma formel $\pi$-adique est équivalent à se donner une famille $(Z_k)_{k\geq 1}$ où $Z_k$ est un $\spec (\O_K/\pi^k\O_K)$-schéma munie d'isomorphismes $Z_{k+1}\otimes \O_K /\pi^k \O_K\iso Z_k$ satisfaisant une condition de cocyle évidente (un schéma formel $\pi$-adique n'est rien d'autre qu'un cas particulier d'ind-schéma). Le schéma  formel associé à une telle famille $(Z_k)_k$ sera noté $\underset{k}{\limi} Z_k$. 
 
\begin{defi}
 On dit que $\mathfrak{Z}$ est sans $\pi$-torsion si le faisceau
 $\O_{\mathfrak{Z}}$ l'est c'est à dire $\O_\mathfrak{Z}\xrig{\;\times
   \pi\;} \O_\mathfrak{Z}$ est un monomorphisme. 
\end{defi} 

Ainsi $\ZZ$ est sans $\pi$-torsion ssi $\forall k\in \N$ le schéma
$\ZZ\otimes \O_K/\pi^k\O_K$ est plat sur $\spec (\O_K/\pi^k \O_K)$. 
On en déduit aussitôt que 
le schéma formel $\pi$-adique $\mathfrak{Z}$ est sans $\pi$-torsion ssi
il possède un
recouvrement affine $(\spf (R_i))_{i\in I}$ tel que $\forall i\; R_i$
soit sans $\pi$-torsion. De plus on en déduit que si $X$ est un schéma
plat sur $\spec (\O_K)$ alors son complété $\pi$-adique est sans
$\pi$-torsion. 

\begin{defi}
Une algèbre $\pi$-adique est une $\O_K$-algèbre séparée complète pour
la topologie $\pi$-adique.
\end{defi}

Ainsi la catégorie des algèbre $\pi$-adiques est équivalente à celle
des schémas formels $\pi$-adiques affines. 

\subsection{Morphismes topologiquement de type fini}

\begin{defi}
Un morphisme de schémas formels $\pi$ adiques $f:\mathfrak{Y}\ldrt
\mathfrak{X}$ est dit localement topologiquement de type fini si le
morphisme induit $\X\otimes \O_K/\pi\O_K\ldrt \YY\otimes \O_K/\pi\O_K$
est localement de type fini. Il est dit topologiquement de type fini
s'il est de plus quasicompact.
\end{defi}

\begin{lemm}
Soit $f:\YY\ldrt \X$. Sont équivalents 
\begin{itemize}
\item $f$ est localement topologiquement de type fini
\item Pour tous ouverts affines $\spf (B)\subset \YY$, $\spf (A)\subset \X$, tels que 
$f (\spf (B))\subset \spf (A)$ le morphisme induit $A\ldrt B$ fait de $B$ une $A$-algèbre topologiquement de type fini c'est à dire isomorphe en tant que $A$-algèbre topologique à $A<T_1,\dots,T_n>/I$ où $A<T_1,\dots, T_n>$ désigne l'anneau des séries formelles strictes,
le complété $\pi$-adique de $A[T_1,\dots,T_n]$, et $I$ est un idéal fermé de $A<T_1,\dots,T_n>$. 
\end{itemize}
\end{lemm}

Les schémas formels $\pi$-adiques localement topologiquement de type fini sur $\spf (\O_K)$ sont étudiés en détails dans \cite{BLI} auquel on renvoie le lecteur. 

\begin{defi}[\cite{BLI}]
Les schémas formels $\pi$-adiques localement topologiquement de type fini sur $\spf
(\O_K)$ sans $\pi$-torsion sont appelés schémas formels admissibles.  
\end{defi}

\subsection{Morphismes topologiquement de présentation finie}

La notion qui suit n'a d'intérêt que pour les schémas formels sans $\pi$-torsion. 

\begin{defi}
Soit $f:\X\ldrt \mathfrak{Y}$ un morphisme entre schémas formels
$\pi$-adiques sans $\pi$-torsion. Il sera dit localement
topologiquement de présentation finie si le morphisme de schémas
induit $\X\otimes \O_K/\pi\O_K\ldrt \mathfrak{Y}\otimes \O_K/\pi\O_K$
est localement de présentation finie. Il sera dit topologiquement de présentation finie
s'il est de plus quasicompact.
\end{defi}

\begin{lemm}
Soit $f:\X\ldrt \mathfrak{Y}$ un morphisme entre schémas formels $\pi$-adiques sans $\pi$-torsion. Sont équivalents 
\begin{itemize}
\item $f$ est localement topologiquement de présentation finie
\item Pour tout ouvert affine  $\spf (B)\subset \YY$, $\spf (A)\subset \X$, tels que 
$f (\spf (B))\subset \spf (A)$ le morphisme induit $A\ldrt B$ fait de $B$ une $A$-algèbre topologiquement de présentation finie c'est à dire isomorphe en tant que $A$-algèbre topologique à $A<T_1,\dots,T_n>/I$ où $I$ est un idéal (saturé) de type fini de $A<T_1,\dots,T_n>$ 
\item Pour tout $k\geq 1$ le morphisme $\X\otimes \O_K/\pi^k \O_K \ldrt \mathfrak{Y}\otimes \O_K/\pi^k\O_K$ est de présentation finie
\end{itemize}
\end{lemm}

Par exemple, si $f:X\ldrt Y$ est un morphisme de présentation finie
entre schémas plats sur $\spec (\O_K)$ 
le morphisme induit entre les complétés $\pi$-adiques est topologiquement de présentation finie.

D'après les résultats de \cite{BLI} un morphisme localement
topologiquement de type fini entre schémas formels admissibles est automatiquement
localement de présentation finie. Cette assertion est en générale
fausse pour les schémas formels plus généraux que nous considérons.

\subsection{Morphismes affines}

\begin{defi}
Un morphisme de schémas formels $\pi$-adiques $\X\ldrt \mathfrak{Y}$ est dit affine si le morphisme de schémas induit  $\X\otimes \O_K/\pi\O_K \ldrt \mathfrak{Y}\otimes \O_K/\pi\O_K$ l'est ou encore de façon équivalente si $\forall k$ le morphisme de schémas induit  $\X\otimes \O_K/\pi^k \O_K \ldrt \mathfrak{Y}\otimes \O_K/\pi^k \O_K$ l'est. 
\end{defi} 

Ainsi si $\X$ est un schéma formel $\pi$-adique la catégorie des $\X$-schémas formels $\pi$-adiques affines est équivalente à la catégories des faisceaux de $\O_\X$-algèbres $\mathcal{A}$ tels que l'application canonique $\mathcal{A}\ldrt \underset{k}{\limp} \mathcal{A}/\pi^k \mathcal{A}$ est un isomorphisme et $\forall k\; \mathcal{A}/\pi^k \mathcal{A}$ est une $\O_{\X\otimes\O_K/\pi^k\O_K}$-algèbre quasi-cohérente. 
\`A $f:\mathfrak{Z} \ldrt \X$ affine on associe la $\O_\X$-algèbre $f_*\O_\mathfrak{Z}$ et dans l'autre sens à $\mathcal{A}$ on associe le $\X$-schéma formel $\underset{k}{\limi} \spec (\mathcal{A}/\pi^k \mathcal{A})$. 

La catégorie des $\X$-schémas formels $\pi$-adiques sans $\pi$-torsion est donnée par les $\O_\X$-algèbres $\mathcal{A}$ sans $\pi$-torsion telles que $\mathcal{A}\iso \underset{k}{\limp} \mathcal{A}/\pi^k \mathcal{A}$ 
et $\mathcal{A}/\pi\mathcal{A}$ soit quasi-cohérente.

\subsection{Morphismes finis}

\begin{defi}
Un morphisme $f:\X\ldrt \mathfrak{Y}$ entre schémas formels $\pi$-adiques est fini si le morphisme induit $\X\otimes \O_K/\pi\O_K \ldrt \mathfrak{Y}\otimes \O_K/\pi\O_K$ l'est.
\end{defi}

Le lemme qui suit n'est qu'une retranscription du lemme de Nakayama.

\begin{lemm}
Sont équivalents pour $f:\X\ldrt \mathfrak{Y}$
\begin{itemize}
\item $f$ est fini
\item Pour tout $k\geq 1$ le morphisme $\X\otimes \O_K/\pi^k\O_K \ldrt \mathfrak{Y}\otimes \O_K/\pi^k \O_K$ est fini
\item $f$ est affine et pour tout ouvert affine $\spf (A)\subset \mathfrak{Y}$ si $f^{-1} (\spf (A)) =\spf (B)$ alors $A$ est un $B$-module de type fini
\end{itemize}
\end{lemm}

\begin{rema}
Si $A$ est une $\O_K$-algèbre $\pi$-adique et $B$
une $A$-algèbre finie alors $B$ est complète pour la topologie
$\pi$-adique. Mais on prendra garde à ce que $B$ n'est pas
nécessairement séparée et n'est donc pas forcément une algèbre
$\pi$-adique ! 
\end{rema}

\subsection{Morphismes topologiquement plats}

\begin{defi}
Un morphisme $f:\X\ldrt \YY$ entre schémas formels $\pi$-adiques  sera dit topologiquement plat si $\forall
k\geq 1$ le morphisme de schémas induit $\X\otimes\O_K/\pi^k \O_K
\ldrt \YY\otimes \O_K/\pi^k \O_K$ l'est.
\end{defi}

Par exemple $\X$ $\pi$-adique est topologiquement plat sur $\spf
(\O_K)$ ssi il est sans $\pi$-torsion. 

\begin{rema}
On prendra garde qu'en général si $\spf (B) \ldrt \spf (A)$ est
topologiquement plat alors $B$ n'est pas nécessairement une
$A$-algèbre plate. Par exemple il n'y a pas de raison pour que si $A$
est une algèbre $\pi$-adique et $f\in A$ le morphisme $A\ldrt
A<\frac{1}{f}>$ soit plat.  
Néanmoins, d'après \cite{BLI}, tout cela est vrai pour des schémas formels
admissibles. 
\end{rema}

\begin{lemm} \label{msfpoivhnDD}
Un morphisme $f:\X\ldrt \YY$ entre schémas formels $\pi$-adiques
sans $\pi$-torsion est topologiquement plat ssi le
morphisme induit $\X\otimes \O_K/\pi\O_K\ldrt \YY\otimes\O_K/\pi\O_K$
est plat. 
\end{lemm}
\dem
C'est une conséquence du lemme 11.3.10.2 de EGA IV. 
\qed

\subsection{Limite projective dans la catégorie des schémas formels $\pi$-adiques}
\label{kchsurpjsds}

\begin{prop}
Soit $(I,\leq)$ un ensemble ordonné co-filtrant et $((\mathfrak{Z}_i)_{i\in I}, (\ph_{ij})_{i\geq j})$ un système projectif de schémas formels $\pi$-adiques tel que les morphismes de transition $\ph_{ij}:\mathfrak{Z}_i \ldrt \ZZ_j$ soient affines. Alors $\underset{i\in I}{\limp} \mathfrak{Z}_i$ existe dans la catégorie des $spf( \O_K)$-schémas formels 
et c'est un schéma formel $\pi$-adiques égal à 
$$
\underset{k\in \N}{\limi} \underset{i\in I}{\limp}
(\mathfrak{Z}_i\otimes \O_K/\pi^k \O_K)
$$
\end{prop}
\dem
La démonstration ne pose pas de problème. On renvoie au chapitre 8 de EGA IV pour les limites projectives de schémas à morphismes de transition affines. 
\qed

\begin{exem}
Si $\mathfrak{Z}_i = \spf (R_i)$ alors $\underset{i}{\limp} \spf (R_i) =\spf (\widehat{R}_\infty)$ où 
$
R_\infty = {\underset{i}{\limi} R_i }
$.
Par exemple si on considère la limite projective de revêtements de Kümmer $\DD \leftarrow \dots \leftarrow \DD \leftarrow \dots$ où $\DD =\spf (\O_K<T,T^{-1}>)$ et les morphismes de transition sont tous $t\mapsto t^p$ alors  
$$
\widehat{R}_\infty =\{ \sum_{\a\in \Z \unp} a_\a T^\a\;|\; a_\a\in \O_K\;\;
 a_\a \underset{|\a| \drt +\infty}{\ldrt} 0 \;\;\;  a_\a \underset{v_p(\a)\drt -\infty}{\ldrt} 0 \;\}
$$
où $a_\a\drt 0$ signifie tendre vers $0$ pour la topologie $\pi$-adique.
\end{exem}

\subsection{Adhérence ``schématique'' de la fibre générique}

Soit $\mathfrak{Z}$ un schéma formel $\pi$-adique. Notons pour tout
$k\geq 1$
$$
\mathcal{I}_k = \underset{i\geq 1}{\limi} \ker \left (
  \O_\mathfrak{Z}/\pi^k\O_\mathfrak{Z} \xrig{\; \times \pi^i\;}
  \O_\mathfrak{Z}/\pi^{k+i}\O_\mathfrak{Z}\right ) \subset
\O_{\mathfrak{Z}}/\pi^{k} \O_\mathfrak{Z}
$$
un faisceau d'idéaux quasicohérent sur $\mathfrak{Z}\otimes
\O_K/\pi^k$. 
Si $\spf (R)\subset  \mathfrak{Z}$ est un ouvert affine
 et $I=\{x\in R\;|\; \exists i\geq
1\; \pi^i x=0 \}$ alors 
$$
\mathcal{I}_k = (I+\pi^k R/\pi^kR)^{\widetilde{\;\;\;\;\;}}
$$
où le tilda signifie ``le faisceau quasicohérent associé''. 
Notons alors 
$$
Z_k = V(\mathcal{I}_k) \subset  \mathfrak{Z}\otimes\O_K/\pi^k\O_K
$$
On a donc
$$
Z_{k+1} \otimes \O_K/\pi^k\O_K = Z_k
$$
Notons alors
$$
\mathfrak{Z}' = \underset{k}{\limi} Z_k
$$
un schéma formel $\pi$-adique. Si $\spf (R)$ est un ouvert affine de
$\mathfrak{Z}$ et $I$ est l'idéal des éléments de $\pi^\infty$-torsion
comme précédemment alors l'ouvert correspondant de $\mathfrak{Z}'$ est
$\spf ( R/\overline{I})$ où $\overline{I}$ désigne l'adhérence de $I$
pour la topologie $\pi$-adique.

\begin{lemm}
Le schéma formel $\pi$-adique 
$\mathfrak{Z}'$ est sans $\pi$-torsion. De plus l'immersion fermée
$\mathfrak{Z}' \hookrightarrow \mathfrak{Z}$ est telle que pour tout
schéma formel $\pi$-adique sans $\pi$-torsion $\mathfrak{Y}$  tout
morphisme 
$\mathfrak{Y} \ldrt \mathfrak{Z}$ se factorise de façon unique via
$\mathfrak{Z}'\hookrightarrow \mathfrak{Z}$. 
\end{lemm}
\dem 
Elle ne pose pas de problème particulier.
\qed

\begin{defi}
Par abus de terminologie 
on appellera $\mathfrak{Z}'$ l'adhérence schématique de la fibre
générique de $\mathfrak{Z}$. On notera $\ZZ^{adh} = \ZZ'$. 
\end{defi}

Le foncteur adhérence de la fibre générique 
 définit donc un adjoint à droite à l'inclusion de la catégorie des schémas formels $\pi$-adiques sans $\pi$-torsion dans celle des schémas formels $\pi$-adiques.
\\
En particulier les produits fibrés existent dans la catégorie des schémas formels $\pi$-adiques sans $\pi$-torsion; il suffit de prendre l'adhérence schématique de la fibre générique du produit fibré usuel en tant que schémas formels.

\begin{rema}
Lorsque $\X$ est admissible d'après les résultats de \cite{BLI} (tout
idéal saturé dans une $\O_K$-algèbre topologiquement de type fini est
de type fini donc fermé)
 $\X^{adh}$ s'obtient en tuant la $\pi^\infty$-torsion dans $\O_\X$. Cela est
 faut en général, il faut quotienter par l'adhérence $\pi$-adique de
 la $\pi^\infty$-torsion.  
\end{rema}

\begin{rema}
\begin{itemize}
\item
Le schéma formel $\ZZ^{adh}$ est quasiséparé car $\ZZ^{adh}\otimes
\O_K/\pi\O_K$ est un sous-schéma fermé de $\ZZ\otimes \O_K/\pi\O_K$. 
\item Si $\mathcal{U}\subset \ZZ$ est un ouvert alors
  $\mathcal{U}^{adh} = \mathcal{U}\times_{\ZZ} \ZZ^{adh}$
\item Si $\YY\ldrt \X$ est un morphisme de schémas formels
  $\pi$-adiques alors $\YY^{adh} = (\YY\times_\X \X^{adh})^{adh}$
\end{itemize}
\end{rema}

Plus généralement que le second point de la remarque précédente on a
le lemme qui suit dont la démonstration est immédiate. 

\begin{lemm}\label{mflspivyyzrRTY3R}
Soit $\YY\ldrt \X$ un morphisme topologiquement plat entre schémas
formels $\pi$-adiques. Alors
$$
\YY^{adh} = \YY\times_\X \X^{adh}
$$
\end{lemm}

\subsection{\'Eclatements formels admissibles}
\subsubsection{Définition et premières propriétés}

\begin{defi}\label{sdghidgpp}
Soit $\mathfrak{Z}$ un schéma formel $\pi$-adique sans $\pi$-torsion
et $\mathcal{I}\subset \O_{\mathfrak{Z}}$ un idéal tel que localement
sur $\mathfrak{Z}$ $\;\exists N\in \N\; \pi^N\O_{\mathfrak{Z}}\subset
\mathcal{I}$ et $\mathcal{I}/\pi^N \O_{\mathfrak{Z}}$ est
quasi-cohérent de type fini. Un tel idéal est dit admissible.
 On appelle éclatement formel admissible
de $\mathcal{I}$ le $\mathfrak{Z}$-schéma formel $\pi$-adique topologiquement de type fini 
$$
\widetilde{\mathfrak{Z}}= \underset{k}{\limi} \text{Proj} \left (
  \bigoplus_{i\geq 0} \mathcal{I}^i/\pi^k \mathcal{I}^i \right )
$$
\end{defi}

\begin{prop}\label{kvhiuazrng}
Avec les notations de la définition précédente 
\begin{itemize}
\item si $\spf (R)\subset \mathfrak{Z}$ est un ouvert affine et 
  $I=\GG(\spf (R),\mathcal{I})$ alors $\widetilde{\mathfrak{Z}}_{|\spf
    (R)}$ s'identifie au complété $\pi$-adique de l'éclatement de
  l'idéal $\widetilde{I}$ de $\spec (R)$ 
\item $\widetilde{\mathfrak{Z}}$ est sans $\pi$-torsion
\item si $\ph : \widetilde{\mathfrak{Z}}\ldrt \mathfrak{Z}$ alors
  $\O_{\widetilde{\mathfrak{Z}}}. \ph^{-1} \mathcal{I}$ est localement
  libre de rang un 
\item $\widetilde{\mathfrak{Z}}$ satisfait à la propriété universelle suivante :
  pour tout $\mathfrak{Z}$-schéma formel $\pi$-adique sans
  $\pi$-torsion $\mathfrak{Y}\xrig{\; \psi\;} \mathfrak{Z}$ tel que
  $\O_{\mathfrak{Y}}.\psi^{-1} \mathcal{I}$ soit localement libre de
  rang un il existe un unique $\mathfrak{Z}$-morphisme
  $\mathfrak{Y}\ldrt \widetilde{\mathfrak{Z}}$
\end{itemize}
\end{prop}
\dem
La première assertion découle de la définition de l'éclatement formel
admissible et du fait que les ``Proj'' commutent aux changements de
bases (EGA II 2.8.10). 
\\
La seconde résulte de la première car $\spec (R)$ étant sans
$\pi$-torsion l'éclatement de l'idéal $\widetilde{I}$ l'est aussi
(l'image réciproque d'un ouvert schématiquement dense par un
éclatement reste schématiquement dense) et donc  son complété $\pi$-adique est encore sans
$\pi$-torsion.
\\
La troisième résulte également de la première. En effet, sur l'éclaté 
de $\widetilde{I}$ dans le schéma $\spec (R)$ l'idéal $\widetilde{I}$
devient localement libre de rang un. Il est donc localement principal sur
le complété $\pi$-adique de ce schéma. Donc $\mathcal{I}$ devient
localement principal sur $\widetilde{\mathfrak{Z}}$. Mais étant donné
que $\mathcal{I}$ contient localement une puissance de $\pi$ et que
$\widetilde{\mathfrak{Z}}$ est sans $\pi$-torsion $\mathcal{I}$
est localement libre de rang un sur $\widetilde{\mathcal{Z}}$.
\\
La dernière assertion résulte aisément de son homologue pour les
schémas (et de la première assertion).
\qed

\begin{rema}
Soit $\mathfrak{Z} = \spf (R)$ $\pi$-adique sans $\pi$-torsion 
 et $\mathcal{I}$ un idéal admissible de $\O_\mathfrak{Z}$.
Il y a alors un idéal 
 $I = (f_1,\dots,f_n)$  de $R$ contenant une puissance de $\pi$ tel
 que  
 $\mathcal{I}$ soit l'image réciproque  de l'idéal quasi-cohérent
$\widetilde{ (I/\pi^N \O_\mathfrak{Z})}$ dans $\O_{\mathfrak{Z}}$ 
pour $N>>0$. 
La description donnée dans le lemme 2.2 de \cite{BLI} de
$\widetilde{\mathfrak{Z}}$ lorsque $\mathfrak{Z}$ est topologiquement
de type fini sur $\O_K$ est en général fausse. On a en fait la
description suivante :
$$
\widetilde{\mathfrak{Z}} = \bigcup_{i=1}^n \mathcal{U}_i
$$
où $\mathcal{U}_i = \spf (A_i)$ est un ouvert affine tel que si
$$
B_i = R<T_1,\dots, \widehat{T_i},\dots, T_n>/(T_j f_i-f_j)_{1\leq j\leq
  n,j\neq i}
$$
et $B'_i =B_i /J_i$ avec $J_i =\{b\in B_i\;|\; \exists k\; \pi^k
b=0\}$ alors $A_i$ est le séparé de $B'_i$, $A_i= B'_i /\cap_{k\geq 0}
\pi^k B'_i$ ou encore $A_i = B_i/\overline{J_i}$ où $\overline{J_i}$ est l'adhérence $\pi$-adique de $J_i$.
\end{rema}

\begin{rema}\label{kvyjapfgft}
Si $\X'\ldrt \spf (A)$ est un éclatement formel admissible et si $A$
est topologiquement de type fini sur $\spf (\O_K)$ alors $\GG
(\X',\O_{\X'})\unpi = A\unpi$ (théorème d'acyclicité de Tate qui
lorsque $\O_K$ est de valuation discrète résulte du théorème de
changement de changement de base propre formel en cohomologie cohérente de EGA III). Mais
pour les schémas formels plus généraux que nous considérons il n'y a
pas de raison pour que cela soit vrai.
\end{rema}

\begin{lemm}\label{sdfjvyutat}
Soit $\mathfrak{Z}$ un schéma formel $\pi$-adique sans $\pi$-torsion et $\mathcal{I}\subset \O_\ZZ$ un idéal admissible.  
\begin{itemize}
\item $\mathcal{I}$ est localement libre de rang un ssi il est localement principal
\item si $(\pi^N)\subset \mathcal{I}$, $\mathcal{I}$ est localement principal ssi $\mathcal{I}/\pi^{N+1}\O_\ZZ$ l'est 
\end{itemize}
En particulier si $(\pi^N)\subset \mathcal{I}$, 
un morphisme $\psi: \YY\ldrt \ZZ$ se relève à l'éclatement de $\mathcal{I}$ ssi 
$(O_\YY.\psi^{-1} \mathcal{I})/\pi^{N+1} \O_\YY$ est localement principal. 
\end{lemm}
\dem
La démonstration ne pose pas de problème.
\qed

\begin{lemm} \label{dfuvhyzrp}
Soit $\mathfrak{Z}$ comme précédemment.
Soient $\mathcal{I}_1, \mathcal{I}_2$ deux idéaux admissibles.
 Soit $\widetilde{\mathfrak{Z}}$ l'éclatement
formel de $\mathcal{I}_1$. 
Alors l'éclatement formel admissible de
$\mathcal{I}_1.\mathcal{I}_2$ s'identifie à l'éclatement formel
admissible de l'image réciproque de $\mathcal{I}_2$ à
$\widetilde{\mathfrak{Z}}$.  
\end{lemm}
\dem
On peut soit le vérifier directement sur la définition de l'éclatement formel admissible, soit en utilisant la propriété universelle des éclatements formels puisque si $A$ est un anneau sans $\pi$-torsion et $I,J$ deux idéaux de type fini de $A$ contenant une puissance de $\pi$ alors $IJ$ inversible $\lssi $ $I$ et $J$ sont inversibles.
\qed

\begin{lemm}
Le composé de deux éclatements formels admissibles d'un schéma formel $\pi$-adique sans $\pi$-torsion quasicompact est un éclatement formel admissible.
\end{lemm}
\dem
On vérifie que 
la démonstration de la proposition 2.5 de \cite{BLI}, qui repose elle même sur le lemme 5.1.4 de \cite{Ray2}, fonctionne en toute généralité (on utilise l'hypothèse faite dès le début que tous nos schémas formels sont quasi-séparés).
\qed

\begin{lemm}\label{kfuspizzz}
Soit $\X$ un schéma formel $\pi$-adique sans $\pi$-torsion quasicompact et $\mathcal{U} \subset \X$ un ouvert quasicompact. Alors tout éclatement formel admissible de $\mathcal{U}$ s'étend en un éclatement formel admissible de $\X$. 
\end{lemm}
\dem La démonstration du lemme 2.6. de \cite{BLI} s'applique.
\qed

\begin{lemm}\label{sdfvyutztr}
Soit $\ph: \widetilde{\ZZ}\ldrt \ZZ$ un éclatement formel admissible. Alors 
$\ph$ est surjectif au niveau des fibres spéciales.
\end{lemm}
\dem
Par définition de l'éclatement formel admissible le morphisme induit en fibre spéciales 
$\widetilde{\ZZ}\otimes \O_K/\pi\O_K \ldrt \ZZ\otimes \O_K/\pi\O_K$
est propre. Son image est donc fermée. Soit $\mathcal{U}\subset \ZZ$
l'ouvert complémentaire de l'image. Alors $\ph^{-1} (\mathcal{U})\ldrt
\mathcal{U}$ est l'éclatement formel admissible de
$\mathcal{I}_{|\mathcal{U}}$. Il suffit donc de voir que pour un
éclatement formel admissible comme dans l'énoncé $\widetilde{\ZZ}\neq
\emptyset$. Soit donc $\spf (R)\subset \ZZ$ un ouvert affine non-vide
et $f:X\ldrt \spec (R)$ l'éclatement de l'idéal $\GG (\spf
(R),\mathcal{I})$ dans $\spec (R)$. D'après la proposition
\ref{kvhiuazrng} $\ph^{-1} (\spf (R))$ s'identifie au complété
$\pi$-adique $\widehat{X}$ de $X$. Mais si $\widehat{X}=\emptyset$ on
a $X=X\otimes_{\O_K} K$. Or le morphisme $f:X\ldrt \spec (R)$ est
propre et est un isomorphisme en fibre générique : $X\otimes K \iso \spec (R\unpi)$. Donc si l'on avait $\widehat{X}=\emptyset$ alors $f(X) = \spec (R\unpi) \subset \spec (R)$ qui serait donc ouvert/fermé dans $\spec (R)$ ce qui est impossible car $R$ est $\pi$-adique sans $\pi$-torsion.
\qed

\subsubsection{Transformée stricte} 

Soit $\ph : \mathfrak{Y}\ldrt \mathfrak{Z}$ un morphisme de schémas
formels $\pi$-adiques sans $\pi$-torsion. Soit $\mathcal{I}\subset
\O_\mathfrak{Z}$ un idéal admissible. Notons $\widetilde{\mathfrak{Z}}\ldrt
\mathfrak{Z}$ l'éclatement formel admissible associé. 

\begin{defi}
On appelle transformé strict de $\mathfrak{Y}$ relativement à
l'éclatement $\widetilde{\mathfrak{Z}}\ldrt
\mathfrak{Z}$ l'adhérence schématique de la fibre générique de 
$\mathfrak{Y}\times_\mathfrak{Z} \widetilde{\mathfrak{Z}}$.
\end{defi}

\begin{prop}\label{subbermpzq}
Le transformé strict de $\mathfrak{Y}$ s'identifie à l'éclatement
formel admissible de l'idéal $\O_\mathfrak{Y}.\ph^{-1}\mathcal{I}$. 
\end{prop}
\dem
Notons $\X$ le transformé strict et $\widetilde{\mathfrak{Y}}$
l'éclatement formel de $\O_\mathfrak{Y}.\ph^{-1}\mathcal{I}$.
Puisqu'il y a une factorisation $ \mathfrak{Y}\times_{\mathfrak{Z}}
\widetilde{\mathfrak{Z}} \ldrt \widetilde{\mathfrak{Z}}\ldrt
\mathfrak{Z}$ 
l'image réciproque à $\mathfrak{Y}\times_{\mathfrak{Z}}
\widetilde{\mathfrak{Z}}$ de l'idéal $\mathcal{I}$ est localement
principal. Donc, puisque cet idéal est admissible, il devient  localement libre de rang un sur $\X \hookrightarrow  \mathfrak{Y}\times_{\mathfrak{Z}}
\widetilde{\mathfrak{Z}}$ (lemme \ref{sdfjvyutat}). 
D'après la propriété universelle de $\widetilde{\mathfrak{Y}}$ il y a
donc un morphisme
$$
\X\ldrt \widetilde{\mathfrak{Y}}
$$
Construisons un morphisme dans l'autre sens. D'après la propriété
universelle de $\widetilde{\mathfrak{Z}}$ le morphisme
$\widetilde{\mathfrak{Y}}\ldrt \mathfrak{Z}$ s'étend en un morphisme
$\widetilde{\mathfrak{Y}}\ldrt \widetilde{\mathfrak{Z}}$. Il fournit
donc un morphisme
$$
\widetilde{\mathfrak{Y}}\ldrt  \mathfrak{Y}\times_{\mathfrak{Z}}
\widetilde{\mathfrak{Z}}
$$
Mais grâce à la propriété caractérisant l'adhérence schématique de la
fibre générique se morphisme se factorise en un morphisme
$$
\widetilde{\mathfrak{Y}}\ldrt \X
$$
On vérifie alors facilement que les deux morphismes précédents sont
inverses l'un de l'autre.
\qed

\begin{exem}
Si $\X'\ldrt \X$ est l'éclatement formel admissible de l'idéal $\mathcal{I}_1$
et $\X''\ldrt \X$ celui de $\mathcal{I}_2$ alors  $(\X'\times_\X \X'')^{adh}$ est l'éclatement formel
admissible de $\mathcal{I}_1.\mathcal{I}_2$.
\end{exem}

\subsubsection{Commutation à la limite projective}

Soit $(I,\leq)$ un ensemble ordonné cofiltrant et
$((\mathfrak{Z}_i)_{i\in I}, (\ph_{ij})_{i\geq j})$ un système
projectif de schémas formels $\pi$-adiques sans $\pi$-torsion tel que les morphismes de
transition $\ph_{ij}$ soient affines. Soit $i_0\in I$ fixé et
$\mathcal{I}$ un idéal admissible de $\O_{\mathfrak{Z}_{i_0}}$. 
Notons pour $i\geq i_0$ $\widetilde{\mathfrak{Z}}_i$ l'éclatement
formel de l'image réciproque de $\mathcal{I}$. On a donc un système
projectif $(\widetilde{\mathfrak{Z}}_i)_{i\in I}$. D'après la
proposition \ref{subbermpzq} les morphismes de transition de cette tour sont
affines.

\begin{prop}\label{vuindyuyszpgyt}
La limite projective $\underset{i\geq i_0}{\limp}
\widetilde{\mathfrak{Z}}_i$ coïncide avec l'éclatement formel de
l'image réciproque à $\underset{i}{\limp} \mathfrak{Z}_i$ de $\mathcal{I}$.
\end{prop}
\dem
Il suffit d'utiliser les propriétés universelles des éclatement et limites projectives pour construire deux morphismes inverses l'un de l'autre.
\qed

\section{La topologie des ouverts admissibles}
\subsection{La catégorie des ouverts admissibles}

\begin{defi}
Soit $\X$ un schéma formel $\pi$-adique sans $\pi$-torsion
quasicompact. Soit $\Ec_\X$ la catégorie dont les objets sont les
éclatements formels admissibles de $\X$ et les morphismes sont les
morphismes de $\X$-schémas formels. 

On note $\OEc_\X$ la catégorie fibrée au dessus de $\Ec_\X$ telle que 
\begin{itemize}
\item la fibre au dessus de $(\X'\ldrt \X)\in \text{Ob} (\Ec_\X)$  est
  la catégorie des ouverts quasicompacts de $\X'$ munie des morphismes
  donnés par l'inclusion
\item le foncteur ``changement de base'' associé à un morphisme
  $\xymatrix@R=5mm@C=5mm{\X' \ar[rr] \ar[rd] && \X'' \ar[ld] \\ & \X}$ dans
  $\Ec_\X$ est le foncteur image réciproque d'un ouvert de $\X''$ à $\X'$ 
\end{itemize} 
\end{defi}

On notera souvent $(\mathcal{U}\subset \X'\ldrt \X)$ un objet de $\OEc_\X$ où $\mathcal{U}$ est un ouvert quasicompact de l'éclatement $\X'$ de $\X$. 
\\
D'après le lemme \ref{dfuvhyzrp} la catégorie $\Ec_\X$ est
cofiltrante. Plus précisément, étant donnés deux objets de $\Ec_\X$ il
existe au plus un morphisme entre deux tels objets et si
$\xymatrix@R=5mm@C=5mm{\X' \ar[rd]  && \X'' \ar[ld] \\ & \X }$ est un
diagramme d'éclatements alors $(\X'\times_\X \X'')^{adh}$ domine ce diagramme. Ainsi la
catégorie $\Ec_\X$ est la catégorie associée à un ensemble ordonnée
dans lequel toute partie finie possède une borne supérieure.

\begin{defi}
On note $\Ad_\X$ la catégorie $\underset{\Ec_\X}{\limi} \OEc_\X$  (cf. section 6
de l'exposé VI de SGA4 tome 2  \cite{SGA4_exp6}). Il s'agit de la catégorie des ouverts admissibles quasicompacts. 
\end{defi}

Le lemme qui suit donne une définition concrète de $\Ad_\X$.

\begin{lemm}
Les objets de $\Ad_\X$ sont les ouverts quasicompacts des éclatements formels admissibles de $\X$. Si $\xymatrix@R=5mm@C=5mm{\X' \ar[rd]  && \X'' \ar[ld] \\ & \X }$ est un diagramme d'éclatements, $\mathcal{U}\subset \X'$ et $\mathcal{V}\subset \X''$ sont deux tels ouverts, si 
$(\X'\times_\X \X'')^{adh}$ désigne l'adhérence schématique de la fibre générique de $\X'\times_\X \X''$
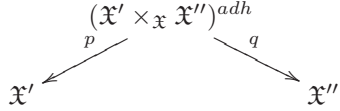
$$
\xymatrix@R=5mm@C=5mm{
& (\X'\times_\X \X'')^{adh} \ar[rd]^q \ar[ld]_{p} \\
\X' && \X''
}
$$ 
on a $\Hom_{\Ad_\X} (\mathcal{U},\mathcal{V}) \neq \emptyset$ ssi $p^{-1} ( \mathcal{U} )\subset q^{-1} ( \mathcal{V})$ et si c'est le cas alors cet ensemble d'homomorphismes est constitué d'un seul élément. 
\end{lemm}
\dem
Tout éclatement dominant $\X'$ et $\X''$ se factorise de façon unique par $(\X'\times_\X \X'')^{adh}$. 
Par définition de la catégorie limite inductive il suffit alors de voir que si 
$\X'''\xrig{\; h\;} (\X'\times_\X \X'')^{adh}$ est tel que $h^{-1} ( p^{-1} (\mathcal{U}) )\subset h^{-1} (q^{-1} (\mathcal{V}))$ alors $p^{-1} ( \mathcal{U} )\subset q^{-1} ( \mathcal{V})$. Mais cela résulte de la surjectivité des morphismes induits par les éclatements au niveau des fibres spéciales c'est à dire le lemme \ref{sdfvyutztr}.
\qed

\subsection{La topologie et le topos admissible}

\begin{defi}
Soit $\ZZ$ un schéma formel $\pi$-adique. On note $|\ZZ|_{qc}$ la
prétopologie des ouverts quasicompacts de $\ZZ$ dont les
recouvrements sont les recouvrements usuels par un nombre fini d'ouverts.
Pour $\mathcal{U}\subset \ZZ$ un ouvert quasicompact on note $\text{Cov}_{|\ZZ|_{qc}} (\mathcal{U})$ ces recouvrements.
\end{defi}

\begin{rema}
Pour vérifier que $|\ZZ|_{qc}$ vérifie bien les axiomes d'une
prétopologie il faut utiliser le fait que $\ZZ$ est quasi-séparé
(hypothèse faite sur tous nos schémas formels). 
\end{rema}

\begin{rema}
Si $|\ZZ|$ désigne la topologie usuelle il y a un foncteur pleinement fidèle continu 
$|\ZZ|_{qc} \ldrt |\ZZ|$ qui induit une équivalence de topos 
$|\ZZ|^{\widetilde{\;}} \iso |\ZZ|_{qc}^{\widetilde{\;}}$ 
puisque les
ouverts quasicompacts forment une famille génératrice de la topologie de $|\ZZ|$
(il s'agit d'une simple application du ``lemme de comparaison'', le
théorème 4.1 de SGA4 tome 1 exposé III).   
\end{rema}

Les familles 
$$
\OEc_\X \ni (\mathcal{U}\subset \X'\ldrt \X) \longmapsto \text{Cov}_{|\X'|_{qc}} (\mathcal{U})
$$
satisfont aux hypothèses du paragraphe 8.3. de \cite{SGA4_exp6} et définissent donc d'après la proposition 8.3.6 de cet exposé une prétopologie sur $\Ad_\X$ appelée prétopologie admissible. 

Concrètement on vérifie avec les définitions de l'exposé 6 de SGA4 que 

\begin{lemm}
Soit $(\mathcal{U}\subset \X'\ldrt \X )\in \text{Ob} (\Ad_\X)$ un
ouvert admissible quasicompact. Alors les recouvrements de $\mathcal{U}$ pour la prétopologie admissible sont les familles de morphismes vers $\mathcal{U}$ dans $\Ad_\X$ isomorphes aux familles finies $(\mathcal{V}_i)_{i\in I}$ d'ouverts quasicompacts d'un éclatement $\X'' \ldrt \X$ dominant $\X'\ldrt \X$, $\X''\xrig{\; h\;} \X'\ldrt \X$, telles que 
$$
h^{-1} (\mathcal{U}) = \bigcup_{i\in I} \mathcal{V}_i
$$
\end{lemm}

De façon encore plus concrète on a 

\begin{lemm}\label{kduvpuzryyt}
Soit $(\mathcal{U}\subset \X'\ldrt \X )\in \text{Ob} \Ad_\X$ un ouvert
admissible quasicompact. Les recouvrements de $\mathcal{U}$ pour la prétopologie admissibles sont les familles finies $(\mathcal{V}_i \subset \X''_i \ldrt \X)_{i\in I}$ telles que si $\mathfrak{Y}$ désigne l'adhérence de la fibre générique de $\X'\times_\X ( \underset{i\in I}{\times_\X} \X_i'')$ (la borne supérieure de tous les éclatements précédents) 
$$
\xymatrix@R=5mm@C=5mm{
 & \mathfrak{Y} \ar[ld]_p \ar[rd]^{q_i} \\
\X' & &  \X''_i
}
$$
alors $p^{-1} ( \mathcal{U}) = \bigcup_{i\in I} q_i^{-1} ( \mathcal{V}_i)$.
\end{lemm}
\dem
Il suffit d'utiliser le lemme précédent couplé au fait que les éclatements sont surjectifs en fibre spéciale (lemme \ref{sdfvyutztr}) qui implique qu'une famille d'ouverts Zariski recouvre un ouvert donné ssi c'est le cas après un éclatement.
\qed

\subsection{Le topos admissible}

Soit pour tout $(\X'\ldrt \X)\in \Ec_\X$ le topos
$|\X'|^{\widetilde{\;}} = |\X' |_{qc}^{\widetilde{\;}}$.
Si $\xymatrix@R=3mm@C=3mm{ \X'' \ar[rd] \ar[rr]^h && \X' \ar[ld] \\ & \X}$ est un morphisme dans $\Ec_\X$ il y
a alors un morphisme de topos
$$(h^*,h_*) : |\X''|^{\widetilde{\;}} \ldrt
|\X'|^{\widetilde{\;}}$$ 
 Lorsque
$\X'$ varie cela donne naissance à un topos fibré au dessus de
la petite catégorie $\Ec_\X$, au sens de l'exposé VI de SGA 4. 

\begin{defi}
On note $(\X^{rig})^{\widetilde{\;}}$ le topos limite projective 
$$
\underset{(\X'\ldrt \X)\in \Ec_\X}{\limp} |\X'|^{\widetilde{\;}} 
$$
au sens de la section 8 de l'exposé VI de SGA 4. 
On l'appelle topos rigide admissible de $\X$. 
\end{defi} 

D'après le théorème 8.2.9 de \cite{SGA4_exp6} ce topos s'identifie
à la catégorie limite projective 
$$
\underset{(\X'\ldrt \X)\in \Ec_\X}{\limp} |\X'|^{\widetilde{\;}} 
$$
où les morphismes de transition sont pour $h: \X''\ldrt \X'$
un morphisme entre éclatements le
morphisme $h_* : |\X''|^{\widetilde{\;}} \ldrt
|\X'|^{\widetilde{\;}}$. 
Concrètement un objet de ce topos est un système de faisceaux
$(\F_{\X'})_{\X'\in \Ec_\X}\in \prod_{\X'\in \Ec_\X} |\X'|^{\widetilde{\;}}$ muni d'isomorphismes $\forall h: \X''\ldrt
\X'$
$$
\a_{h} : h_* \F_{\X''}\iso \F_{\X'}
$$
vérifiant une condition de cocyle évidente.

\begin{exem}\label{ffkfdfiqdzzr}
Si $\X$ est topologiquement de type fini sur $\spf (\O_K)$ alors
d'après le théorème d'acyclicité de Tate ce topos est muni d'un
``faisceau structural'' $\O_{\X^{rig}}$ dont la ``composante'' sur l'éclatement $\X'$ 
est le faisceau $\O_{\X'}\unpi \in |\X'|^{\widetilde{\;}}$.  En général pour des $\X$ avec lesquels
nous travaillons cela est faux (cf. remarque \ref{kvyjapfgft}). Plus généralement,
toujours si $\X$ est topologiquement de type fini sur $\spf (\O_K)$, 
 si $\F$ est un faisceau cohérent sur $\X$ alors
$(h: \X'\ldrt \X) \mapsto h^* \F\unpi$ définit un faisceau cohérent de
$\O_{\X^{rig}}$-modules sur $(\X^{rig})^{\widetilde{\;}}$. 
\end{exem}

\begin{prop}
Le topos $(\X^{rig})^{\widetilde{\;}}$ s'identifie au topos des faisceaux sur $\Ad_\X$ 
muni de la topologie admissible.
\end{prop}
\dem
C'est une conséquence des résultats de la section 8.3 de \cite{SGA4_exp6}.
\qed

\subsection{Topos limite projective contre topos total}

On renvoie à la section 8.5 de \cite{SGA4_exp6} pour les généralités
concernant le lien entre une limite projective de topos fibré et son
topos total.

Traduisons les résultats de cette section 8.5 de \cite{SGA4_exp6} dans
le cas du topos rigide admissible.

\begin{defi}
Soit le topos fibré sur la catégorie $\Ec_\X$ qui à $(\X'\ldrt \X)\in
\Ec_\X$ associe $|\X'|^{\widetilde{\;}}$. On note
$\mathcal{T}_{\X^{rig}}$ le topos total associé. 
\end{defi}

Concrètement $\mathcal{T}_{\X^{rig}}$ est la catégorie des
$(\F_{(\X'\drt \X)})\in \prod_{(\X'\drt \X)\in \Ec_\X}
|\X'|^{\widetilde{\;}}$ munis de morphismes 
$$
\forall \xymatrix@R=3mm@C=3mm{\X''
\ar[rr]^h \ar[rd] && \X' \ar[ld] \\ & \X}\;\;\; \a_h: \F_{\X'\drt \X}
\ldrt h_* \F_{(\X''\drt \X)}
$$
vérifiant une certaine condition de cocyle.
Comparé au topos limite projective on relâche donc la condition pour
$\a_h$ d'être un isomorphisme. Il y a un morphisme de topos
$$
Q: (\X^{rig})^{\widetilde{\;}}\ldrt \mathcal{T}_{\X^{rig}}
$$
tel que 
$$
Q_*\left  ( (\F_{(\X'\drt \X)}), (\a_h)\right ) = \left  ( (\F_{(\X'\drt \X)}), (\a_h)\right )
$$
Voici une description du foncteur image réciproque

\begin{lemm}
Soit $\F= \left  ( (\F_{(\X'\drt \X)}), (\a_h)\right )\in
\mathcal{T}_{\X^{rig}}$. Alors $Q^*\F = \left  (
  (\mathcal{G}_{(\X'\drt \X)}), (\b_h)\right )$ où 
$$
\G_{(\X'\drt \X)} = \underset{(\X''\xrig{\; h\;} \X') }{\limi} h_*
\F_{(\X"\drt \X)}
$$
\end{lemm}
\dem 
Il s'agit d'une application de la proposition 8.5.3 de
\cite{SGA4_exp6} (pour tout $h:\X''\ldrt \X'$ $\;\; h_*
:|\X''|^{\widetilde{\;}}\ldrt |\X'|^{\widetilde{\;}}$ commute aux
petites limite projectives filtrantes (d'après le théorème 5.1 de
\cite{SGA4_exp6}, $h$ étant un morphisme cohérent). 
\qed

\begin{exem}\label{dvvhyaephyt}
Il y a un faisceaux de $\O_K$-algèbres $\mathcal{A}$ défini par $\;((\X'\ldrt \X)\longmapsto
\O_{\X'})\in \mathcal{T}_{\X^{rig}}$. Alors $Q^*\mathcal{A}$ est le
faisceau sur $(\X^{rig})^{\widetilde{\;}}$ qui à un ouvert admissible
$(\mathcal{U}\subset \X'\ldrt \X)$ associe
$$
\underset{(\mathcal{U}'\drt \mathcal{U})\in \Ec_\mathcal{U}}{\limi} \GG (\mathcal{U}', \mathcal{O}_{\mathcal{U'}})
$$
Il s'agit du faisceau associé au préfaisceau qui à un ouvert
admissible $(\mathcal{U}\subset \X'\ldrt \X)$ associe
$\GG(\mathcal{U},\O_{\mathcal{U}})$. 
Lorsque $\X$ est topologiquement de type fini sur $\spf (\O_K)$
$\O_{\X^{rig}}:=\mathcal{A}\unpi$ est le faisceau structural de
l'espace rigide (qui d'après la remarque  \ref{ffkfdfiqdzzr} se
définit sans $Q^*$, c'est déjà un faisceau sur la limite projective)
tandis que $\mathcal{A}$ est le sous-faisceau $\O_{\X^{rig}}^+$ des
fonctions rigides de normes infini $\leq 1$ (cf. \cite{Hu1} dans le contexte
des espaces adique).  Ainsi lorsque $\X$ n'est pas forcément
topologiquement de type fini sur $\spf (\O_K)$ et que le théorème
d'acyclicité de Tate n'est pas vérifié le faisceau $\mathcal{A}\unpi$
définit un faisceau structural sur $(\X^{rig})^{\widetilde{\;}}$ qui
d'après le lemme précédent s'obtient en ``forçant'' le théorème
d'acyclicité de Tate.  
\end{exem}

\subsection{Fonctorialité de la topologie et du topos admissible}

\subsubsection{Espaces rigides}

\begin{defi}
On appelle catégorie des espaces rigides quasicompacts le localisé de la catégorie des schémas formels $\pi$-adiques sans $\pi$-torsion quasicompacts relativement aux flèches données par les éclatements formels admissibles.

Si $\X$, $\YY$ sont deux schémas formels $\pi$-adiques sans $\pi$-torsion on note $\X^{rig}$ l'objet associé dans la catégorie des espaces rigides.
\end{defi}

On vérifie grâce aux propriétés des éclatements formels admissibles que cet ensemble de flèches permet un calcul des fractions à gauche.
On a donc
$$
\Hom (\X^{rig}, \YY^{rig}) =\underset{(\X'\drt \X)\in\Ec_\X}{\limi} \Hom (\X', \YY)
$$

\begin{exem}
Le foncteur $\Ec_\X\ldrt (\text{Espaces rigides}/\X^{rig})$ qui à
$(\mathcal{U}\subset \X'\ldrt \X)$ associe $\mathcal{U}^{rig}\ldrt
\X^{rig}$ est pleinement fidèle. Cela résulte du lemme
\ref{kfuspizzz}. Cela permet de voir la catégorie des ouverts
admissibles quasicompacts de $\X$ comme une sous-catégorie pleine de la catégorie $(\text{Espaces rigides}/\X^{rig})$. 
\end{exem}

\begin{rema}
Pour les espaces non-quasicompacts 
la bonne définition de la 
 catégorie des espaces rigides (quasiséparés) n'est pas le localisé
de la catégorie des schémas formels $\pi$-adiques sans $\pi$-torsion
relativement aux éclatements formels admissibles. Si $\X$ est un tel
schéma formel non-quasicompact il faut voir l'espace rigide associé à
$\X$ comme un faisceau sur le gros site admissible 
limite inductive $\underset{\mathfrak{U}\subset \X}{\limi}
\mathfrak{U}^{rig}$ où $\mathfrak{U}$ parcourt les ouvert quasicompacts
de $\X$.
\end{rema}

\subsubsection{Fonctorialité en les schémas formels}

Soit $f:\YY\ldrt \X$ un morphisme  de schémas formels $\pi$-adiques sans $\pi$-torsion quasicompacts. Les schémas formels $\X$ et $\YY$ étant quasiséparés quasicompacts un tel morphisme est quasicompact. 
 Il induit un foncteur
$$
\Ec (f) : \Ec_\X \ldrt \Ec_\YY
$$
qui à l'éclatement $(\X'\ldrt \X)$ associe le transformé strict de $\YY$ c'est à dire  $(\YY\times_\X \X')^{adh}$.

Il induit également un morphisme cartésien de catégories fibrés 
$$
\xymatrix{
\OEc_\X \ar[r]^{\OEc (f)}  \ar[d] & \OEc_\YY \ar[d] \\
\Ec_\X \ar[r]^{\Ec (f)} & \Ec_\YY
}
$$
qui en fait un morphisme de sites fibrés lorsqu'on munit les fibres de la topologie des ouverts Zariski quasi-compacts (cf. section précédente).
Il induit donc par passage à la limite inductive un morphisme de sites
$$
\Ad (f) : \Ad_\YY\ldrt \Ad_\X
$$
entre les sites admissibles de $\X$ et $\YY$.
D'où un morphisme de topos
$$
(f^*,f_*) : (\YY^{rig})^{\widetilde{\;}} \ldrt (\X^{rig})^{\widetilde{\;}}
$$
qui est le morphisme de topos induit par passage à la limite projective entre les topos fibrés au dessus de $\Ec_\YY$ et $\Ec_\X$.

\subsubsection{Fonctorialité en les espaces rigides}

\begin{prop}
Soit $f:\X'\ldrt \X$ un éclatement formel admissible. Alors $\Ad (f) : \Ad_{\X'}\ldrt \Ad_\X$, resp. $(f^*,f_*) : (\X'^{rig})^{\widetilde{\;}} \ldrt (\X^{rig})^{\widetilde{\;}}$, est
une équivalence de sites, resp. de topos.
\end{prop}
\dem
Tout a été fait pour puisque la sous-catégorie pleine $\Ec_{\X'} \subset \Ec_\X$ est co-finale.
\qed

\begin{coro}
Le site admissible et le topos admissible sont ``fonctoriels'' dans la
catégorie des espaces rigides quasicompacts :  
pour $X$ un espace rigide quasicompact on peut définit $|X|^{\widetilde{\;\;}}$ son topos
admissible, son site des ouverts admissibles quasicompacts $\Ad_X$ et un morphisme
d'espaces rigides induit un morphismes de topos et de sites.
\end{coro}

\subsection{Commutation des topos admissible à la limite projective}

Soit $(I,\leq)$ un ensemble ordonné cofiltrant et $(\X_i)_{i\in I}$ un système projectif de schémas formels $\pi$-adiques sans $\pi$-torsion quasicompacts tel que les morphismes de transition $\forall i\geq j\; \ph_{ij}: \X_i\ldrt \X_j$ soient affines. Soit
$$
\X_\infty =\underset{i\in I}{\limp} \X_i
$$
(cf. section \ref{kchsurpjsds}). 

\subsubsection{Rigidité des éclatements}

\begin{lemm}\label{kuvpahtny}
Soit $\X$ un schéma formel $\pi$-adique sans $\pi$-torsion, $\YY$ un
$\X$-schéma formel $\pi$-adique sans $\pi$-torsion et $\X'\ldrt \X$
l'éclatement formel admissible d'un idéal $\mathcal{I}$ tel que
$\pi^N\O_\X \subset \mathcal{I}$. Alors l'application de réduction
modulo $\pi^{N+1}$ 
$$
\Hom_\X (\YY, \X') \iso \Hom_{\X\otimes \O_K/\pi^{N+1}}( \YY\otimes \O_K/\pi^{N+1}
, \X'\otimes \O_K/\pi^{N+1})
$$
est une bijection.
\end{lemm}
\dem
C'est une conséquence du lemme \ref{sdfjvyutat} puisque si $f:\YY\ldrt
\X$ alors 
$\Hom_\X ( \YY,\X')$ est non-vide ssi l'idéal $\O_\YY.f^{-1}
\mathcal{I}/\pi^{N+1} \O_\YY$ est localement principal, et si c'est le
cas $\Hom_\X ( \YY,\X')$  est constitué d'un seul élément.
\qed

\subsubsection{Décomplétion des éclatements et des morphismes entre eux}

\begin{prop}\label{ksfjuzGJZR90234Tg}
\begin{enumerate}
\item
Soit $\X_\infty'\ldrt \X_\infty$ un éclatement formel admissible. Il existe alors $i_0 \in I$ et un éclatement formel admissible $\X'_{i_0} \ldrt \X_{i_0}$ tel que 
$$
\X'_\infty =\underset{i\geq i_0}{\limp} (\X'_{i_0}\times_{\X_{i_0}} \X_i )^{adh}
$$
où $ (\X'_{i_0}\times_{\X_{i_0}}\X_i )^{adh}$ désigne le transformé strict de $\X_i\ldrt \X_{i_0}$.
\item Soit 
$$
\xymatrix@R=4mm@C=4mm{
\X'_\infty \ar[rd] \ar[dd]_{f_\infty} \\
 &  \X_\infty  \\
\X''_\infty \ar[ru]
}
$$
un morphisme entre éclatements. Il existe alors $i_0\in I$ et un morphisme entre éclatements
$$
\xymatrix@R=4mm@C=4mm{
\X'_{i_0}\ar[rd] \ar[dd]_{f_{i_0}} \\
 & \X_{i_0}  \\
\X''_{i_0} \ar[ru]
}
$$
tel que $\X'_{i_0}\ldrt \X_{i_0}$, resp.  $\X''_{i_0}\ldrt \X_{i_0}$, induise l'éclatement
$\X'_\infty \ldrt \X_\infty$, resp.  $\X''_{\infty}\ldrt \X_{\infty}$, au sens du point précédent et $f_{i_0}$ induise $f_\infty$. 
\end{enumerate}
\end{prop}
\dem
Démontrons le premier point. Soit $\mathcal{I}_\infty \subset \O_{\X_\infty}$ l'idéal admissible définissant l'éclatement $\X'_\infty \ldrt \X_\infty$. Soit $N\in\N$ tel que $\pi^N \O_{\X_\infty} \subset \mathcal{I}_\infty$. \'Etant donné que l'idéal $\mathcal{I}_\infty /\pi^N \O_{\X_\infty}$ est de type fini, $\X_\infty$  quasicompact, $\X_{\infty}\otimes \O_K/\pi^N\O_K = \underset{i\in I}{\limp} \X_i\otimes \O_K/\pi^N \O_K$, il existe $i_0\in I$ et un idéal $\mathcal{I}_{i_0}$ admissible de $\O_{\X_{i_0}}$ tel que $\pi^N \O_{\X_{i_0}} \subset \mathcal{I}_{i_0}$ et 
$$
\O_{\X_\infty}.\mathcal{I}_{i_0} = \mathcal{I}_\infty
$$
Le premier point résulte alors de la proposition \ref{vuindyuyszpgyt}.

Démontrons le second point. On peut appliquer la proposition
\ref{kuvpahtny}, ou reprendre sont argument ce que nous faisons.
Soient $\mathcal{I}_\infty$, resp. $\mathcal{J}_\infty$, les idéaux admissibles définissant $\X'_\infty \ldrt \X_\infty$, resp.  $\X''_\infty \ldrt \X_\infty$. Supposons que $\pi^N \O_{\X_\infty}\subset \mathcal{I}_\infty$ et  $\pi^N \O_{\X_\infty}\subset \mathcal{J}_\infty$. Soient $i_1\in I$, $\mathcal{I}_{i_1}\subset \O_{\X_{i_1}}$, resp.  $\mathcal{J}_{i_1}\subset \O_{\X_{i_1}}$, tels que $\mathcal{I}_\infty = \O_{\X_\infty}.\mathcal{I}_{i_1}$, resp. $\mathcal{J}_\infty = \O_{\X_\infty}.\mathcal{J}_{i_1}$, comme dans le point précédent.
\\
Pour $i\geq i_1$ soit $\X'_i\ldrt \X_i$, resp. $\X''_i\ldrt \X_i$, l'éclatement de l'idéal $\O_{\X_i}.\mathcal{I}_{i_1}$, resp. $\O_{\X_i}.\mathcal{J}_{i_1}$. Il existe un morphisme $f_i$
$$
\xymatrix@R=4mm@C=4mm{
\X'_{i}\ar[rd] \ar@{-->}[dd]_{f_{i}} \\
 & \X_{i}  \\
\X''_{i} \ar[ru]
}
$$
ssi $\O_{\X'_i}.J_i = \O_{\X'_i}.J_{i_1}$ est localement libre de rang un (et si un tel morphisme existe il est unique et induit nécessairement $f_\infty$). 
Mais d'après le lemme \ref{sdfjvyutat} cela est équivalent à ce que $\O_{\X'_i}.J_{i_1}/\pi^{N+1} \O_{\X'_i}$ soit localement principal. 
Or l'existence de $f_\infty$ implique que $\underset{i\geq i_1}{\limi} \O_{\X'_i}. J_{i_1}/\pi^{N+1} \O_{\X'_i}$ est localement principal. L'idéal $J_{i_1}/\pi^{N+1} \O_{\X_{i_1}}$ étant de type fini et $\X'_\infty$ quasicompact on en déduit l'existence de $i_0\geq i_1$ tel que $\O_{\X'_{i_0}}.\mathcal{J}_{i_1}/\pi^{N+1}\O_{\X'_{i_0}}$ soit principal.
\qed

\begin{coro}\label{mdjhyvypart}
Soit $\mathcal{P}$ la catégorie fibrée au dessus de $I$ de fibre en $i\in I$ $\; \Ec_{\X_i}$ et telle que pour $i\geq j$ le foncteur ``changement de base'' $\Ec_{\X_j}\ldrt \Ec_{\X_i}$ soit l'application ``transformée strict''. Alors le foncteur naturel 
$$
\underset{I}{\limi} \mathcal{P} \ldrt \Ec_{\X_\infty}
$$
défini par la proposition \ref{vuindyuyszpgyt} induit une équivalence de catégories. 
\end{coro}

\subsubsection{La limite inductive des ouverts admissibles quasicompacts
  sont les ouverts admissibles quasicompacts 
 de la limite projective}

Reprenons les notations de la section précédente.

Lorsque $i$ varie dans $I$ les catégories des ouverts admissibles $\Ad_{\X_i}$ forment une
catégorie fibrée sur $I$. 

On s'intéresse à la catégorie
$$
\underset{i\in I}{\limi} \Ad_{\X_i} = \underset{i\in I}{\limi}
\underset{(\X'\drt \X)\in \Ec_{\X_i}}{\limi} |\X'|_{qc}
$$
Pour cela introduisons la catégorie $\mathcal{C}$ sont les objets sont
les couples $(i,\X'_i\ldrt \X_i)$ où $i\in I$ et $(\X'\ldrt
\X_i)\in\Ec_{\X_i}$ et les morphismes sont les diagramme commutatifs
pour $i\ldrt j$
$$
\xymatrix{
\X'_i \ar[r] \ar[d] & \X'_j  \ar[d] \\
\X_i \ar[r] & \X_j
}
$$
ou de façon équivalente un $\X_i$-morphisme de $\X'_i$ vers
$(\X_i\times_{\X_j}\X'_j )^{adh}$. Il résulte de cette dernière
description que cette catégorie est rigide et on vérifie aisément
qu'elle est cofiltrante. Soit la catégorie fibrée $\mathcal{D}$ au dessus de
$\mathcal{C}$ définit par la fibre au dessus de $(i,\X'_i\ldrt \X_i)$
est $|\X'_i|_{qc}$ la catégorie des ouverts quasicompacts de $\X'_i$
et les foncteurs de ``changement de base'' entre fibres sont les
foncteurs évidents. Alors 
$$
\underset{i\in I}{\limi} \Ad_{\X_i} = \underset{\mathcal{C}}{\limi} \mathcal{D}
$$
Pour tout $i\in I$ il y a un foncteur évident 
$$
\Ad_{\X_i} \ldrt \Ad_{\X_\infty}
$$
qui induit un foncteur cartésien de la catégorie fibrée des
$\Ad_{\X_i}$ lorsque $i$ varie vers la catégorie fibrée
 $I\times \Ad_{\X_\infty}$ et donc
d'après la propriété universelle des limites inductives (SGA 5 exposé
6 proposition 6.2) un foncteur
$$
\underset{i\in I}{\limi} \Ad_{\X_i} \ldrt \Ad_\X
$$
De même il y a un foncteur cartésien
$$
\mathcal{D} \ldrt \mathcal{C} \times \Ad_{\X_\infty}
$$
d'où un foncteur 
$$
\underset{C}{\limi} \mathcal{D} \ldrt \Ad_{\X_\infty}
$$
compatible à l'identification précédente de $\underset{i\in I}{\limi}
\Ad_{\X_i} = \underset{\mathcal{C}}{\limi} \mathcal{D}$.

\begin{prop}\label{vhkqpuzg}
Le foncteur $\underset{i\in I}{\limi}
\Ad_{\X_i} = \underset{C}{\limi} \mathcal{D} \ldrt \Ad_{\X_\infty}$
induit une équivalence de catégories.
\end{prop}
\dem 
La démonstration résulte facilement de la description concrète de la
catégorie limite inductive, de la proposition \ref{mdjhyvypart} et du
lemme qui suit.
\qed

\begin{lemm}\label{kicpjzr}
Soit $(Z_i)_{i\in I}$ un système projectif de schémas quasicompacts
quasiséparés 
dont les morphismes de transition sont affines et $Z_\infty =
\underset{\i\in I}{\limp} Z_i$. On a alors une équivalence de catégories $\underset{i\in I}{\limi}
|Z_i|_{qc} \iso |Z_\infty|_{qc} $. 
\end{lemm}

\subsubsection{La limite inductive des sites admissibles est le site
  admissible de la limite projective}

Lorsque $i$ varie les $\Ad_{\X_i}$ munis de leur site
admissible forment un site fibré. De plus la famille de prétopologies
sur les $\Ad_{\X_i}$ satisfont aux hypothèses de la section 8.3 
de \cite{SGA4_exp6} et définissent donc d'après la proposition 8.3.6
de \cite{SGA4_exp6} une prétopologie sur $\underset{i\in I}{\limi}
\Ad_{\X_i}$. 

Du point de vue du site fibré $\mathcal{D}$ introduit dans la section
précédente, lorsque $(i,\X'_i\ldrt \X) \in \text{Ob} (\mathcal{C})$
varie la
prétopologie $|\X'_i|_{qc}$ satisfait aux hypothèses de la section 8.3
de \cite{SGA4_exp6}
 et définit donc une prétopologie sur
$\underset{\mathcal{C}}{\limi} \mathcal{D}$. Via l'identification
entre
$\underset{i\in I}{\limi}
\Ad_{\X_i}$
et $\underset{\mathcal{C}}{\limi} \mathcal{D}$ ces deux prétopologies
coïncident.

\begin{prop}\label{udouiapdui}
L'équivalence de catégories la proposition \ref{vhkqpuzg} induit une équivalence
entre la limite inductive des sites admissibles $\Ad_{\X_i}$ et le
site admissible $\Ad_{\X_\infty}$. Plus précisément les prétopologies
définies précédemment sur ces deux catégories coïncident via cette
équivalence de catégories. 
\end{prop}
\dem
Il suffit d'utiliser la description concrète des recouvrements pour la
prétopologie donnée  par le lemme \ref{kduvpuzryyt} couplé au lemme
qui suit.
\qed

\begin{lemm}
Soit $(Z_i)_{i\in I}$ un système projectif de schémas quasicompacts
quasiséparés 
dont les morphismes de transition sont affines et $Z_\infty =
\underset{\i\in I}{\limp} Z_i$. Via l'équivalence de catégories du
lemme \ref{kicpjzr} 
la prétopologie sur  $\underset{i\in I}{\limi}
|Z_i|_{qc}$ limite inductive des prétopologies
telle que définie dans la section 8.3 de l'exposé 6 de SGA 4 coïncide
avec  celle sur $|Z_\infty|_{qc}$.
\end{lemm}
\dem
Il suffit concrètement de montrer que si 
$i_0\in I$, $(U_{\a})_\a$ est une famille finie d'ouverts
quasicompacts de $Z_{i_0}$ et $V$ un ouvert quasicompact de $Z_{i_0}$
tels que si $p:Z_\infty \ldrt X_{i_0}$ on ait 
$$
p^{-1} ( V) =\bigcup_{\a} p^{-1} ( U_\a)
$$
il existe alors $i\geq i_0$ tel que si $\ph_{i i_0} : Z_i\ldrt
Z_{i_0}$ on ait
$$
\ph_{i i_0}^{-1} ( V) =\bigcup_{\a} \ph_{i i_0}^{-1} ( U_\a)
$$
Cela ne pose aucun problème.
\qed

\subsubsection{La limite projective des topos admissibles est le topos admissible de la limite projective}

\begin{prop}\label{jfyyusouzrpgoie53TU5}
Il y a une équivalence de topos $$\underset{i\in I}{\limp} (\X_i^{rig})^{\widetilde{\;}} \simeq (\X_\infty^{rig})^{\widetilde{\;}}$$
\end{prop}
\dem
C'est une conséquence de la proposition \ref{udouiapdui}.
\qed

\section{Le point de vue spectral sur la topologie admissible}\label{REsdldidr345sf}

\subsection{Rappels sur les espaces spectraux}

On renvoie à \cite{Hoch} pour plus de détails.

\begin{defi}[\cite{Hoch}]
Un espace topologique $X$ est dit spectral s'il est quasicompact sobre et possède une base de sa topologie stable par intersections finies formée d'ouverts quasicompacts.
\end{defi}

Une définition équivalente est sobre, possède une base d'ouverts quasicompacts et une intersection finie d'ouverts quasicompacts est quasicompacte.

\begin{exem}
Soit $X$ un schéma quasicompact quasiséparé. Alors $|X|$ est spectral.
\end{exem}

On renvoie à \cite{Hoch} pour les propriétés de base des espaces spectraux. On
retiendra particulièrement les quelques faits suivants pour $X$ spectral :
\begin{itemize}
\item Par définition un ensemble constructible dans $X$ est un élément de
  l'algèbre de Boole engendrée par les ouverts quasicompacts (par
  exemple un fermé constructible est un fermé dont le complémentaire
  est quasicompact).
\item Soit $X_{cons}$ l'ensemble $X$ muni de la topologie engendrée
  par les ouverts constructibles.
Ses ouverts sont les ensembles ind-constructibles et ses fermés sont
les ensembles pro-constructibles. Alors $X_{cons}$ est compact !
\item Par exemple si $f:X\ldrt Y$ est une application continue quasicompacte entre espaces
  spectraux alors $f:X_{cons}\ldrt Y_{cons}$ est continu et donc
  l'image d'un ensemble pro-constructible est proconstructible.
\item Si $Z$ est proconstructible alors $\overline{Z} = \cup_{z\in Z}
  \overline{\{ z\} }$, les spécialisations d'éléments de $Z$, en particulier $Z$ est fermé ssi il est stable
  par spécialisation.
\end{itemize}

\subsection{Prétopologie quasicompacte sur les espaces spectraux et
  passage à la limite projective}
\label{jcutzrpz36UEPd}

Soit $X$ un espace topologique spectral. L'intersection d'un nombre fini d'ouverts quasicompacts de $X$ est quasicompacte. On note alors $X_{qc}$ la catégorie des ouverts quasicompacts de $X$ et on la munie d'une prétopologie en posant 
$$
\forall U\in X_{qc} \;\; \text{Cov} (U)=\{ \text{familles finies } (V_\a)_\a \text{ d'ouverts qc. tq. } U=\cup_\a V_\a \;\}
$$
On note encore $X_{qc}$ pour le site associé. 

Le foncteur pleinement fidèle 
$$
X_{qc} \ldrt X
$$
induit d'après le théorème 4.1 de l'exposé III de SGA4  une équivalence de topos 
$$
(X_{qc})^{\widetilde{\;}}\iso X\widetilde{\;}
$$

\begin{prop}\label{kdvyyzrp46EGodg}
Soit $(I,\leq)$ un ensemble ordonnée cofiltrant et $(X_i)_{i\in I}$ un système projectif d'espaces spectraux dont les morphismes de transition sont quasicompacts.
 Alors $X_\infty=\underset{i\in I}{\limp} X_i$ est un espace spectral. De plus il y a une équivalence de sites
$$
\underset{i\in I}{\limi} (X_i)_{qc}\iso (X_\infty)_{qc}
$$
(plus précisément les prétopologies des ouverts quasicompacts se correspondent au sens de la
section 8.3 de l'exposé 6 de SGA 4) et de topos
$$
X_\infty^{\widetilde{\;}} \iso
\underset{i\in I}{\limp} X_i^{\widetilde{\;}} 
$$
\end{prop}
\dem
Le fait que $X_\infty$ soit spectral  est démontré dans \cite{Hoch}. 

Le résultat sur les prétopologies se montre de la façon
suivante. Soit pour tout $i\in I$ $\; p_i : X_\infty \ldrt X_i$ et
$\forall i\geq j\; \ph_{ij}:X_i\ldrt X_j$. Par définition une base
d'ouvert quasicompacts de $X_\infty$ est formée des intersections
finies d'ensembles de la forme $p_i^{-1} (U_i)$ où $U_i$ est un ouvert
quasicompact de $X_\infty$. L'ensemble ordonné $(I,\leq)$ étant
cofiltrant une base stable par intersection et union finie 
est donc formée des $p_i^{-1} (U_i)$ avec $U_i$
quasicompact dans $X_i$
De cela on déduit que tout ouvert
quasicompact de $X_\infty$ provient par image réciproque d'un ouvert
quasicompact en niveau fini. Donc le foncteur
$$
\underset{i\in I}{\limi} (X_i)_{qc} \ldrt (X_{\infty})_{qc}
$$
est essentiellement surjectif.
\\
Montrons la pleine fidélité.
 Si $U,V\subset X_i$ sont deux
ouverts quasicompacts il s'agit de voir que 
$$
p_i^{-1} (U)\subset
p_j^{-1} ( V) \limpl \exists j\geq i\; \ph_{ji}^{-1} ( U)\subset
\ph_{ji}^{-1} (V)
$$
Si $p_i^{-1} (U)\subset p_i^{-1} (V)$ alors 
$$
\bigcap_{j\geq i} \ph_{ji} ( \ph_{ji}^{-1} (U)) \subset V
$$
Notons $\forall j\geq i\; K_j =  \ph_{ji} ( \ph_{ji}^{-1} (U))$. 
D'après les rappels faits sur les espaces spectraux les  $K_j$ 
sont compacts pour la topologie constructible sur $X_i$ et $V$ étant
quasicompact il est ouvert pour la topologie constructible. Donc
$\cap_{j\geq i}
K_j \subset V$ et le fait que $(I,\leq)$
soit cofiltrant implique qu'il existe $j\geq i$ tel que $K_j \subset
V$. Cela implique $\ph_{ji}^{-1} ( U )\subset \ph_{ji}^{-1}
(V)$. D'où la pleine fidélité.

Le fait que les prétopologies limite inductive et celle  de $(X_\infty
)_{qc}$ coïncident s'en déduit aisément.

Pour l'équivalence de sites on utilise une fois de plus les résultats
de la section 8.3 de \cite{SGA4_exp6} couplés à quelques manipulations élémentaires de
topologies.
L'équivalence de topos s'en déduit alors toujours d'après les résultats de
l'exposé 6 de SGA 4.
\qed

\subsection{Application au topos admissible}

\begin{defi}
Soit $\X$ un schémas formel $\pi$-adique sans $\pi$-torsion. On note
$$
|\X^{rig}| = \underset{(\X'\drt \X)\in \Ec_\X}{\limp} |\X'|
$$
comme espace topologique.
\end{defi}

Le corollaire suivant résulte de la proposition précédente.

\begin{coro}
L'espace topologique $|\X^{rig}|$ est spectral. 
Il y a une équivalence de topos $$(\X^{rig})^{\widetilde{\;}} \simeq
|\X^{rig}|^{\widetilde{\;}}$$
entre le topos admissible et la catégorie des faisceaux sur l'espace
spectral $|\X^{rig}|$. 
\end{coro}

Ainsi pour un espace rigide $\X^{rig}$ l'espace topologique
$|\X^{rig}|$ est l'espace des points du topos admissible
associé. L'espace topologique $|\X^{rig}|$ et le topos
$(\X^{rig})^{\widetilde{\;}}$ se déterminent mutuellement et la
catégorie des ouverts admissibles quasicompacts de $\X^{rig}$  
est équivalente à celle des ouverts quasicompacts de $|\X^{rig}|$. 

\begin{defi}\label{msfo7izr24}
On note $sp : |\X^{rig}|\ldrt |\X|$ l'application continue naturel
que l'on appellera application de spécialisation. D'après le lemme \ref{sdfvyutztr}
elle est surjective. 
\end{defi}

\subsection{Description de l'espace $|\X^{rig}|$ comme espace de
  Zariski-Riemann : le point de vue de Huber et Fujiwara}
\subsubsection{Rappels sur les anneaux $I$-valuatifs d'après Fujiwara}

Soit $A$ un anneau et $I$ un idéal de $A$. Dans la section 3.1 de
\cite{Fuji2} Fujiwara dit que $A$ est $I$-valuatif si 
\begin{itemize}
\item $I=(t)$ est principal où $t$ est régulier 
\item Tout idéal de type fini dans $A$ contenu une puissance de $I$
  est principal (donc inversible)
\end{itemize}

Le résultat principal de la section 3.1 de \cite{Fuji2} peut alors
s'énoncer ainsi. 

\begin{prop}[\cite{Fuji2}, section 3]
Soit $A$ local $I$-valuatif avec $I\subset \text{Rad} A$, $I=(t)$. Alors
\begin{itemize}
\item $\mathfrak{P}=\bigcap_{n\geq 0} I^n$ est un idéal premier
\item $A[\frac{1}{t}]$ est un anneau local d'idéal maximal
  $\mathfrak{P}[\frac{1}{t}]$
\item $A/\mathfrak{P}$ est un anneau de valuation de corps des
  factions $A/\mathfrak{P} [\frac{1}{t}]$, le corps résiduel de
  l'anneau local $A[\frac{1}{t}]$
\end{itemize}
\end{prop}

\subsubsection{Application au fibres du faisceau structural d'un espace
  rigide}

Soit $\X$ un schéma formel $\pi$-adique sans $\pi$-torsion.
On note $\O_{\X^{rig}}^+ \subset \O_{\X^{rig}}$ les faisceaux sur
$|\X^{rig}|$ définis dans l'exemple \ref{dvvhyaephyt}. 

\begin{coro}
Soit $x\in |\X^{rig}|$ alors le séparé $\pi$-adique de
$\O_{\X^{rig},x}^+$, $V= \O_{\X^{rig},x}^+/\cap_{n\geq 0} \pi^n \O_{\X^{rig},x}^+ $, est un anneau
de valuation sans $\pi$-torsion, $\O_{\X^{rig},x}$ est un anneau local d'idéal maximal 
$(\cap_{n\geq 0} \pi^n\O_{\X^{rig},x}^+)\unpi$ et 
de corps résiduel $V\unpi$.
\end{coro}
\dem
Il suffit de vérifier les hypothèses de la proposition précédente. Le
point $x\in |\X^{rig}|$ correspond à une collection de points
$(x_{\X'\drt \X})_{(\X'\drt \X)\in \Ec_\X}$, $x_{\X'\drt \X} \in \X'$
telle que si $h:\X''\ldrt \X'$ est un morphisme entre éclatements
alors $h(x_{\X''\drt\X})=x_{\X'\drt \X}$. Alors 
$$
\O_{\X^{rig},x}^+ = \underset{(\X'\drt \X)\in \Ec_\X}{\limi}
\O_{\X',x_{\X'\drt \X}}\;\;\;\; \O_{\X^{rig},x} = \O_{\X^{rig},x}^+\unpi
$$
Il est clair que $\pi$ est régulier dans $\O_{\X^{rig},x}^+$. Si $J$
est un idéal de type fini de $\O_{\X^{rig},x}^+$ contenant une
puissance de $\pi$ il provient par image
réciproque d'un idéal admissible  $J'$  de $\O_{\mathcal{U}}$
où $(\mathcal{U}\subset \X'\drt \X)$ est un ouvert admissible
quasicompact. L'éclatement formel admissible de l'idéal $J'$ de
$\mathcal{U}$ s'étend en un éclatement formel admissible de $\X'$ et
donc de $\X$. On en déduit donc que $J$ est principal.
\qed

\begin{defi}
Soit $x\in |\X^{rig}|$. On note $k(x)$ le corps résiduel de
$\O_{\X^{rig},x}$, un corps valué, $k(x)^0 =
\O_{\X^{rig},x}^+/\cap_{n\geq 0} \pi^n \O_{\X^{rig},x}^+$ son anneau
de valuation et $\widetilde{k(x)}$ le corps résiduel de $k(x)^0$. 
\end{defi}

\subsubsection{Anneaux de valuation rigides}

On présente ici les anneaux de valuation qui vont nous intéresser.

Soit $V$ un anneau de valuation $\pi$-adique sans $\pi$-torsion. En
particulier on suppose que la topologie de la valuation est définie
par la topologie $\pi$-adique. Cela signifie que si $v:V\ldrt
\GG\cup \{+\infty \}$ est la valuation sur $V$ alors $v(\pi)\neq \infty$
et 
$v(\pi^n)\underset{n\drt +\infty}{\ldrt} \infty$. Le corps des
fractions d'un tel anneau de valuation est $V\unpi$. On remarquera
qu'un morphisme entre deux tels anneaux de valuation est injectif.  

\begin{defi}
On appellera de tels anneaux de valuation des anneaux de valuation rigides.
\end{defi}

On supposera toujours que les valeurs de $v$ engendrent $\GG$. 
Pour un tel anneau de valuation $(V,v)$ on a une bijection entre
\begin{itemize}
\item Les idéaux premiers de $V$ non-nuls
\item Les idéaux premiers de $V/\pi V$ 
\item Les sous-groupes convexes de $\GG$
\item Les anneaux de valuation de $V\unpi$ contenant $V$
\end{itemize}
\`A un idéal premier $\mathfrak{P}\subset V$ on associe le sous-groupe
convexe $H=v(V\setminus \mathfrak{P}) \cup - v(V\setminus
\mathfrak{P})= v (V_{\mathfrak{P}}^\times)$ et l'anneau de valuation $V_\mathfrak{P}$.
Le sous-groupe $H$ étant convexe $\GG/H$ est strictement ordonné et
c'est le groupe des valeurs de la valuation de $V_\mathfrak{P}$.  
 On a le
diagramme 
$$
\xymatrix{
V \ar[r]^(.4)v \ar@{^(->}[d] & \GG\cup \{ \infty \} \ar@{->>}[d] \\
V_\mathfrak{P} \ar[r]^(.3){v'} & \GG/H\cup \{ \infty \}
}
$$
La valuation $v' : V\unpi \ldrt \GG/H\cup \{ \infty \}$ est une
générisation de la valuation $v$. Deux telles valuations définissent
la même topologie sur $V\unpi$. 

Les sous-groupes convexes de $\GG$ sont strictement ordonnés, comme le
sont les idéaux premiers de $V$ et il existe un plus petit idéal
premier $\mathfrak{P}$ non-nul dans $V$, i.e. $\mathfrak{P}$ est de
hauteur $1$. L'anneau de valuation associé $V_\mathfrak{P}$ est alors
de hauteur 1 c'est à dire que sa valuation est définie par une
valuation à valeurs dans $\R$ muni de l'ordre usuel qui est une
générisation maximale de la valuation $v$. Ainsi la topologie de
$V\unpi$ est définie par une valuation à valeurs dans $\R$. 
\\

Dans l'autre sens, étant donné $(V,v)$ comme précédemment les anneaux
de valuation de $V\unpi$ contenus dans $V$ correspondent aux anneaux
de valuation du corps résiduel de $V$, $V/\mathfrak{m}_V$. \`A
$V'\subset V$ on associe l'image de $V'$ dans
$V/\mathfrak{m}_V$. Ainsi avec les notations précédentes, le corps
résiduel de $V_\mathfrak{P}$ est $\text{Frac} (V/\mathfrak{P})$ et 
$V/\mathfrak{P} V \subset \text{Frac} (V/\mathfrak{P})$ correspond au
sous-anneau de valuation $V$ de $V_\mathfrak{P}$.

\subsubsection{Points rigides}

Soit $V$ un anneau de valuation rigide.
Tout idéal de type fini dans $V$ est principal, donc d'après la
propriété universelle des éclatements 
tout point $x:\spf (V)\ldrt \X$ s'étend de manière unique en un système compatible de points 
$$
\xymatrix@R=3mm@C=3mm{
 & \X' \ar[d] \\
\spf (V) \ar[ru] \ar[r] & \X
}
$$
pour $(\X'\ldrt\X)\in \Ec_\X$. L'image du point fermé de $\spf (V)$
définit dont un élément de $|\X^{rig}|$. 

Réciproquement étant donné $x\in \X^{rig}$ d'après la section
précédent fournit un morphisme $\spf ( \widehat{k(x)^0})\ldrt \X$ qui
redonne le point $x$ grâce à la construction précédente. 

On obtient ainsi la description

\begin{prop}\label{mspivhsbbryy}
Il y a une bijection entre $|\X^{rig}|$ et les classes de points
 $\spf (V) \ldrt \X$ où $V$ est un anneau de valuation rigide et deux
 points $\spf (V_1)\ldrt \X, \; \spf (V_2) \ldrt \X$ sont équivalents
 si il existe un diagramme 
$$
\xymatrix@R=4mm@C=4mm{
\spf (V_1) \ar[rd] \\
 & \spf (V_3) \ar[r] & \X  \\
\spf (V_2) \ar[ru] 
}
$$
où les morphismes entres spectres formels d'anneaux de valuation
envoient le point fermé sur le point fermé (i.e. il s'agit de la
relation de domination des anneaux de valuation).
\end{prop}

En particulier

\begin{coro}
Soit $\X=\spf (A)$. Alors $|\X^{rig}|$ est l'ensemble des classes de valuations
$v:A\ldrt \GG\cup \{ \infty \}$ telles que 
\begin{itemize}
\item $v(A)$ engendre le groupe des valeurs de la valuation
\item $v(\pi)\neq \infty$
\item $v(\pi^n)\underset{n\drt +\infty}{\ldrt} \infty$ i.e. $\forall
  \gamma\in \GG\;\exists n\in \N\; n v(\pi)\geq \gamma$
\end{itemize}
et deux valuations $(v_1,\GG_1)$ et $(v_2,\GG_2)$ sont équivalentes si
il existe un isomorphisme de groupes ordonnés $\a:\GG_1\iso \GG_2$ tel
que $\a\circ v_1= v_2$.  
\end{coro}

\begin{defi}
Pour $A$ comme précédemment on notera $\text{Spa} (A)$ l'ensemble de valuations précédentes. 
On munit $\Spa (A)$ de la topologie induite par la topologie $|\spf (A)^{rig}|$. 
\end{defi}

Avec les notations de Huber, lorsque $A$ est topologiquement de type fini, c'est ce que Huber note $\text{Spa} ( A^\rhd, A^+)$ où $A^\rhd = A\unpi$ et $A^+$ est la fermeture intégrale de $A$ dans $A\unpi$. 

Enfin notons le lemme suivant qui donne une définition plus pratique
de $|\X^{rig}|$.

\begin{lemm}
Il y a une bijection entre $|\X^{rig}|$ est les classes d'équivalences
de points $\spf (V)\ldrt \X$ où $V$ est un anneau de valuation rigide
et deux points $\spf (V_1)\ldrt \X,\; \spf (V_2)\ldrt \X$ sont
équivalents ssi il existe un diagramme
$$
\xymatrix@R=3mm@C=4mm{
 & \spf (V_1) \ar[rd] \\
\spf (V') \ar[ru] \ar[rd] && \X \\
 & \spf (V_2) \ar[ru] 
}
$$
où $V'$ est un anneau de valuation rigide 
tel que les morphismes $\spf (V')\ldrt \spf (V_1)$, $\spf (V')\ldrt
\spf (V')\ldrt \spf (V_2)$ envoient le point fermé sur le point fermé.
\end{lemm}
\dem
La condition est évidemment suffisante. Réciproquement, avec les
notations de la proposition \ref{mspivhsbbryy} il suffit de voir qu'il
existe un élément de $|(\spf (V_1)\times_{\spf (V_3)} \spf
(V_2))^{rig}|$ s'envoyant sur les points fermés dans $|\spec (V_1/\pi
V_1)|\times_{|\spec (V_3/\pi V_3)|} |\spec (V_2)|$. Mais cela résulte
de la surjectivité de $|\spec (V_1/\pi V_1)\times_{\spec (V_2/\pi
  V_2)} \spec (V_3/\pi V_3)|\ldrt |\spec (V_1/\pi
V_1)|\times_{|\spec (V_3/\pi V_3)|} |\spec (V_2)|$ et de la
surjectivité du morphisme de spécialisation (lemme \ref{msfo7izr24}) 
$sp : |(\spf (V_1)\times_{\spf (V_3)} \spf
(V_2))^{rig}| \ldrt |\spf (V_1)\times_{\spf (V_3)}\spf (V_2)|$.
\qed
 
\subsubsection{Topologie de l'espace de Zariski-Riemann}

\begin{defi}
Soit $A$ un anneau $\pi$-adique sans $\pi$-torsion. Soit $g,f_1,\dots,f_n\in A$ une famille finie  telle que l'idéal engendré dans $A$ contienne une puissance de $\pi$. Soit $X=\text{Spa} (A)$. On note 
$$
X\left <\frac{f_1,\dots, f_n}{g}\right > =\{ v\in \text{Spa} (A)\;|\; \forall i\; v(f_i)\geq v (g)\}
$$
\end{defi}

Ainsi si l'idéal  $(f_1,\dots, f_n)$ de $A$ contient une puissance de $\pi$ on a un recouvrement 
$$
\Spa (A) = \bigcup_{1\leq i\leq n} \Spa (A) \left < \frac{f_1,\dots,\hat{f}_i,\dots, f_n}{f_i} \right >
$$

\begin{prop}
Soit $X =\Spa (A)$. Alors lorsque $(f_0,f_1,\dots,f_n)$ varie parmi les familles finies engendrant un idéal contenant une puissance de $\pi$ les ensembles $X\left <\frac{f_1,\dots, f_n}{f_0}\right >$ forment une base d'ouverts quasicompacts stable par intersections finies de le topologie de $\Spa (A)$. 
\end{prop}
\dem
Soit $\X=\spf (A)$.
Une base d'ouverts quasicompacts de $|\X^{rig}|$ est donnée par les 
$|\mathcal{U}^{rig}|\subset |\X^{rig}|$ où $(\mathcal{U}\subset \X'\ldrt \X)$ est un ouvert admissible de $\X$. 

Soit $I$ un idéal de type fini dans $A$ contenant une puissance de $\pi$. Soit $X=\spec (A)$ et $X'=\text{Proj} (\oplus_{k\geq 0} I^k )\ldrt X$ l'éclatement de l'idéal $I$ dans $X$. L'éclatement formel admissible associé est $\X'= (X')^{\widehat{\;}}$, le complété $\pi$-adique de $X'$. Une base d'ouverts quasicompacts de $X'$ est donnée par 
$$ \forall n\geq 0\;
\forall f\in I^n\;\; \{ x\in X'\; |\; f\text{ engendre l'idéal inversible } \O_{X'}.\widetilde{I}^n\text{ en } x \;\}
$$
lorsque $n$ et $f$ varient. Les ouverts de $\X'$ sont les traces d'ouverts de $X'$ sur la fibre spéciale $ X'\otimes \O_K/\pi \O_K = \X'\otimes \O_K/\pi\O_K$. Pour $V$ un anneau de valuation rigide il y a une bijection entre les points de $\X^{rig}$, $\spf (V)\ldrt \X'$, et les morphismes $\spec (V)\ldrt X'$. Si $U\subset X'$ est un ouvert et $\mathcal{U}\subset \X'$ est l'ouvert associé de la fibre spéciale un morphisme $\spf (V) \ldrt \X'$ se factorise par 
$\mathcal{U}$ ssi le morphisme associé $\spec (V)\ldrt X'$ se factorise par $U$. 

Maintenant si $n\geq 0$, $I^n=(f_1,\dots,f_n)$ et $g\in I^n$ alors un morphisme 
$\a:A\ldrt V$ et donc $\spec (V)\ldrt X'\ldrt X=\spec (A)$ 
se factorise par l'ouvert où $g$ engendre $\O_{X'}.\widetilde{I}^n$ ssi $(\a (g))$ engendre $(\a(f_1),\dots, \a (f_n))$ dans $V$ c'est à dire ssi $\forall i\; v(\a ( f_i))\geq v(\a (g))$.
\qed
\\

Lorsque $A$ est topologiquement de type fini sur $\spf (\O_K)$ on
retrouve donc bien la topologie de l'espace valuatif au sens de Huber
(\cite{Hu7} et \cite{Hu6}). 

Voici le lemme qui permet de recoller la description précédente dans le cas affine. 

\begin{lemm}
Soit $\X$ un schéma formel $\pi$-adique sans $\pi$-torsion
quasicompact et $\X=\bigcup_{i\in I} \mathcal{U}_i$ un recouvrement
 fini de $\mathcal{X}$ par des ouverts quaiscompacts. Alors l'espace
 topologique $|\X^{rig}|$ est obtenu par recollement des
 $|\mathcal{U}_i^{rig}|$, $i\in I$, le long des ouverts $|(\mathcal{U}_i\cap\mathcal{U}_j)^{rig}|$, $i,j\in I$.  
\end{lemm}
\dem 
C'est une conséquence du lemme \ref{kfuspizzz}.
\qed

\subsection{Ouverts surconvergents et espace analytique de Berkovich}\label{kdvhuzr46HrhfG}

\subsubsection{Générisations dans l'espace de Zariski-Riemann}

On vérifie avec les notations de la section \ref{jcutzrpz36UEPd} que
si $X_\infty =\underset{i\in I}{\limp} X_i$ est une limite projective
d'espaces spectraux à morphismes de transition quasicompacts alors 
$$
\forall x=(x_i)_{i\in I}, y=(y_i)_{i\in I}\in X_\infty\;\; x\succ y
\lssi \forall i\; x_i \succ y_i
$$
Soit $\X$ un schéma formel $\pi$-adique sans $\pi$-torsion et
$(x_{\X'})_{\X'\drt \X}, (y_{\X'})_{\X'\drt \X}\in |\X^{rig}|$. Alors
$x\succ y \lssi \forall (\X'\ldrt \X)\in \Ec_\X\; x_{\X'}\succ
y_{\X'}$ dans $|\X'|$ (et il suffit de le vérifier dans un ensemble
cofinal d'éclatements puisque les morphismes de transition entre
ceux-ci sont propres). 

\begin{lemm}
Soit $\spf (V)\xrig{\; x\;} \X$ un point rigide. Les générisations de
$[x]\in |\X^{rig}|$ sont les points $[x_\mathfrak{P}]$ où 
$x_\mathfrak{P}:\spf (\widehat{V_{\mathfrak{P}}})\ldrt \spf (V)\ldrt
\X$ (complétion $\pi$-adique de $V_\mathfrak{P}$) avec $\mathfrak{P}$ un idéal premier non-nul dans $V$.
\end{lemm}
\dem
Si $y\succ [x]$ où $y\in |\X^{rig}|$ alors d'après les considérations
précédentes $\forall k\geq 1\; \O_{\X^{rig},y}^+/(\pi^k)$ est obtenu
comme une localisation de $\O_{\X^{rig},[x]}^+/(\pi^k)$ puisque c'est
le cas en niveau fini pour chaque éclatement de $\X$.
\qed

\begin{coro}
Si $x\in |\X^{rig}|$ l'ensemble des générisations de $x$ est
totalement ordonné et possède un unique élément maximal associé à une
valuation de rang $1$. Si $\X=\spf (A)$ et $x$ correspond à la
valuation $v:A\ldrt \GG\cup \{+\infty \}$ alors l'ensemble des
générisations de $x$ est en bijection avec l'ensemble totalement
ordonné des sous-groupes convexes propres de $\GG$. \`A $H\subset \GG$
convexe est associé $v/H : A\xrig{\; v\;} \GG\cup \{+\infty \}\ldrt
\GG/H\cup \{+\infty \}$.
\end{coro}

Dans le langage de Huber ``toutes les générisations dans $\text{Spa}
(A)$ sont des générisations secondaires'' (\cite{Hu7} et \cite{Hu6}). 

\subsubsection{L'espace analytique de Berkovich associé}

\begin{defi}
Soit $\X$ $\pi$-adique sans $\pi$-torsion quasicompact.
On note $|\X^{an}|$ le plus grand quotient séparé de $|\X^{rig}|$ que
l'on muni de la topologie quotient.
\end{defi}

L'ensemble $|\X^{an}|$ s'identifie à l'ensemble des générisations
maximales de $|\X^{rig}|$. Ses ouverts correspondent aux ouverts de
$|\X^{rig}|$ stables par spécisalisation (les ouverts
surconvergents). Les points de $|\X^{an}|$ sont donc donnés par les
points rigides $\spf (\O_L)\ldrt \X$ où $L|K$ est une extension valuée
compléte pour une valuation $v:L\ldrt \R\cup\{+\infty \}$. 

Comme dans les lemmes 8.1.5 et 8.1.8. de \cite{Hu1} on peut vérifier que
l'ensemble des ouverts surconvergents de $|\X^{rig}|$ sont les 
intérieurs des ensembles fermés constructibles de $|\X^{rig}|$. De
plus lorsque $\X=\spf (A)$ la topologie sur $|X^{an}|$ est celle engendrée 
par les ouverts du type
$$
\{ x\in |\X^{an}|\;|\; |f(x)|<|g(x)|\} \text{ où } f,g\in A
$$
Donc, $|\X^{an}|$ s'identifie à l'espace topologie de Berkovich
$\mathcal{M} (A\unpi)$ associé à l'algèbre de Banach $A\unpi$
(cf. \cite{Berkspectral}).

\section{\'Etude des morphismes finis localement libres rig-étales entre schémas formels}

\subsection{Morphismes finis localement libres}

Soit $f:\X\ldrt \mathfrak{Y}$ un morphisme entre schémas formels $\pi$-adiques sans $\pi$-torsion. 

\begin{defi}
Le morphisme $f$ est fini localement libre si le morphisme induit entre schémas $\X\otimes \O_K/\pi\O_K\ldrt \mathfrak{Y}\otimes\O_K/\pi\O_K$ l'est. 
\end{defi}

\begin{lemm}
Sont équivalents
\begin{itemize}
\item $f$ est fini localement libre
\item $f$ est fini topologiquement de présentation finie et
  topologiquement plat
\item Pour tout entier $k\geq 1$ le morphisme induit $\X\otimes \O_K/\pi^k \O_K\ldrt \mathfrak{Y} \otimes\O_K/\pi^k \O_K$ est fini localement libre 
\end{itemize}
\end{lemm}
\dem
Appliquer le lemme \ref{msfpoivhnDD}. 
\qed

\begin{lemm}
Soit $A$ un anneau $\pi$-adique sans $\pi$-torsion, $M$ un $A$-module de type fini $\pi$-adiquement séparé (i.e. $M$ est $\pi$-adique et $M/\pi M$ est de type fini) tel que $\forall k\; M/\pi^k M$ soit un $A/\pi^k A$-module projectif. Alors $M$ est un $A$-module projectif.  En particulier $M$ est de présentation finie.
\end{lemm}
\dem
Soit $u:A^n\twoheadrightarrow M$ une surjection. On cherche une section $\e:M\ldrt A^n$ telle que $u\circ \e =Id$. Pour tout entier $k\geq 1$ notons $X_k$ l'ensemble des section de $u_k= u\text{ mod } \pi^k$, $u_k : (A/\pi^k A)^n \twoheadrightarrow M/\pi^k M$,
$$
X_k= \{ \e_k\in \Hom_{A/\pi^k A} (M/\pi^k M, (A/\pi^k A)^n)\;|\; u_k\circ \e_k= Id \}
$$
Par hypothèse $\forall k\; X_k\neq \emptyset$. De plus $X_k$ est un $\Hom_{A/\pi^k A} ( M/\pi^k M, \ker u_k)$-torseur. 

\'Etant donné que dans la suite exacte
$$
0\ldrt \ker u_{k+1} \ldrt (A/\pi^{k+1}A)^n \xrig{\; u_{k+1}\;} M/\pi^{k+1} M \ldrt 0
$$
les modules sont projectifs on obtient par application de $-\otimes_{A/\pi^{k+1} A} A/\pi^k A$ l'égalité  $\ker u_k = \ker u_{k+1}/\pi^k \ker u_{k+1}$. De cela on déduit que $\forall k\geq 1$ l'application
$$
\Hom_{A/\pi^{k+1}A} (M/\pi^{k+1} M, \ker u_{k+1}) \ldrt \Hom_{A/\pi^k A} ( M/\pi^k M, \ker u_k)
$$
est surjective. En effet, si $v\in \Hom_{A/\pi^k A} ( M/\pi^k M, \ker u_k)$ l'existence de $v'$ dans le diagramme qui suit
$$
\xymatrix@R=7mm@C=7mm{
M/\pi^{k+1} M \ar@{..>}[r]^{v'} \ar@{->>}[d] & \ker u_{k+1} \ar@{->>}[d] \\
M/\pi^k M \ar[r]^v & \ker u_k
}
$$
résulte de la projectivité de $M/\pi^{k+1} M$ comme $A/\pi^{k+1}A$-module.

Donc $\underset{k}{\limp} X_k\neq \emptyset$.
\qed

On montre également aisément le lemme qui suit.

\begin{lemm}
Le lemme précédent reste valable en remplaçant module projectif par
module libre.
\end{lemm}

\begin{coro}
Les assertions suivantes sont équivalentes pour $f:\X\ldrt
\mathfrak{Y}$ comme précédemment.
\begin{itemize}
\item  Le morphisme $f$ est fini localement libre
\item $f$ est
affine et pour tout ouvert $\spf (A)\subset \mathfrak{Y}$ si $f^{-1} (
\spf (A))=\spf (B)$ alors $B$ est un $A$-module projectif de type fini
\item $f$ est affine et il existe un recouvrement affine
  $(\mathcal{U}_i)_{i}$ de $\mathfrak{Y}$, $\mathfrak{U}_i =\spf (
  A_i)$ tel que si $f^{-1} ( \mathcal{U}_i) = \spf (B_i)$ alors $B_i$
  est un $A_i$-module libre
\end{itemize}
\end{coro}

\begin{coro}
Soit $f:\X\ldrt \mathfrak{Y}$ un morphisme fini localement libre entre
schémas formels $\pi$-adiques sans $\pi$-torsion et un morphisme 
$\mathfrak{Y}'\ldrt \mathfrak{Y}$ où $\mathfrak{Y}'$ est également
$\pi$-adiques sans $\pi$-torsion. Alors le changement de base $\X\times_{\mathfrak{Y}} \mathfrak{Y}'$ 
est sans $\pi$-torsion et est fini localement libre au dessus de $\YY'$.
 Par exemple si $\mathfrak{Y}'\ldrt
\mathfrak{Y}$ est un éclatement formel admissible le transformé strict
de $\X\ldrt \mathfrak{Y}$ est $\X\times_{\mathfrak{Y}} \mathfrak{Y}'$.
\end{coro}
\dem
Il suffit de remarquer que si $A$ est $\pi$-adique sans $\pi$-torsion
et $B$ une $A$-algèbre projective de type fini alors pour tout
changement de base $A\ldrt A'$ avec $A'$ $\pi$-adique sans
$\pi$-torsion
$$
B\otimes_A A' = B\widehat{\otimes}_A A'
$$
En effet, $B$ étant un $A$-module de type fini $B\otimes_A A'$ est
$\pi$-adiquement complet, et étant projectif il est séparé (il se
plonge dans un $(A')^r$ qui est séparé). 
\qed

\subsection{Morphismes finis localement libres rig-étales}

\subsubsection{Algèbres  finies de présentation finie rigidifiées :
rappels}\label{sfkviyzpvba}

On rappelle ici une notion introduite dans le chapitre III.2 de
\cite{Elkik} auquel on renvoie pour plus de détails.

Soit $A$ un anneau et $M$ un $A$-module de présentation finie muni
d'une présentation
$$
A^p\xrig{\; L\;} A^q \ldrt M\ldrt 0
$$
où $L$ est une matrice $q\times p$. Soit le foncteur sur la catégorie
des $A$-algèbres qui à la $A$-algèbre $B$ associe l'ensemble des
structures de $B$-algèbres rigidifiées sur $M\otimes_A B$ relativement à la présentation
$L$  au sens de \cite{Elkik} (section
III.2). Rappelons que toute
structure d'algèbre de sur $M\otimes_A B$ possède une telle rigidification.
Ce foncteur est représentable par un $A$-schémas affine de
présentation finie $W\ldrt \spec (A)$. Soit $V\subset W$ l'ouvert où
l'algèbre universelle est étale au dessus de $W$.
On vérifie également comme dans \cite{Elkik} que $V$ est lisse au
dessus de $\spec (A)$. 

\subsubsection{Définition et propriété des morphismes finis localement
  libres rig-étales}

\begin{lemm}
Soit $A$ un anneaux $\pi$-adique sans $\pi$-torsion et $X\ldrt \spec
(A)$ un $A$-schéma affine. Soit $U\subset X$ un ouvert et $\e$ une
section
$$
\xymatrix{
X\ar[d] \\
\spec (A) \ar@(ru,lu)[u]_(.26)\e
}
$$
Supposons qu'il existe un recouvrement formel $\spf (A) =
\bigcup_{i\in I} \spf (A <\frac{1}{f_i} >)$ tel que $\forall i\;$ la
section $\e$ sur $\spec (A<\frac{1}{f_i}>)$ se factorise par $U$ en
dehors de $V(\pi)$
$$
\xymatrix{
 & &  U\ar@{^(->}[r]  & X\ar[d] \\
\spec (A<\frac{1}{f_i}>\unpi) \ar@{^(->}[r] \ar@(ur,lu)[rru]
 & \spec (A<\frac{1}{f_i}>)
\ar[r] & \spec (A[\frac{1}{f_i}]) \ar@{^(->}[r] & \spec (A)\ar@(ru,lu)[u]_(.26)\e
}
$$
Alors $\e$ se factorise à travers $U$ en dehors de $V(\pi)$.
\end{lemm}
\dem
Soit $X=\spec (D)$, $I\subset D$ un idéal tel que $U= V(I)^c$ et
$J\subset D$ l'idéal définissant la section $\e$, $J=\ker ( \e^* :D
\twoheadrightarrow A)$. 

L'assertion $\e$ se factorise par $U$ hors de $V(\pi)$ est équivalente
à 
$$
\exists k\in \N\;\; \pi^k \in I+J
$$
L'anneau étant $\pi$-adique elle est encore équivalente à ce que 
$$
\exists k\in \N\;\; \pi^k \in I+J+\pi^{k+1} D/\pi^{k+1} D \subset
D/\pi^{k+1} D
$$
Soit $\forall k\in \N$  $\;\Delta_k$ l'idéal $I+J+\pi^{k+1}
D/\pi^{k+1} D$ et $\widetilde{\Delta}_k$ le faisceau d'idéaux sur
$\spec (D/\pi^{k+1} D)$. 

Par hypothèse 
$$
\forall i\in I\; \exists k\in \N\;\; (\pi^k) \subset
\widetilde{\Delta}_{k| D(f_i)} \text{ où } D(f_i)\subset \spec
(D/\pi^{k+1} D)
$$
Donc, $\exists k\;\forall i \;\; (\pi^k) \subset
\widetilde{\Delta}_{k|D (f_i)}$ ce qui implique puisque $\spec
(D/\pi^{k+1} D) = \bigcup_i D(f_i)$ que $(\pi^k)\subset \Delta_k$.
\qed

\begin{prop}
Soit $f: \X\ldrt \mathfrak{Y}$ un morphisme fini localement libre
entre schémas formels $\pi$-adiques sans $\pi$-torsion. Sont
équivalents
\begin{itemize}
\item Il existe un recouvrement affine $\mathcal{Y} = \bigcup_i
  \mathcal{U}_i$ tel que si $\mathcal{U}_i= \spf (A_i)$, $ f^{-1} (
  \mathcal{U}_i)=\spf (B_i)$, alors $B_i\unpi$ est une $A_i
  \unpi$-algèbre étale.
\item Pour tout ouvert affine $\mathcal{U} \subset \mathfrak{Y}$,
  $\mathcal{U} =\spf (A)$, si $f^{-1} ( \mathcal{U}) = \spf (B)$ alors
  $B\unpi$ est une $A\unpi$-algèbre étale.  
\end{itemize}
\end{prop}
\dem
Si $A$ est un anneau $\pi$-adique et $B$ une $A$-algèbre finie
projective alors $B$ est une $A$-algèbre $\pi$-adique. En effet, étant
donné que $B$ est un $A$-module de type fini $B$ est $\pi$-adiquement
complète, et étant donné que $B$ est projective $\exists n\; B\subset
A^n$ et est donc séparée puisque $A^n$ l'est.
\\
Donc pour une telle $A$-algèbre $B$ et pour $f\in A\; B\otimes_A
A<\frac{1}{f}>$ vérifiant les hypothèses précédentes en remplaçant $A$
par $A<\frac{1}{f}>$ et $B$ par $B\otimes_A A<\frac{1}{f}>$ on obtient
$$
B\otimes_A A<\frac{1}{f}> 
 = B<\frac{1}{f}>
$$ 
et est une $A<\frac{1}{f}>$-algèbre
finie projective.

Venons en maintenant à la démonstration. 
Soit $A$ une algèbre $\pi$-adique sans $\pi$-torsion et $B$ une
$A$-algèbre finie projective étale en dehors de $V(\pi)$. Alors
d'après ce qui précède $\forall f\in A\;\; B<\frac{1}{f}>$ est une
$A<\frac{1}{f}>$-algèbre 
 étale
en dehors de $V(\pi)$. 

Il s'agit maintenant de voir que si $A$ est $\pi$-adique sans
$\pi$-torsion et que $B$ est une $A$-algèbre finie projective telle
qu'il existe un recouvrement formel $spf (A) =\bigcup_i \spf
(A<\frac{1}{f_i}>)$ tel que $\forall i\; B<\frac{1}{f_i}>$ soit une
$A<\frac{1}{f_i}>$-algèbre étale en dehors de $V(\pi)$ alors $B$ est
étale hors de $V(\pi)$. Mais, après avoir choisi une présentation
finie de $B$ comme $A$-module puis une rigidification de la structure
de $A$-algèbre sur $B$ relativement à cette présentation, 
 cela résulte du lemme précédent appliqué 
au schéma $X=W$ de la section \ref{sfkviyzpvba} et son ouvert $V$.
\qed

\begin{defi}
Un morphisme satisfaisant aux conditions de la proposition précédente
sera dit fini localement libre rig-étale. 
\end{defi}

\begin{lemm}
Soit $f:\X\ldrt \mathfrak{Y}$ un morphisme fini localement libre
rig-étale de schémas formels $\pi$-adiques sans $\pi$-torsion et
$\mathfrak{Y}'\ldrt \mathfrak{Y}$ un morphisme où $\mathfrak{Y}'$ est
$\pi$-adique sans $\pi$-torsion. Alors $\X\times_\mathfrak{Y}
\mathfrak{Y}' \ldrt \mathfrak{Y}'$ est fini localement libre
rig-étale.

En particulier le transformé strict d'un morphisme fini localement
libre rig-étale en est un. 
\end{lemm}
\dem
Cela résulte de ce que si $A$ est $\pi$-adique sans $\pi$-torsion, $B$
une $A$-algèbre finie projective étale en dehors de $V(\pi)$, $C$ une
$A$-algèbre $\pi$-adique sans $\pi$-torsion alors, comme expliqué dans
la démonstration précédente, $B$ étant projective finie sur $A$
$$
B\otimes_A C = B\widehat{\otimes}_A C
$$
\qed

\begin{lemm}
Le composé de deux morphismes finis localement libres rig-étales en est encore un.
\end{lemm}

\subsection{Rigidité}

\begin{prop}[Fujiwara]\label{lsdspvuuhntyhg}
Soit $\X$ un schéma formel $\pi$-adique sans $\pi$-torsion quasicompact et 
$\mathfrak{Y}$ un $\X$-schéma formel $\pi$-adique sans $\pi$-torsion fini localement libre rig-étale. Il existe alors un entier $N (\mathfrak{Y})$ tel que pour tout $\X$-schéma formel $\pi$-adique sans $\pi$-torsion $\mathfrak{Y}'$ et deux morphismes 
$$
\xymatrix{
\mathfrak{Y}' \ar[rd] \ar@<.6ex>[rr]^{f_1} \ar@<-.6ex>[rr]_{f_2} && \mathfrak{Y} \ar[ld] \\ 
 & \X
}
$$
vérifiant $f_1 \equiv f_2 \text{ mod } \pi^{N(\mathfrak{Y})}$ on ait $f_1=f_2$.

De plus si $\mathfrak{Z}\ldrt \mathfrak{X}$ est un morphisme de changement de base avec $\mathfrak{Z}$ $\pi$-adique sans $\pi$-torsion alors on peut prendre $N(\mathfrak{Y}\times_\X \mathfrak{Z}) = N(\mathfrak{Y})$. 
\end{prop}
\dem
Lorsque tous les schémas formels sont affines 
c'est une conséquence de la proposition 2.1.1 de \cite{Fuji2}.
Plus précisément 
lorsque $\mathfrak{Y} = \spf (B)$ et $\X=\spf (A)$ on peut prendre $N(\mathfrak{Y})$ tel que 
$\pi^{N(\mathfrak{Y})} \Omega^1_{B/A} =0$. Le où les schémas
formels ne sont plus affines s'en déduit facilement.
L'assertion  sur le changement de base se déduit aussitôt de la borne
explicite donnée dans le cas affine. 
\qed

\begin{theo}\label{klfkjdiuappfjt}
Soient $\X,\mathfrak{Y}$ comme dans la proposition précédente. Il existe alors un entier $N (\mathfrak{Y})$ tel que pour tout $\mathfrak{Y}'$ comme précédemment 
l'application de réduction 
$$
\Hom_{\X} ( \mathfrak{Y'} ,\mathfrak{Y}) \ldrt \Hom_{\X \otimes\O_K/\pi^{N(\mathfrak{Y})} } ( \mathfrak{Y'}\otimes\O_K/\pi^{N(\mathfrak{Y})} , \mathfrak{Y}\otimes\O_K/\pi^{N(\mathfrak{Y})} )
$$
soit une bijection. De plus, comme précédemment, on peut choisir l'entier $N(\mathfrak{Y})$  invariant pas changement de base.
\end{theo}
\dem
Supposons d'abord que tous nos schémas formels soient affines
$$
\X=\spf (A),\; \mathfrak{Y} =\spf ( B),\; \YY'=\spf (B')
$$
La $A$-algèbre $B$ est de présentation finie. \`A chaque choix d'un
ensemble fini de générateurs et relations $\mathcal{S}$ de cette
$A$-algèbre i.e. un isomorphisme
$$
B\simeq A[X_1,\dots,X_N]/J\text{ et } J=(f_1,\dots,f_q)
$$
Elkik associe dans \cite{Elkik} un idéal $H_{B/A} ( \mathcal{S})$ de
$A[X_1,\dots,X_N]$ tel que l'ouvert 
$$
V ( H_{B/A} ( \mathcal{S}))^c \cap V(J) \subset \spec (B)
$$
soit l'ouvert de lissité de $\spec (B)\ldrt \spec (A)$. De plus si $\ph
: A\ldrt A'$ est une $A$-algèbre alors $\mathcal{S}$ détermine une
présentation $\ph (\mathcal{S})$ de $B\otimes_A A'$ 
$$
B\otimes_A A ' \simeq A'[X_1,\dots X_N]/ (f'_1,\dots, f'_q) \text{ où
} \forall i \; f'_i \text{ est l'image de } f_i \text{ dans } A'[\underline{X}]
$$
et on a 
$$
A'[X_1,\dots,X_N]. H_{B/A} ( \mathcal{S}) \subset H_{B\otimes_A A'} (
\ph (\mathcal{S})
$$
Fixons un tel système de générateurs/relations de $B/A$. Par
hypothèse, $\spec (B)\ldrt \spec (A)$ étant lisse en dehors de
$V(\pi)$,
$$
\exists N_0\in \N\;\; \pi^{N_0}\in H_{B/A} ( \mathcal{S})+J
$$
D'après le lemme 1 de \cite{Elkik} si $N_1>2 N_0+1$ et $\overline{\e}$ 
$$
\xymatrix{
\spec ( B/\pi^{N_1} B) \ar[d] \\
\spec (A/\pi^{N_1} A) \ar@(ru,lu)[u]_(.26){\overline{\e}}
}
$$
est une section approchée modulo $\pi^{N_1}$ de $\spec (B)\ldrt \spec
(A)$ alors $\exists \e$ une section
$$
\xymatrix{
\spec (  B) \ar[d] \\
\spec ( A) \ar@(ru,lu)[u]_(.26){\e}
}
$$
telle que $\e \equiv \overline{\e}\; [ \pi^{N_1 -N_0}]$.

Appliquons cela à 
$$
\xymatrix{
\spec (B)\times_{\spec (A)} \spec (B') \ar[d] \\
\spec (B')
}
$$
D'après les rappels précédents sur le fait que l'idéal $H_{B/A} (
\mathcal{S})$ ne peut que grandir par changement de base on en
déduit que si $N_1>2N_0$  et $\overline{f}$
$$
\xymatrix{
\spec ( B'/\pi^{N_1 B'}) \ar[rr]^{\overline{f}}\ar[rd] & & \spec ( B/\pi^{N_1}
B) \ar[ld] \\
&\spec (A/\pi^{N_1} A)
}
$$
est un morphisme approché mod $\pi^{N_1}$ alors $\exists f$
$$
\xymatrix{
\spec (B') \ar[rr]^{f} \ar[rd] && \spec (B) \ar[ld] \\
& \spec (A)
}
$$
tel que $f\equiv \overline{f} \; [\pi^{N_1 -N_0}]$. 

Choisissons $N_1$ tel que $N_1>2N_0$ et $N_1-N_0$ soit strictement
supérieur à l'entier de la proposition précédente. On obtient alors
que 
$$
\Hom_{\spec (A)} ( \spec (B'),\spec (B))\iso \Hom_{\spec ( A/\pi^{N_1}
  A)} ( \spec ( B'/\pi^{N_1} B'), \spec ( B/\pi^{N_1} B))
$$
Le cas où nos schémas formels ne sont pas affines s'en déduit car
grâce à l'assertion d'injectivité (la proposition précédente) on voit
que les morphismes construits entre ouverts affines se recollent
automatiquement. 

L'assertion concernant le fait que l'on peut choisir l'entier
$N(\mathfrak{Y})$ invariant pas changement de base résulte de ce que
l'idéal $H_{B/A} ( \mathcal{S})$ grandit par changement de base et de
l'invariance par changement de base de l'entier de 
la proposition précédente.
\qed

\subsection{Décompletion des schémas formels finis localement libres
  rig-étales}

\subsubsection{Décompletion des algèbres}

\begin{lemm}
Soit $(A,\mathcal{I})$ un couple hensélien et $M$ un $\widehat{A}$-module projectif de type fini. Il existe alors un $A$-module projectif de type fini $N$ et un isomorphisme $M\simeq N\otimes_A \widehat{A}$. 
\end{lemm}
\dem
Soit $n\in \N$ et $\mathcal{M}_{n/A}=\spec (A[X_{ij}]_{1\leq i,j\leq n})$ le schéma en $A$-algèbres des matrices carrées de taille $n$.  
Le sous-schémas de $\mathcal{M}_{n/A}$ classifiant les idempotents $\e\in \mathcal{M}_{n/A}$, $\e^2 =\e$, est étale. En effet, si $B$ est un anneau, $J$ un idéal de carré nul dans $B$, $\e\in M_n (B)$ un élément vérifiant $\e^2 = \e + x$ avec $x\in M_n (J)$ alors $\e x =x\e$ et donc si $\e' = \e + (1- 2\e) x$ on a 
$$
\e'^2 =\e' \text{ et } \e' \equiv \e \text{ mod } J
$$
Si $\e \in M_n (\widehat{A})$ est un idempotent il existe donc $\e'\in M_n (A)$ tel que 
$$
\e'² = \e' \text{ et } \e' \equiv \e \text{ mod } \mathcal{I}
$$
Maintenant si $B$ est une algèbre munie d'un idéal de carré nul $J$, $\e$ et $\e'$ sont deux idempotents de $M_n (B)$ et $\e' = \e + x$ où $x\in M_n (J)$ alors 
$$
\e' = (Id + x) \e (Id + x)
$$
Donc si $\e,\e'\in M_n (\widehat{A})$ sont deux idempotents tels que $\e'\equiv \e \text{ mod
 } \mathcal{I}$ alors il existe $u,v\in \GL_n (\widehat{A})$ tels que 
$$
\e' = u\e v
$$
D'où le résultat. 
\qed

\begin{rema}
Lorsque $A$ est noethérien le lemme précédent est bien sûr beaucoup plus faible que le théorème 3 de \cite{Elkik}. Le théorème 3 de \cite{Elkik} s'étend cependant aux cas où on ne fait pas d'hypothèse de noethérianité sur $A$, $\mathcal{I} =(t)$ avec $t$ régulier. C'est ce cas que nous utiliserons mais dans notre situation nous supposerons, avec les notations du théorème 3 de \cite{Elkik}, que le module $\overline{M}$ est localement libre sur $\widehat{A}$ et pas seulement en dehors de $V(\mathcal{I})$. 
 Dans ce cas là, comme on vient de le voir, la preuve est nettement plus simple. C'est pourquoi nous l'avons donnée. 
\end{rema}

\begin{prop}\label{mfpsufnuytrzpp}
Soit $(I,\leq)$ un ensemble ordonné filtrant et $(R_i)_{i\in I}$ un système inductif d'algèbres $\pi$-adiques sans $\pi$-torsion. Notons $R_\infty =\underset{i\in I}{\limi} R_i$ et $\widehat{R}_\infty$  le complété $\pi$-adique de $R_\infty$. Soit $B_\infty$ une $\widehat{R}_\infty$ algèbre $\pi$-adique sans $\pi$-torsion telle que $\spf ( B_\infty )\ldrt \spf ( \widehat{R}_\infty )$ soit fini localement libre rig-étale (i.e. $B_\infty$ est un $\widehat{R}_\infty$-module projectif de type fini et $B_\infty \unpi$ une $\widehat{R}_\infty\unpi$-algèbre étale). Il existe alors $i_0\in I$ et une $R_{i_0}$-algèbre $\pi$-adique sans $\pi$-torsion $B_{i_0}$ telle que $\spf (B_{i_0}) \ldrt \spf (R_{i_0})$ soit fini localement libre rig-étale et un isomorphisme
$$
B_\infty \simeq B_{i_0}\otimes_{R_{i_0}} \widehat{R}_\infty  = B_{i_0}\widehat{\otimes}_{R_{i_0}} \widehat{R}_\infty
$$
\end{prop}
\dem
Le couple $(R_\infty,\pi R_\infty)$ est hensélien.
D'après le lemme précédent il existe un $R_\infty$-module projectif de type fini $M$ 
et un isomorphisme de $\widehat{R}_\infty$-modules 
$$
B_\infty \simeq M\otimes_{R_\infty} \widehat{R}_\infty
$$
Fixons un tel isomorphisme.
Fixons une présentation finie $L$ de $M$. Soit $W$ le schéma affine au dessus de $\spec (R_\infty)$ classifiant les structures d'algèbres rigidifiées sur $M$ relativement à la présentation $L$ (cf. section \ref{sfkviyzpvba}). Soit $V\subset W$ l'ouvert où l'algèbre universelle est étale. 

La présentation $L$ fournit grâce à l'isomorphisme précédent $B_\infty \simeq M\otimes_{R_\infty} \widehat{R}_\infty$ 
 une présentation de $B_\infty$ comme $\widehat{R}_\infty$-module. Fixons une structure d'algèbre rigidifiée relativement à cette présentation de $B_\infty$. Cela nous fournit donc une section
$\e$ de $W$ au dessus de $\widehat{R}_\infty$ qui par hypothèse se factorise par $V$ en dehors de $V(\pi)$ 
$$
\xymatrix{
 & V\ar@{^(->}[r] & W\ar[d] \\
\spec (\widehat{R}_\infty \unpi) \ar@{^(->}[r] \ar[ru] & \spec (\widehat{R}_\infty )\ar[r] \ar[ru]^\e & \spec (R_\infty)
}
$$
On peut maintenant appliquer le théorème 2 bis de \cite{Elkik} pour
conclure l'existence pour tout entier $n$ d'une section $\e'$ de $W$  au dessus de $\spec (R_\infty)$ se
factorisant par l'ouvert $V$ en dehors de $V(\pi)$ telle que 
$$
\e'\equiv \e \; [\pi^n]
$$
et donc le diagramme suivant privé de ses termes extrêmes gauche et droite commute mod $\pi^n$ 
$$
\xymatrix{
 & & V\ar@{^(->}[d] \\
 & &  W\ar[d] \\
\spec (\widehat{R}_\infty \unpi) \ar@{^(->}[r] \ar@(ru,l)[rruu] & \spec
(\widehat{R}_\infty )\ar[r] \ar[ru]^\e & \spec (R_\infty)
\ar@(ru,lu)[u]_(.26){\e'} & \ar@{_(->}[l] \spec (R_\infty \unpi) \ar@(ru,r)[luu]
}
$$
Pour un entier $n$ donné il existe donc une $R_\infty$-algèbre
$B'_\infty$ telle que $B'_\infty \simeq M$ comme $R_\infty$-module,
$B'_\infty\unpi/R_\infty \unpi$ est étale et un isomorphisme
$$
B'_\infty/\pi^n B'_\infty \simeq B_\infty /\pi^n B_\infty
$$
L'algèbre $B'_\infty$ est donc un $R_\infty$-module projectif fde type
fini. Donc $B'_\infty \otimes_{R_\infty} \widehat{R}_\infty=\widehat{B}'_\infty$ et $\spf
(\widehat{B}'_\infty)\ldrt \spf (\widehat{R}_\infty)$ est fini
localement libre rig-étale. D'après le théorème \ref{klfkjdiuappfjt}
pour $n$ choisi suffisamment grand (ne dépendant que de
$B_\infty/\widehat{R}_\infty$) il existe un morphisme relevant
l'identité
$$
\xymatrix{
\spf ( \widehat{B}'_\infty )\ar[rr]^f \ar[rd] && \ar[ld] \spf (B_\infty) 
 \\
& \spf (\widehat{R}_\infty)
}
$$
tel que $f\equiv Id \text{ mod }\pi^n$ via l'isomorphisme
$B'_\infty/\pi^n B'_\infty \simeq B_\infty /\pi^n B_\infty$. D'après
le lemme qui suit c'est un isomorphisme.

Il est maintenant aisé de vérifier qu'il existe $i_0\in I$ et
$B_{i_0}$ une $R_{i_0}$-algèbre projective finie étale en dehors de
$\pi$ telle que $B'_\infty \simeq B_{i_0}\otimes_{R_{i_0}} R_\infty$.
\qed

\begin{lemm}
Soit $A$ un anneau $\pi$-adique et $P_1$ $\; P_2$ deux $A$-module
projectifs de type fini. Un élément $f\in \Hom_A (P_1,P_2)$ est un
isomorphisme ssi $f\text{ mod }\pi : P_1/\pi P_1\ldrt P_2/\pi P_2$ en
est un.
\end{lemm}
\dem
Les modules étant projectifs de type fini ils sont $\pi$-adique. Il
suffit donc de montrer que $\forall n\; f\text{ mod } \pi^n$ est un
isomorphisme. D'après le lemme de Nakayama $\forall n\; f \text{ mod
}\pi^n$ est surjectif. De plus $\forall n$ les fonctions localement
constantes rang des modules projectifs $P_1/\pi^n P_1$,
resp. $P_2/\pi^n P_2$, sont égales sur $\spec (A/\pi^n A)= \spec
(A/\pi A)$ puisque $f\text{ mod }\pi$ est un isomorphisme. Donc
$f\text{ mod }\pi^n$ est un isomorphisme.
\qed

\begin{theo}\label{skfivyzptuu}
Soit $(I,\leq)$ un ensemble ordonné cofiltrant et $(\X_i)_{i\in I}$ un
système projectif de schémas formels $\pi$-adiques sans $\pi$-torsion
quasicompacts sont les morphismes de transition sont affines. Posons
$\X_\infty =\underset{i\in I}{\limp} \X_i$. Pour tout schéma formel
$\pi$-adique sans $\pi$-torsion $\mathfrak{Z}$ notons
$\mathcal{R}_\ZZ$ le catégorie des $\ZZ$-schémas formels
($\pi$-adiques sans $\pi$-torsion) finis localement libres
rig-étales. Le système $I\ni i\longmapsto \mathcal{R}_{\X_i}$ forme
une catégorie fibrée via les applications de changement de base. De
plus lorsque $i$ varie il y a un système compatible de foncteurs
``changement de base'' $\mathcal{R}_{\X_i}\ldrt
\mathcal{R}_{\X_\infty}$ d'où un foncteur
$$
\underset{i\in I}{\limi} \mathcal{R}_{\X_i} \ldrt \mathcal{R}_{\X_\infty}
$$
Ce foncteur induit une équivalence de catégories.
\end{theo}
\dem
Commençons par la pleine fidélité de ce foncteur. Soient $i_0\in I$ et
$\YY',\YY$ deux $\X_{i_0}$-schémas formels dans
$\mathcal{R}_{\X_{i_0}}$. Soit $N =N(\YY)$ l'entier associé à $\YY$
fourni par le théorème \ref{klfkjdiuappfjt} (qui rappelons le peut
être choisi invariant par changement de base). Alors
\begin{eqnarray*}
&& \underset{i\geq i_0}{\limi} \Hom_{\X_i} ( \YY'\times_{\X_{i_0}} \X_i,
\YY\times_{\X_{i_0}} \X_i ) \\
&=& \underset{\i\in I}{\limi}
\Hom_{\X_i\otimes \O_K/\pi^N} ( \YY'\otimes \O_K/\pi^N
\times_{\X_{i_0}\otimes \O_K/\pi^N } \X_i\otimes \O_K/\pi^N,  \YY\otimes \O_K/\pi^N
\times_{\X_{i_0}\otimes \O_K/\pi^N } \X_i\otimes \O_K/\pi^N  ) \\
&=& 
\Hom_{\X_\infty \otimes \O_K/\pi^N } (  \YY'\otimes \O_K/\pi^N
\times_{\X_{i_0}\otimes \O_K/\pi^N } \X_\infty \otimes \O_K/\pi^N,  \YY\otimes \O_K/\pi^N
\times_{\X_{i_0}\otimes \O_K/\pi^N } \X_\infty\otimes \O_K/\pi^N  ) \\
&=& \Hom_{\X_\infty} ( \YY'\times_{\X_{i_0}} \X_\infty,
\YY\times_{\X_{i_0}} \X_\infty )
 \end{eqnarray*}
L'avant dernière égalité résultant de ce que modulo $\pi^N$ $\;\YY$ et
$\YY'$ sont de présentation finie sur $\X_{i_0}\otimes\O_K/\pi^N\O_K$,
et la dernière par une nouvelle application du théorème
\ref{klfkjdiuappfjt}.

Reste à vérifier la surjectivité essentielle de notre foncteur. La
proposition \ref{mfpsufnuytrzpp} l'affirme lorsque les schémas formels
$\X_i$ sont affines. Lorsqu'ils ne le sont pas il suffit de choisir
$i_0\in I$ et une décomposition $\X_{i_0} =\bigcup_\a \mathcal{U}_\a$
en un nombre fini d'ouverts affines. Alors si $\YY\ldrt \X_\infty$ est
un élément de $\mathcal{R}_{\X_\infty}$ d'après la proposition
\ref{mfpsufnuytrzpp}, étant donné qu'il n'y a qu'un nombre fini
d'ouverts affines, il existe $i_1\geq i_0$ et des éléments
$(\mathcal{V}_\a)_\a$ de $\mathcal{R}_{\mathcal{U}_\a
  \times_{\X_{i_0}} \X_{i_1}}$ et des isomorphismes  
$$
\forall \a\; \mathcal{V}_\a \times_{\X_{i_1}}  \X_{\infty}
\iso \YY\times_{\X_\infty} (\mathcal{U}_\a \times_{\X_{i_0}} \X_\infty )
$$
Grâce à l'égalité entre les $\Hom$ démontrée précédemment (la pleine
fidélité), quitte à augmenter $i_1$ en $i_2\geq i_1$ on peut supposer que les
$\mathcal{V}_\a$ sont munis d'une donnée de recollement au dessus de
$\X_{i_2}$ relativement
aux $(\mathcal{U}_\a \cap \mathcal{U}_\b)\times \X_{i_2}$. 
Toujours grâce à l'égalité entre les $\Hom$ et quitte à encore
augmenter l'indice $i_2$ en $i_3\geq i_2$ on peut supposer que ces données de
recollement satisfont à la condition de cocyle permettant de les
recoller en un schéma formel au dessus de $\X_{i_3}$.
\qed

\section{\'Etude d'une certaine catégorie de morphismes rig-étales}
\subsection{Définitions}

\begin{defi}
Un morphisme $f:\YY \ldrt \X$ entre schémas formels $\pi$-adiques sera
dit étale si le morphisme induit entre schémas $\YY\otimes \O_K/\pi\O_K
\ldrt \X\otimes \O_K/\pi \O_K$ l'est.
\end{defi}

Bien sûr comme d'habitude si $\X$ est un schéma formel $\pi$-adiques
l'application de réduction modulo $\pi$ induit une équivalence entre
les $\X$-schémas formels $\pi$-adiques étales et les $\X\otimes \O_K/\pi\O_K$-schémas étales.

\begin{defi}\label{mspfihztyzr}
Un morphisme $f:\YY \ldrt \X$ entre schémas formels $\pi$-adiques sans
$\pi$-torsion quasicompacts sera dit de type $(\mathcal{E})$ s'il se factorise en un
composé de morphismes de schémas formels $\pi$-adiques sans
$\pi$-torsion
$$
\YY \ldrt \mathfrak{S} \ldrt \mathfrak{T} \ldrt \mathfrak{W} \ldrt \X
$$
où
\begin{itemize}
\item $\mathfrak{W}\ldrt \X$ est isomorphe à un éclatement formel admissible de
  $\X$
\item $\mathfrak{T}\ldrt \mathfrak{W}$ est étale quasicompact
\item $\mathfrak{S}\ldrt \mathfrak{T}$ est fini localement libre rig-étale
\item $\YY\ldrt \mathfrak{S}$ est isomorphe à un éclatement formel
  admissible de $\mathfrak{S}$
\end{itemize}
\end{defi}

\begin{exem}
Par exemple un ouvert admissible quasicompact de $\X$, $(\mathcal{U}\subset
\X'\ldrt \X)$ est un morphisme de type $(\mathcal{E})$.
\end{exem}

\begin{lemm}\label{jfuzptzszr}
Soit $\YY\ldrt \X$ un morphisme de type $(\E)$ et $\ZZ\ldrt \X$ un
morphisme avec $\ZZ$ $\pi$-adique sans $\pi$-torsion quasicompact alors le
morphisme $(\YY\times_\X \ZZ)^{adh} \ldrt \ZZ$ est de type $(\E)$. 
\end{lemm}
\dem 
Utilisant le lemme \ref{mflspivyyzrRTY3R} on vérifie que si 
$\YY \ldrt \mathfrak{S} \ldrt \mathfrak{T} \ldrt \mathfrak{W} \ldrt
\X$ est une décomposition comme dans la définition \ref{mspfihztyzr}
alors
$$
(\YY\times_\X \ZZ)^{adh} \ldrt (\mathfrak{S}\times_\X \ZZ)^{adh} \ldrt
(\mathfrak{T}\times_\X \ZZ)^{adh} \ldrt (\mathfrak{W}\times_\X \ZZ)^{adh} \ldrt
\X
$$
en est encore une. 
\qed

\subsection{Les morphismes de type $(\E)$ engendrent la topologie
  étale des espaces rigides usuels}

\begin{theo}\label{rhfueo5pjgf}
Soit $\YY\ldrt \X$ un morphisme rig-étale entre deux schémas formels
admissibles quasicompacts. Il existe alors
un diagramme
$$
\xymatrix{
\YY' \ar[rr]^f \ar[rd]^g && \YY \ar[ld]\\
 & \X
}
$$
où $g$ est de type $(\E)$ et $|f^{rig}|: |\YY'^{rig}|\ldrt |\YY^{rig}|$ est
surjectif au niveau des espaces de Zariski-Riemann tels que définis dans la
section \ref{REsdldidr345sf}.
\end{theo}
\dem
D'après le théorème de platification de Raynaud-Gruson (\cite{BLII}) il
existe un éclatement formel admissible $\X'\ldrt \X$ tel que le
transformé strict
$$
\ZZ =( \X'\times_\X\YY)^{adh} \ldrt \X'
$$
soit plat (au sens de \cite{BLI}). Ce morphisme étant rig-étale plat il est
quasi-fini en fibre spéciale : le morphisme de schémas 
$$
\ZZ\otimes\O_K/\pi\O_K \ldrt \X'\otimes\O_K/\pi\O_K
$$
est quasi-fini. Il existe donc un diagramme de $\O_K/\pi\O_K$-schémas 
$$
\xymatrix{
W\ar[d]_\a \ar[r]^(.3)\b & \ZZ\otimes\O_K/\pi\O_K \ar[d] \\
T \ar[r]_(.3)\gamma & \X'\otimes \O_K/\pi\O_K
}
$$
où $\a$ est plat fini, $\gamma$ est étale et $\beta$ est étale
surjectif.
Les morphismes $\b$ et $\gamma$ étant étales ce diagramme se relève de
façon unique en un diagramme de schémas formels admissibles
$$
\xymatrix{
\mathfrak{W} \ar[r]^{\widetilde{\b}} \ar[d]_{\widetilde{\a}} & \ZZ
\ar[d]^h \\
\mathfrak{T} \ar[r]^{\widetilde{\gamma}} & \X'
}
$$
où $\widetilde{\a}$ est fini localement libre. \'Etant donné que
$\widetilde{\gamma}^{rig}\circ \widetilde{\a}^{rig}$ est rig-étale et
$\widetilde{\gamma}^{rig}$ est rig-étale, $\widetilde{\a}$ est
rig-étale. De plus $\widetilde{\b}$ étant étale surjectif
$\widetilde{\b}^{rig}$ est surjectif au niveau des espaces
spectraux. D'où le résultat.
\qed

\subsection{Rigidité}

\begin{prop} \label{DGTeztJTUtuu2258}
Soit $\YY\ldrt \X$ un morphisme de type $(\E)$ entre schémas formels
$\pi$-adiques sans $\pi$-torsion. Il existe alors un entier $N$ tel
que pour tout $\X$-schéma formel $\pi$-adique sans $\pi$-torsion
l'application de réduction modulo $\pi^N$
$$
\Hom_\X ( \YY',\YY)\iso \Hom_{\X\otimes \O_K/\pi^N} ( \YY'\otimes
\O_K/\pi^N, \YY\otimes \O_K/\pi^N)
$$
soit une bijection.
\end{prop}
\dem
Il suffit d'empiler les différentes assertions de rigidité concernant
les morphismes qui entrent dans la définition d'un morphisme de type
$(\E)$ : le lemme \ref{kuvpahtny} pour les éclatements, la proposition 
\ref{lsdspvuuhntyhg} pour le morphisme fini localement libre
rig-étale, quant aux morphisme étale c'est évident.
\qed

\subsection{Décompletion}\label{kfusyttgtnbgg}

\begin{lemm}\label{vvmpzfinggz}
Soit $(I,\leq)$ un ensemble ordonné cofiltrant, $(\X)_{i\in I}$ un
système projectif de schémas formels $\pi$-adiques
et $\X_\infty =\underset{i\in I}{\limp }\X_i$. 
 Si $\ZZ_{\et}$ 
désigne la catégorie des  morphismes étales quasicompacts vers $\ZZ$
alors 
$$
\underset{i\in I}{\limi} (\X_i)_{\et} \iso (\X_\infty)_{\et}
$$
est une équivalence
\end{lemm}
\dem
\'Etant donné que $\ZZ_\et$ est équivalente à $(\ZZ\otimes
\O_K/\pi\O_K)_{\et}$ cela se ramène à un énoncé classique sur les
schémas (cf. EGA IV).
\qed

\begin{theo}\label{DGufbrzgzr3277Jne}
Soit $(I,\leq)$ un ensemble ordonné cofiltrant et $(\X_i)_{i\in I}$ un
système projectif de schémas formels $\pi$-adiques sans $\pi$-torsion quasicompacts
dont les morphismes de transition sont affines. Soit $\X_\infty =
\underset{i\in I}{\limp} \X_i$. 
 Pour tout $i\in I$
soit $\E_{\X_i}$ la catégorie des morphismes de type $(\E)$ au dessus
de $\X_i$. Lorsque $i$ varie les applications de changement de base
(cf. lemme \ref{jfuzptzszr}) définissent une catégorie fibrée. Il y a
un foncteur 
$$
\underset{i\in I}{\limi} \E_{\X_i} \ldrt \E_{\X_\infty}
$$
Ce foncteur induit une équivalence de catégories.
\end{theo}
\dem
C'est un exercice d'empilage du 
 corollaire \ref{mdjhyvypart} du théorème
\ref{skfivyzptuu} et du lemme \ref{vvmpzfinggz}.
\qed
\\

Soit $\X$ un schéma formel $\pi$-adique sans $\pi$-torsion. Les
morphismes de type $(\E)$ $\YY\ldrt \X$ sont stables par éclatement
formel admissible : si $\YY'\ldrt \YY$ est un éclatement formel
admissible le composé $\YY'\ldrt \YY\ldrt \X$ est encore de type
$(\E)$. On peut donc définir $\E_{\X^{rig}}$ la catégorie des
morphismes de type $(\E)$ vers $\X$ localisée relativement aux
éclatements. On démontre de la même façon en utilisant une fois
de plus le corollaire \ref{mdjhyvypart} et le théorème précédent.

\begin{theo}\label{mdpojkzrnjznhfbeyty24}
Avec les notations du théorème précédent on a une équivalence de
catégories
$$
\underset{i\in I}{\limi} \E_{\X_i^{rig}} \ldrt \E_{\X_\infty^{rig}}
$$
\end{theo}

\section{Le topos rig-étale d'un schéma formel $\pi$-adique quasicompact}
\subsection{Sur un point concernant les topologies de Grothendieck}

Nous utiliserons le théorème clef suivant. 

\begin{theo}[SGA4 exposé III, théorème 4.1]\label{odpmmlfjhjutrr}
Soit $\mathcal{C}$ une petite catégorie, $\mathcal{C}'$ un site et
$u:\mathcal{C}\ldrt \mathcal{C}'$ un foncteur pleinement
fidèle. Supposons que tout objet de $\mathcal{C}'$ puisse être recouvert
par des objets provenant de $\mathcal{C}$. Munissons $\mathcal{C}$ de
la topologie induite. Alors
\begin{itemize}
\item Le foncteur $u$ est continu et cocontinu
\item Si $u_s :\mathcal{C}'^{\widetilde{\;\;}} \ldrt
  \mathcal{C}^{\widetilde{\;\;}}$ est le foncteur $\F\mapsto \F\circ u$
  et $u^s : \mathcal{C}^{\widetilde{\;}}\ldrt
  \mathcal{C}'^{\widetilde{\;\;}}$ le foncteur tel que
$$
\forall \G\in \mathcal{C}^{\widetilde{\;}}\; u^s\G\text{ est le faisceau associé au préfaisceau }
Y\longmapsto \underset{(X,m) \atop {X\in \mathcal{C} \atop m: Y\ldrt
  u(X)}}{\limi} \G (X)
$$
ils induisent une
équivalence de topos
$$
\xymatrix{
\mathcal{C}'^{\widetilde{\;\;}} \ar@<.6ex>[r]^{u_s} &
\ar@<.6ex>[l]^{u^s} \mathcal{C}^{\widetilde{\;\;}}
}
$$
\end{itemize}
\end{theo}

Le corollaire qui suit est scandaleusement absent de SGA4. 

\begin{coro}\label{if2458ihzrmpjgf}
Sous les hypothèses du théorème précédent une famille $(U_\a\ldrt
X)_\a$ de morphismes de $\mathcal{C}$ est couvrante ssi la famille 
$(u(U_\a) \ldrt u(X))_\a$ l'est. 
\end{coro}
\dem
D'après la proposition 1.6. de l'exposé III de SGA4 si $(U_\a\ldrt
X)_\a$ est couvrante alors $(u(U_\a)\ldrt u(X))_\a$ l'est. 
Réciproquement si $(u(U_\a)\ldrt u(X))_\a$ est couvrante soit $R$ le
crible couvrant de $u(X)$ engendré par cette famille. Le foncteur $u$
étant cocontinu la famille de morphismes $Y\ldrt X$ telle qu'il existe
$\a$ et une factorisation $u(Y) \ldrt u(U_\a)\ldrt u(X)$ est un crible
couvrant de $X$. Mais $u$ étant pleinement fidèle l'existence d'une
telle factorisation est équivalente à l'existence d'une factorisation 
$Y\ldrt U_\a \ldrt X$ c'est à dire $Y\ldrt X$ appartient au crible
engendré par $(U_\a\ldrt X)$ qui est donc couvrant. 
\qed

\subsection{Définitions}

{\it Convention : Désormais on fixe une petite sous-catégorie pleine de la
  catégorie des schémas formels $\pi$-adiques quasicompacts, telle que
  tous les schémas formels avec lesquels nous travaillerons soient
  dans
  cette petite catégorie. On vérifie aisément que cela est possible
  puisque nous travaillerons avec des limites projectives 
indexées par des ensembles fixés de schémas
  formels topologiquement de type fini, des
  morphismes topologiquement de type fini au dessus de ces
  éclatements... En particulier les catégories sous-jacentes de 
tous les sites avec lesquels nous
  travaillerons seront petites.
}

\begin{defi}
Soit $\X$ $\pi$-adique sans $\pi$-torsion quasicompact et $\E_{\X^{rig}}$ la
catégorie des morphismes de type $(\E)$ vers $\X$ localisée
relativement aux éclatements formels admissibles (cf. section
\ref{kfusyttgtnbgg}). On munit $\E_{\X^{rig}}$ de la topologie
engendrée par les familles finies de morphismes $(\mathfrak{U}^{rig}_\a
\ldrt \YY^{rig})_\a$
$$
\xymatrix@R=4mm@C=4mm{
\mathfrak{U}_\a^{rig} \ar[rd] \ar[rr] && \YY^{rig} \ar[ld] \\
 & \X^{rig}
}
$$
telles que l'application continue
$\coprod_\a |\mathfrak{U}_\a^{rig}|\ldrt |\YY^{rig}|$
entre les espaces de Zariski-Riemann soit surjective. On note $\X_{\E-rig-\et}$ le site associé. 
\end{defi}

\begin{rema}\begin{itemize}
\item
On remarquera qu'il s'agit d'une topologie qui n'est pas définie à
partir d'une prétopologie. En particulier les familles couvrantes pour
cette topologie n'ont à priori aucune description concrète ! Le problème vient que les morphismes dans $\X_{\E-rig-\et}$ ne sont pas forcément quarrables.
\item Si $f:\YY\ldrt \X$ est un morphisme de schémas formels quasicompacts $\pi$-adiques sans $\pi$-torsion alors $f$ induit un foncteur entre les catégories sous-jacentes de celle de $\X_{\E-rig-\et}$ vers celle de $\YY_{\E-rig-\et}$. Ce foncteur transforme familles couvrantes en familles couvrantes. Néanmoins, il n'y a pas de raison pour que ce foncteur soit continu ! Et même s'il l'est il n'y a pas de raison pour qu'il induise un morphisme de topos c'est à dire le foncteur $f^*$ induit au niveau des faisceaux commute aux limites projectives finies !
\end{itemize}
\end{rema}

\begin{defi}
Soit $\X$ un schéma formel admissible quasicompact. On note $\X_{rig-\et}$ le site
des espaces  rigides
quasicompacts
rig-étales au dessus
de $\X^{rig}$ muni de la 
topologie associée à la prétopologie  telle que $\text{Cov} ( \YY)$ consiste en les familles 
finies $(\mathfrak{U}^{rig}_\a \ldrt \YY^{rig})_\a$  telles que
$\coprod_\a|\mathfrak{U}^{rig}_\a|\ldrt |\YY^{rig}|$ soit surjectif.
\end{defi}

D'après l'exposé II de SGA4 les familles couvrantes de 
$\X_{rig-\et}$ sont les familles $(\mathfrak{U}^{rig}_\a\ldrt
\YY^{rig})_{\a}$ telles qu'il existe une sous-famille finie
surjective au niveau des espaces de Zariski-Riemann.

On renvoie à \cite{BLIII}, \cite{Hu1} ou \cite{Fuji2} pour les propriétés de base des morphismes étales
entre espaces rigides usuels. On retiendra en particulier que l'image
d'un morphisme étale est un ouvert quasicompact de l'espace de
Zariski-Riemann (cela peut se vérifier aisément en utilisant le
théorème de platification de Raynaud-Gruson).
 De cela on déduit que les familles couvrantes de 
$\X_{rig-\et}$ sont les familles $(\mathfrak{U}^{rig}_\a\ldrt
\YY^{rig})_{\a}$ qui induisent une surjection au niveau des espaces de Zariski-Riemann.
 En particulier
le topos associé est cohérent.

\subsection{Lien entre les sites $\X_{\E-rig-\et}$ et $\X_{rig-\et}$ pour $\X$ admissible}

\begin{prop}
Soit $\X$ admissible. Le foncteur d'inclusion  
$$u: \X_{\E-rig-\et}\hookrightarrow \X_{rig-\et}$$
est continu et est tel que la topologie induite par $u$  sur la catégorie 
sous-jacente à $\X_{\E-rig-\et}$ soit la topologie de $\X_{\E-rig-\et}$. 
Les familles couvrantes de $\X_{\E-rig-\et}$ sont les familles de morphismes
$(\mathfrak{U}_\a^{rig}\ldrt \YY^{rig})_\a$ telles qu'il existe une sous-famille finie 
induisant une surjection au niveau des espaces de Zariski-Riemann. Ce foncteur induit une équivalence de topos
$$
(u^s,u_s): (\X_{\E-rig-\et})^{\widetilde{\;\;}} \iso (\X_{rig-\et})^{\widetilde{\;\;}}
$$
\end{prop}
\dem
D'après le théorème \ref{rhfueo5pjgf} le foncteur $u$ satisfait aux hypothèses du théorème \ref{odpmmlfjhjutrr}. Soit $\mathcal{T}$ la topologie sur la catégorie sous-jacente à $\X_{\E-rig-\et}$ induite par $u$ et la topologie de $\X_{rig-\et}$. D'après le corollaire \ref{if2458ihzrmpjgf} les familles couvrantes pour $\mathcal{T}$ sont les familles dont on peut extraire une sous-famille finie surjective au niveau des espaces de Zariski-Riemann. Les cribles engendrés par de telles familles sont les cribles contenant les cribles engendrés par les familles finies surjectives au niveau des espaces de Zariski-Riemann, c'est à dire ceux utilisés pour définir la topologie engendrée de $\X_{rig-\et}$. Donc $\mathcal{T}$ coïncide avec la topologie de $\X_{\E-rig-\et}$. 
\qed

\subsection{Le théorème principal sur la décompletion des topos
  rig-étales}

Soit $(I,\leq)$ un ensemble ordonné cofiltrant et $(\X_i)_{i\in I}$
un système projectif de schémas formels admissibles quasicompacts dont les
morphismes de transition sont affines. On note $\X_\infty 
=\underset{i\in I}{\limp} \X_i$. 

\begin{theo}\label{jkfjdgyuzr345Yfg}
Les familles couvrantes de $(\X_\infty)_{\E-rig-\et}$ sont les familles $(\mathfrak{U}_\a^{rig}\ldrt \YY^{rig})_\a$ possédant une sous-famille finie induisant une surjection au niveau des espaces de Zariski-Riemann. De plus il y a une équivalence de topos
$$
(\X_\infty)_{\E-rig-\et}^{\widetilde{\;\;}} \simeq \underset{i\in I}{\limp} (\X_i)_{rig-\et}^{\widetilde{\;\;}}
$$
\end{theo}
\dem
Le site fibré $i\longmapsto (\X_i)_{rig-\et}$ satisfait aux hypothèses de la section 8.3.1 de \cite{SGA4_exp6}. On en déduit d'après la proposition 8.3.6 de \cite{SGA4_exp6} que l'on a une prétopologie explicite sur le site $\underset{i\in I}{\limi} (\X_i)_{rig-\et}$ obtenue à partir des prétopologies sur les $(\X_i)_{rig-\et}$ lorsque $i$ varie. 
\\
Si $Y\in \underset{i\in I}{\limi} (\X_i)_{rig-\et}$,  pour cette
prétopologie $\text{Cov} ( Y)$ consiste en les familles finies de
morphismes dans $\underset{i\in I}{\limi} (\X_i)_{rig-\et}$ isomorphes
dans $\underset{i\in I}{\limi} (\X_i)_{rig-\et}$ à une famille finie
de morphismes dans un $(\X_i)_{rig-\et}$ qui forme un recouvrement
dans le site $(\X_i)_{rig-\et}$.
\\
Les familles couvrantes de la topologie de $\underset{i\in I}{\limi} (\X_i)_{rig-\et}$ sont donc les familles de morphismes dans $\underset{i\in I}{\limi} (\X_i)_{rig-\et}$ possédant une sous-famille finie isomorphe à une famille finie couvrante d'un $(\X_i)_{rig-\et}$ pour un $i\in I$. 
\\
Par exemple pour $i\in I$ 
 une famille finie de morphismes dans $(\X_i)_{rig-\et}$, $(U_\a\ldrt Y)_\a$, devient couvrante dans $\underset{i\in I}{\limi} (\X_i)_{rig-\et}$ ssi il existe $j\geq i$ tel que la famille tirée en arrière
$(U_\a\times_{\X_i^{rig}} \X_j^{rig} \ldrt Y\times_{\X_i^{rig}} \X_j^{rig} )_\a$ devient couvrante dans $(\X_j)_{rig-\et}$. 
\\

Rappelons que l'on note $\E_{\X^{rig}}$ la catégorie sous-jacente au site $\X_{\E-rig-\et}$. Il y a un foncteur pleinement fidèle
$$
u: \underset{i\in I}{\limi} \E_{\X_i^{rig}} \ldrt \underset{i\in I}{\limi} (\X_i)_{rig-\et}
$$
auquel on peut appliquer le théorème \ref{odpmmlfjhjutrr} et le corollaire \ref{if2458ihzrmpjgf}. Soit $\mathcal{T}$ la topologie induite sur $\underset{i\in I}{\limi} \E_{\X_i^{rig}}$ par $u$ et la topologie précédente sur $\underset{i\in I}{\limi} (\X_i)_{rig-\et}$.
D'après le corollaire \ref{if2458ihzrmpjgf} et la description précédente des familles couvrantes de $ \underset{i\in I}{\limi} (\X_i)_{rig-\et}$ on a une description concrète des familles couvrantes du site $(\underset{i\in I}{\limi} \E_{\X_i^{rig}},\mathcal{T})$. Il s'agit maintenant de vérifier que via l'équivalence de catégories du théorème \ref{mdpojkzrnjznhfbeyty24}
$$
\underset{i\in I}{\limi} \E_{\X_i^{rig}} \iso \E_{\X_\infty^{rig}}
$$
les familles couvrantes se correspondent. Mais cela résulte de la proposition qui suit.
\qed

\begin{prop}
Soit $i\in I$ et $(U_\a\ldrt Y)_\a$ une famille finie de morphismes dans $\E_{\X_i^{rig}}$. Supposons que la famille
$$
(U_\a\times_{\X_i^{rig}}\X_\infty^{rig} \ldrt Y\times_{\X_i^{rig}} \X_\infty^{rig})_\a
$$
induise une surjection au niveau des espaces de Zariski-Riemann
$$
\coprod_\a |U_\a\times_{\X_i^{rig}}\X_\infty^{rig} |\twoheadrightarrow | Y\times_{\X_i^{rig}} \X_\infty^{rig}|
$$
Il existe alors $j\geq i$ tel que 
$$
\coprod_\a |U_\a\times_{\X_i^{rig}} \X_j^{rig} | \ldrt |Y\times_{\X_i^{rig}} \X_j^{rig} |
$$
soit surjectif.
\end{prop}
\dem
On a 
$$
|Y\times_{\X_i^{rig}} \X_\infty^{rig}| = \underset{j\geq i}{\limp} |Y\times_{\X_i^{rig}}|
$$
(utiliser la description des espaces de Zariski-Riemann en termes de
points à valeurs dans des anneaux de valuation rigides). 
Soit 
$$
\forall j\geq i\; \psi_{ji} : |Y\times_{\X_i^{rig}} \X_{j}^{rig}| \ldrt |Y|\text{ et }
f:\coprod_\a |U_\a|\ldrt |Y|
$$
Par hypothèse  $\dpt{\bigcup_{j\geq i} \text{Im } \psi_{ji} \subset \text{Im }f}$. Le morphisme $\coprod_\a U_\a \ldrt Y$ étant rig-étale $\text{Im }f$ est ouvert quasicompact. De plus
$\forall j\; \text{Im }\psi_{ji}$ est pro-constructible. On conclut
comme dans la démonstration de la proposition \ref{kdvyyzrp46EGodg} que
$$
\exists j\geq i \; \text{Im }\psi_{ji} \subset \text{Im } f
$$
\qed

\begin{coro}
Le topos limite projective $\underset{i\in I}{\limp}
(\X_{i})_{rig-\et}^{\widetilde{\;\;}}$ ne dépend que de
$\X_\infty$. De plus si $\La$ est un anneau il y a un isomorphisme
canonique
$$
\underset{i\in I}{\limi} R\GG ( (\X_i)_{rig-\et}, \La ) \simeq R\GG ( (\X_\infty)_{\E-rig-\et},\La)
$$
dans $\DD^+ ( \La-\text{Mod})$ : la limite inductive de la cohomologie
étale de la fibre générique des $\X_i$ ne dépend que
$\X_\infty$. 
\end{coro}
\dem
L'assertion concernant la cohomologie résulte du théorème 8.7.3 de
\cite{SGA4_exp6} dont les hypothèses sont vérifiées grâce aux
propriétés de cohérence des topos rigides étales usuels.
\qed 

\begin{coro}\label{DSSBeztez46I}
Soit $\X$ un schéma formel $\pi$-adique sans $\pi$-torsion quasicompact. La famille de morphisme $(\mathfrak{V}_\a^{rig}\ldrt \YY^{rig})_\a$ dans $\X_{\E-rig-\et}$ est couvrante ssi on peut en extraire une famille finie induisant une surjection au niveau des espaces de Zariski-Riemann.
\end{coro}
\dem
Si $\X$ est affine il est facile de vérifier qu'il s'écrit comme limite projective cofiltrante de schémas formels affines topologiquement de type fini sur $\spf (\O_K)$. Dans ce cas là le résultat est donc une conséquence du théorème \ref{jkfjdgyuzr345Yfg}. Pour $\X$ général soit $(\mathfrak{U}_i)_{i\in I}$ un recouvrement affine fini de $\X$. La famille $(\mathfrak{U}_i^{rig}\ldrt \X^{rig})_{i\in I}$ forme une famille couvrante de l'objet final de $\X_{\E-rig-\et}$. Une famille $(\mathfrak{V}_\a^{rig} \ldrt \YY^{rig})_\a$ dans $\X_{\E-rig-\et}$ est donc couvrante ssi $\forall i\in I \; (\mathfrak{V}_\a^{rig}\times_{\X^rig} \mathfrak{U}_i^{rig})_\a$ l'est. On est donc ramené au cas affine car l'application
$$
\coprod_{i\in I} |\YY^{rig}\times_{\X^{rig}} \mathfrak{U}_i^{rig} | \twoheadrightarrow |\YY^{rig}|
$$
est surjective.
\qed

\begin{rema}
L'auteur de cet article met son lecteur au défi de trouver une démonstration du corollaire précédent sans passer par la théorie des espaces rigides des schémas formels topologiquement de type fini sur $\spf (\O_K)$. 
\end{rema}

\section{Le topos rig-étale d'un schéma formel $\pi$-adique
  non-quasicompact}
\subsection{Le topos}

Soit $\X$ un schéma formel $\pi$-adique sans $\pi$-torsion
non nécessairement quasicompact.

On note $\E_{\X}$ la catégorie des morphismes $\YY\ldrt \X$ qui
peuvent s'écrire sous la forme
$$
\YY\xrig{\; f\; }\mathfrak{U} \hookrightarrow \X
$$
où $\mathfrak{U}$ est un ouvert quasicompact de $\X$ et $f$ est de
type $(\mathcal{E})$. On appellera les morphisme $\YY\ldrt \X$ dans
$\E_\X$ les morphismes de type $(\E)$. 

\begin{lemm}
Soit $\YY\ldrt \mathfrak{U}$ un morphisme de type  $(\mathcal{E})$
vers un ouvert quasicompact de $\X$ et $\mathfrak{V}$ un ouvert
quasicompact de $\X$ contenant $\mathfrak{U}$. Alors le composé 
$\YY\ldrt \mathfrak{U} \hookrightarrow \mathfrak{V}$ est de type  $(\mathcal{E})$.
\end{lemm}
\dem
Utiliser le lemme \ref{kfuspizzz}.
\qed

On vérifie alors que $\E_\X$ est équivalente à 
$$
\underset{\mathfrak{U}\subset \X\atop \mathfrak{U} \text{
    quasicompact}}{\limi} \E_\mathfrak{U}
$$
limite inductive sur les ouverts quasicompacts de $\X$ des morphismes
de type $(\mathcal{E})$ au dessus d'un tel ouvert.

Plus généralement si $\E_{\X^{rig}}$ désigne la catégorie $\E_\X$
localisée relativement aux éclatements formels admissibles alors 
$$
\E_{\X^{rig}}\simeq 
\underset{\mathfrak{U}\subset \X\atop \mathfrak{U} \text{
    quasicompact}}{\limi} \E_{\mathfrak{U}^{rig}}
$$

\begin{defi}
On note $\X_{\E-rig-\et}$ le site donc la catégorie sous-jacente est
$\E_{\X^{rig}}$
et dont la topologie est engendrée par les familles finies de
morphismes surjectives au niveau des espaces de Zariski-Riemann. 
\end{defi}

Si $(\mathfrak{V}_\a\ldrt \YY)_\a$ est une famille de morphismes dans
$\E_{\X^{rig}}$ et si le morphismes structurel $\YY^{rig}\ldrt \X^{rig}$ se factorise par
l'ouvert quasicompact $\mathfrak{U}\subset \X$ il en est de même des
$\mathfrak{V}_\a^{rig}\ldrt \X^{rig}$. De cela et du corollaire
\ref{DSSBeztez46I} on déduit  le
lemme qui suit.

\begin{lemm}
Une famille de morphismes dans $\X_{\E-rig-\et}$ est couvrante ssi
il existe une sous-famille finie surjective au niveau des espaces de Zariski-Riemann.
\end{lemm}

De là on déduit aisément.

\begin{prop}
Le topos $(\X^{rig})_{\E-rig-\et}\top$ s'identifie au topos
$$
\underset{\mathfrak{U}\subset \X\atop \mathfrak{U} \text{
    quasicompact}}{\limp} \mathfrak{U}\top_{\E-rig-\et}
$$
Si $(\mathfrak{U}_i)_{i\in I}$ est un recouvrement de $\X$ par des
ouverts quasicompacts 
alors le topos $(\X^{rig})_{\E-rig-\et}\top$ s'identifie au topos recollé
des $(\mathfrak{U}_i^{rig})\top_{\E-rig-\et}$ le long des
$\mathfrak{U}_i\cap \mathfrak{U}_j$, $i,j\in I$ i.e. aux objets
cartésiens du topos total du diagramme 
$$
\xymatrix{
\coprod_{i,j,k} (\mathfrak{U}_i\cap \mathfrak{U}_j\cap
\mathfrak{U}_k)\top_{\E-rig-\et} \ar@<.8ex>[r]\ar[r] \ar@<-.8ex>[r]
& \coprod_{i,j}  (\mathfrak{U}_i\cap \mathfrak{U}_j)\top_{\E-rig-\et}
\ar@<.6ex>[r] \ar@<-.6ex>[r] & \coprod_{i} (\mathfrak{U}_i)\top_{\E-rig-\et}
}
$$
\end{prop}

\subsection{Cohomologie à support compact}

On définit ici la cohomologie à support compact de $\X^{rig}$. On
prendra garde que la définition donnée bien que suffisante pour nos
besoins n'est pas la bonne en général. Par exemple pour les schémas
formels admissibles ce n'est la bonne définition que pour ceux dont
les composantes irréductibles de la fibre spéciale sont propres, par
exemple le schéma formel de Deligne-Drinfeld $\widehat{\Omega}$.

\begin{defi}
Soit $\F\in \X\top_{\E-rig-\et}$. On note $\GG_! (\X^{rig},\F)$ les
sections de $\GG ( \X^{rig},\F) = \underset{\mathfrak{U} 
  \text{ quasicompact}}{\limp} \GG ( \mathfrak{U}^{rig},\F)$ dont le
support est contenu dans un ouvert quasicompact de $\X^{rig}$. 
\end{defi}

Si $|\X^{rig}| = \underset{\mathfrak{U}\text{ q.c.}}{\limi}
|\mathfrak{U}^{rig}|$ désigne l'espace de Zariski-Riemann de $\X$ il y
a alors une projection du topos rig-étale vers le topos
de l'espace topologique $|\X^{rig}|$, $\pi:
(\X^{rig})\top_{\E-rig-\et} \ldrt |\X^{rig}|\top$. Alors $\GG_!
(\X^{rig},-) =\GG_c ( |\X^{rig}|,-)\circ \pi_*$ où $\GG_c (
|\X^{rig}|,-)$ désigne les sections dont le support est un fermé
quasicompact de l'espace de Zariski-Riemann.

\begin{prop}
Soit $\F\in \La-(\X^{rig})\top_{\E-rig-\et}$. 
Soit $(\mathfrak{U}_i)_{i\in I}$ un recouvrement localement fini de $\X$ par des
ouverts quasicompacts. Notons $\forall \a\subset I\; \mathfrak{U}(\a)
= \cap_{i\in \a} \mathfrak{U}_i$. 
 Il y a alors une suite spectrale 
$$
E^{pq}_1= \bigoplus_{\a\subset I\atop |\a|=p+1} H^q (
\mathfrak{U}(\a)^{rig}, \F) \limpl H^{p+q}_! ( \X^{rig},\F)
$$
\end{prop}
\dem
Considérons le morphisme 
$$
j: \coprod_{i\in I} \mathfrak{U}_i \ldrt \X
$$
Si $\G$ est un faisceau sur $\X_{\E-rig-\et}$ il donné naissance à une
résolution co-simpliciale $\G\ldrt C^\bullet (\G)$ fonctorielle en
$\G$ où
$$
C^p (\G) = (j_*j^*)^{p+1} (\G) = \prod_{\a\subset I\atop |\a|=p+1} j_{\a*} (\G_{|\mathfrak{U}(\a)^{rig}})
$$
avec $j_\a: \mathfrak{U} (\a)^{rig}\hookrightarrow \X^{rig}$. De plus
si $\G$ est injectif alors $C^p (\G)$ l'est également. De l'égalité 
$$
\GG_! ( \X^{rig}, C^p (\G)) = \bigoplus_{\a\subset I\atop |\a|=p+1}
\GG ( \mathfrak{U} ( \a)^{rig},\G)
$$
on conclut facilement.
\qed

\begin{coro}\label{DGkygygez380EZDP}
Tout faisceau flasque sur $\X_{\E-rig-\et}$ est $\GG_! (\X^{rig},-)$-acyclique.
\end{coro}
\dem
Le schéma formel $\X$ étant quasi-séparé il existe un recouvrement
comme dans la proposition précédente.
Dans la suite spectrale précédente les morphismes $\mathfrak{U}
(\a)^{rig}\ldrt \X^{rig}$ sont alors dans $\X_{\E-rgi-\et}$.
\qed

On retiendra donc que l'on peut calculer la cohomologie à support
compact à l'aide de résolutions flasques. 

\section{Le formalisme des faisceaux équivariants lisses}\label{lsfkjsut25tggpoi}

\subsection{Hypothèses}\label{lfizyr25si}

Soit $G$ un groupe topologique possédant un sous-groupe ouvert
profini. Soit $\La$ un anneau. 
\\
Soit $\mathcal{R}$ une catégorie et $X$ un objet de $\mathcal{R}$, ou
plus généralement de $\mathcal{R}^{\widehat{\;\;}}$, la catégorie des
préfaisceaux sur $\mathcal{R}$, 
 muni
d'une action du groupe ``abstrait'' $G$.
\\
On suppose donné une sous-catégorie pleine $\mathcal{C}$ de la
catégorie $\mathcal{R}/X$ des morphismes vers $X$ et une topologie sur
$\mathcal{C}$ faisant donc de $\mathcal{C}$ un site. On supposera
également que 
tout objet de $\CC$ est quasicompact. On suppose de plus que $\forall
g\in G$ et $(U\ldrt X)\in \CC$  le produit cartésien du diagramme
$$
\xymatrix@R=4mm@C=4mm{
 & U \ar[d] \\
X\ar[r]^g & X
}
$$ 
existe dans $\mathcal{R}$ et est dans $\CC$. On notera $g^{-1} (U)$ ce produit cartésien.
 On suppose que $U\mapsto g^{-1} (U)$ 
 transforme familles couvrantes en familles couvrantes et définit donc
 un isomorphisme du site $\CC$ dans lui même. 
 Ainsi
il y a une ``action'' de $G$ sur le site $\CC$.  On notera $(g^*,g_*)$
le morphisme de topos associé où $g^* \F (U) = \F (g (U))$. 
 On a des identifications canoniques
$\forall g,g'\in G\; (gg')^* = g'^* g^*$. 
Nous ferons l'hypothèse suivante : 
\\

{\it Hypothèses de continuité :  
\`A tout $(U\ldrt X)\in \mathcal{C}$ est associé ``canoniquement'' un sous-groupe compact ouvert $K_U\subset G$ 
et un relèvement $(\beta_k)_{k\in K_U}$ de l'action de $K_U$ sur $X$ à $U$
$$
\forall k\in K_U\; 
\xymatrix@R=3mm@C=6mm{
U\ar[r]^{\beta_k} \ar[d] & U\ar[d] \\ 
X \ar[r]^k & X
}
$$
où $\beta_k\circ \beta_{k'}=\beta_{kk'}$. Le mot ``canoniquement'' signifie que si 
$\xymatrix@R=3mm@C=3mm{V\ar[rr]^f \ar[rd] && U\ar[ld]\\ & X}$ est un morphisme dans $\CC$  il existe alors
un sous-groupe ouvert de $K_U\cap K_V$ en restriction auquel $f$ est compatible aux relèvements  
de l'action à $U$ et $V$. 
}

\subsection{$G$-faisceaux lisses}

Rappelons que par définition un $G$-faisceau (resp. préfaisceau)
 sur $\mathcal{C}$ est un
faisceau $\F\in \CC\top$  (resp. $G$-préfaisceau $\F\in \CC^{\widehat{\;\;}}$)
muni d'isomorphismes 
$$
\forall g\in G\; c_g : g^*\F\iso \F
$$
vérifiant la condition de cocyle $\forall g,g'\in G\;\; c_{g'}\circ
g'^* c_g = c_{gg'}$. 

Si $\F$ est un $G$-faisceau (resp. préfaisceau) et $U\in \mathcal{C}$ 
$$
\forall k\in K_U\;\; (k^* \F)(U) = \F (k (U)) = \F (U)
$$
et donc $K_U$ agit sur $\F(U)$. De plus si $V\ldrt U$ est un morphisme
dans $\mathcal{C}$ l'application de restriction $\F (U) \ldrt \F(V)$
est $K$-équivariante pour $K$ compact ouvert suffisamment petit.

\begin{defi}
Un $G$-faisceau (resp. préfaisceau) $\F$ sera dit lisse si $\forall
U\in\mathcal{C}$ l'action de $K_U$ sur $\F (U)$ est lisse, i.e. le
stabilisateur de tout élément est un sous-groupe ouvert.
\end{defi}

\begin{defi}
On note $\CC\top_G$, resp. $\CC\top_{G-\lss}$, la catégorie des
$G$-faisceaux, resp. des $G$-faisceaux lisses. On note
$\La  -\CC\top_{G}$, $\La -\CC\top_{G-\lss}$ les catégories de faisceaux de
$\La$-modules associées.
\end{defi}

\begin{exem}
Si l'action de $G$ sur $X$ est triviale $\La-\CC\top_G$
s'identifie à la catégorie $\La[G]-\CC\top$ et $\La-\CC\top_{G-\lss}$
à la sous-catégorie de $\La[G]-\CC\top$ des faisceaux $\F$ tels que
$\forall U\in \CC\; \F(U)$ soit un $G$-module lisse. Lorsque $\mathcal{R}$ est
la catégorie à un seul élément et un seul morphisme on retrouve les
$G$-modules lisses utilisés dans la théorie de la représentation des
groupes $p$-adiques ou en cohomologie galoisienne.
\end{exem}

\subsection{Les différentes opérations reliant $G$-faisceaux,
  $G$-faisceaux lisses et faisceaux}

On notera
$$
\xymatrix@R=5mm@C=5mm{
\CC\top_{G-\lss} \ar@{^(->}[rr]^i \ar[rd]_k && \CC\top_G \ar[ld]^j \\
& \CC\top
}
$$
les foncteurs évidents où $k$ et $j$ sont l'oubli de l'action de $G$
et $i$ est le plongement tautologique. On notera de la même façon le
diagramme associé pour les faisceaux de $\La$-modules.

\subsubsection{Faisceau associé à un préfaisceau}

Bien sûr si $\F$ est un $G$-préfaisceau le faisceau associé est un
$G$-faisceau.

\begin{lemm}\label{kfjsfiz24694zf}
 Le faisceau associé
à un $G$-préfaisceau lisse est un $G$-faisceau lisse.
\end{lemm}
\dem
Par définition (cf. SGA4 exposé II section 3) le faisceau associé à
$\F$ est $L(L\F)$ où l'opération $\F\mapsto L\F$ est
$$
L\F (U) =\underset{R\in J(U)}{\limi} \Hom (R,\F)
$$
où $J(U)$ désigne l'ensemble des cribles couvrants de $U$. Mais $U$
étant quasicompact l'ensemble des cribles engendrés par des familles
couvrantes finies est cofinal dans $J(X)$. La limite inductive
précédente peut donc se calculer sur ces cribles là. Or si le crible
$R$ est engendré par la famille finie $(V_\a\ldrt U)_\a$ alors
$$
\Hom (R,\F) \subset \prod_\a \F (V_\a)
$$
et étant donné $(s_\a)_\a\in  \prod_\a \F (V_\a)$ 
il existe un sous-groupe ouvert agissant discrètement sur tous les
$s_\a$.
\qed

\begin{exem}
Sous l'hypothèse du lemme précédent le faisceau constant
$\underline{\La}$ est lisse.
\end{exem}

Il est clair que la catégorie $\La-\CC\top_G$ est abélienne et que le
foncteur $j:\La-\CC\top_G\ldrt \La-\CC\top$ est exact i.e. les noyaux et
conoyaux de morphismes de $G$-faisceaux sont les noyaux et conoyaux
des morphismes de faisceaux sous-jacents qui sont naturellement munis
d'une structure de $G$-faisceau. 

\begin{coro}\label{kfuyzr35fzfr46}
La catégorie $\La-\CC\top_{G-\lss}$ des $G$-faisceaux de
$\La$-modules lisses est une sous-catégorie abélienne de
$\La-\CC\top_G$ la catégorie de $G$-faisceaux de $\La$-modules.
\end{coro}
\dem
Il est clair qu'un sous-$G$-faisceau d'un $G$-faisceau lisse est lisse
et que donc les noyaux de morphismes dans $\CC\top_G$ entre deux $G$-faisceaux lisses dans
sont lisses. Pour les conoyaux il suffit d'appliquer le lemme
précédent.
\qed

\subsubsection{Inductions brutales}

Soit $\F$ un faisceau. On note 
$$
\text{Ind} \F= \prod_{g\in G} \F \;\;\;\;\;\; \text{ind}\F= \bigoplus_{g\in
  G} \F
$$
les $G$-faisceaux munis de l'action de $G$ par permutations (la somme
directe définissant $\text{ind} \F$ est prise dans la catégorie des
faisceaux, c'est le faisceau associé au préfaisceau somme directe). 
Cela définit des foncteurs $\xymatrix{\CC\top
  \ar@<.6ex>[r]^{\text{Ind}} \ar@<-.6ex>[r]_{\text{ind}} &
  \CC\top_G}$. 
Ils forment avec $j$ un triplet de foncteurs adjoints
\begin{eqnarray*}
\forall \F\in \CC\top\;\forall \G\in \CC\top_G \;\; \Hom_{\CC\top}
(j\G,\F) &\simeq & \Hom_{\CC\top_G} (\G, \text{Ind}\F) \\
\Hom_{\CC\top} ( \F, j\G) &\simeq & \Hom_{\CC\top_G} ( \text{ind} \F, \G)
\end{eqnarray*}
Il en est de même avec les faisceaux de $\La$-modules.
En particulier on déduit de ces adjonctions que $\text{Ind}$ et $j$
transforment faisceaux injectifs de $\La$-modules en faisceaux
injectifs. De plus $\text{Ind}$ est exact à gauche (mais pas exact en
général). De cela on déduit que la catégorie abélienne $\La-\CC\top_G$
 possède suffisamment d'injectifs (bien sûr le lecteur
avisé sait que $(\CC\top_G,\underline{\La})$ est un topos annelé et
que donc $\La-
\CC\top_G$ possède suffisamment d'injectifs, mais $\text{Ind}$
permet de construire explicitement de tels objets injectifs à partir
de ceux de $\CC\top$). 
\\

Si $\F\in \La-\CC\top_G$ et $e$ désigne un objet de $\CC\top$ auquel
l'action de $G$ se relève i.e. un élément de $\CC\top_G$ (par exemple l'objet final) 
le complexe $R\GG (e, j\F) \in \DD^+ (\La)$ est muni d'une action de
$G$. Cependant si 
$$
F= \GG(e,-)\circ j  : \La -\CC\top_G \ldrt \La [G]-\text{Mod}
$$
étant donné que $j$ est exact et transforme injectifs en injectifs
$RF=R\GG(e,-)\circ j$ et donc le complexe $R\GG(e,j\F)$ provient
canoniquement d'un complexe dans $\DD^+ (\La [G])$ via le
foncteur $\DD^+ (\La[G]) \ldrt G-\DD^+ (\La)$, $ G-\DD^+ (\La)$
désignant les objets de la catégorie dérivée $\DD^+ (\La)$ munis d'une
action de $G$.

L'un des buts de  ce chapitre est de faire de même avec les faisceaux lisses.

\subsubsection{Lissification}

Soit $\F\in \CC\top_G$. On note $\F^{\lss}$ le préfaisceau
$$
U\longmapsto \{ \text{sections lisses de } \F (U) \}
$$
dont on vérifie facilement en utilisant l'hypothèse de quasicompacité
des objets de $\CC$ 
que c'est un faisceau. Cela définit un foncteur 
$$
(-)^{\lss}: \CC\top_G \ldrt \CC\top_{G-\lss} 
$$
adjoint à droite au foncteur $i$ 
$$
\forall \F\in \CC\top_{G-\lss}\; \forall \G\in \CC\top_G\;\;
\Hom_{\CC\top_G} ( i\F,\G) \simeq \Hom_{\CC\top_{G-\lss}} ( \F,\G^{lss})
$$
De cela on déduit que pour les faisceaux de $\La$-modules $(-)^{\lss}$
transforme les injectifs en injectifs et est exact à gauche.
On en déduit le lemme suivant.

\begin{lemm}\label{kfusfufuzr45r}
La catégorie abélienne $\La-\CC\top_{G-\lss}$
possède suffisamment d'injectifs.
\end{lemm}

En résumé il y a un diagramme
$$
\xymatrix@R=9mm@C=9mm{
\CC\top_{G-\lss} \ar@{^(->}[rr]^i \ar[rd]_k && \CC\top_G \ar[ld]|j \ar@(lu,ru)[ll]_{(-)^{\lss}}
\\
& \CC\top  \ar@/^.9pc/[ru]|(.4){\text{Ind}} \ar@/_1.3pc/[ru]|{\text{ind}}
}
$$

\begin{rema}
Le foncteur composé $(-)^{\lss}\circ \text{Ind}$ permet donc de construire suffisamment d'injectifs dans $\La-\CC\top_{G-\lss}$. Néanmoins en général l'image par $k$ d'un tel objet injectif n'est pas un faisceau  injectif.
\end{rema}

\subsubsection{Induction/restriction}\label{kfjduytb35sfo}

Soit $F\in \CC^{\widehat{\;\;}}$ un préfaisceau. Rappelons (SGA4 exposé III section 5) que l'on a un site $\CC/F$ dont la topologie est induite par celle de $\CC$ via le foncteur canonique $\CC/F\ldrt \CC$. Rappelons également que $\CC\top/F$ s'identifie à $(\CC/F)\top$. Il y a un triplet de foncteurs adjoints
$$
\xymatrix@C=13mm{
\CC\top/F = (\CC/F)\top \ar@<1.5ex>[r]^(.7){j_*} \ar@<-1.5ex>[r]_(.7){j_!} & \ar[l]|(.3){j^*} \CC\top
}
$$
où $j^*$ est le foncteur de restriction aux objets au dessus de $F$ et $j_! \F$ est le faisceau associé au préfaisceau
$$
U\longmapsto \underset{m\in F(U)}{\limi} \F (U,m)
$$
\\

Soit maintenant $K$ un sous-groupe compact ouvert de $G$ et $F\in \CC^{\widehat{\;\;}}_K$ un $K$-préfaisceau i.e. l'action de $K$ sur $X$ s'étend en une action sur $F$ via $F\ldrt X$. On prendra par exemple $U\in \CC$, un sous-groupe compact $K$ tel que l'action de $K$ se relève à $U$  et $F=h_U$ le préfaisceau représenté par $U$. Supposons de plus que $F$ soit un $K$-préfaisceau lisse. Alors la catégorie $\CC/F$ vérifie les mêmes hypothèses que celle de $\CC$ : pour tout $U \ldrt F$ (i.e. un élément de $F(U)$) l'action de $K$ sur $F$ se relève ``canoniquement'' à $U$  sur un sous-groupe ouvert de $K$ et de plus quitte à se restreindre à des sous-groupes ouverts plus petits cette action relevée est fonctorielle dans $\CC/F$. 
\\

On peut donc définir la catégorie $(\CC/F)\top_{K-\lss}$ des
$K$-faisceaux lisses au dessus de $F$. 
 Il y a alors un triplet de foncteurs adjoints 
$$
\xymatrix@C=13mm{
(\CC/F)\top_{K-\lss} \ar@<1.5ex>[r]^(.55){j_{K*}} \ar@<-1.5ex>[r]_(.55){j_{K!}} & \ar[l]|(.45){j_K^*} \CC\top_{G-\lss}
}
$$
où
\begin{itemize}
\item Le foncteur $j_K^*$ est le foncteur de restriction du faisceau aux objets au dessus de $F$ et de l'action de $G$ à $K$
\item Le foncteur $j_{K!}$ est celui qui à $\F\in K-(\CC/F)\top_{disc}$ associe  
$$
\bigoplus_{\Omega} j_! \F
$$
où $\Omega$ désigne un ensemble de représentants dans $G$ de  $K\bc G$, 
dont on vérifie aisément qu'il est muni d'une action lisse de $G$ par ``permutations'' des éléments de $\Omega$.
\item Le foncteur $j_{K*}$  associe à $\F$ le $G$-faisceau lisse 
$
\left ( \prod_{\Omega} j_*\F \right )^{\lss}
$
où $\Omega$ est comme précédemment.
\end{itemize}

\begin{defi}
On notera $j_K^*(-) = \text{Res}^{G/X}_{K/F} (-)$ la restriction, $j_{K!}(-) = \text{c-ind}^{G/X}_{K/F} (-)$ l'induction compacte et
$j_{K*} = \text{Ind}^{G/X}_{K/F} (-)$ l'induite lisse.
\end{defi}

On a les formules usuelles de réciprocité de Frobenius qui traduisent les adjonctions entre nos trois foncteurs 
\begin{eqnarray*}
\Hom_{\CC\top_G}(\text{c-ind}^{G/X}_{K/F} \F, \G)& \simeq &\Hom_{(\CC/F)\top_K } ( \F, \text{Res}^{G/X}_{K/F} \G)  \\
\Hom_{\CC\top_G } (\G, \text{Ind}^{G/X}_{K/F} \F ) &\simeq & \Hom_{(\CC/F)\top_K } (  \text{Res}^{G/X}_{K/F} \G, \F )
\end{eqnarray*}
et les propriétés de transitivité usuelles
lorsque $K'\subset K$, $F'\ldrt F$, l'action de $K'$ se relève à $F'$ en une action discrète et est compatible à celle de $K$ sur $F$, 
 qu'on laisse en exercice au lecteur.

\begin{exem}
Considérons les faisceaux de $\La$-modules. 
Les $G$-faisceaux lisses $\text{c-ind}^{G/X}_{K/U} (\underline{\La})$, pour $U$ dans $\CC$ et $K$ agissant sur $U$ relevant l'action de $G$ sur $X$, vérifient
$$
\Hom_{G} ( \text{c-ind}^{G/X}_{K/U} (\underline{\La}), \F) \simeq \F (U)^K
$$
Ainsi si $\F$ est lisse, $\F(U) = \underset{K}{\limi} \F(U)^K$ et donc les objets 
$ \text{c-ind}^{G/X}_{K/U} (\underline{\La})$ forment un système de générateurs de la catégorie abélienne $\La-\CC\top_{G-\lss}$. 
\end{exem}

\begin{exem}
Soit $(U_i\ldrt Y)_{i\in I}$ une famille finie de morphismes dans $\CC$. On vérifie aisément que le produit fibré $F=\dpt{\underset{i}{\times}_{h_Y}} h_{U_i} \in \CC^{\widehat{\;\;}}$ est un $K$-préfaisceau lisse pour $K$ suffisamment petit. Alors $\text{c-Ind}^{G/X}_{K/F} \underline{\La}$ vérifie
$$
\Hom_G ( \text{c-Ind}^{G/X}_{K/F} \underline{\La}, \F ) \simeq \Hom ( F ,\F)^K
$$
\end{exem}

\subsection{Les G-faisceaux lisses forment un topos}

\begin{prop}
La catégorie des $G$-faisceaux lisses $\CC\top_{G-\lss}$ est un topos.
\end{prop}
\dem
On applique le critère de Giraud (théorème 1.2 de l'exposé IV de SGA4). 
L'existence de limites projectives finies ne pose pas de problème puisque les 
limites projectives finies de $G$-faisceaux lisses dans la catégorie des faisceaux sont encore des $G$-faisceaux lisses. L'existence de sommes directes quelconques ainsi que de quotients par les relations d'équivalence est une conséquence du lemme \ref{kfjsfiz24694zf}. 
L'existence d'une petite famille génératrice se déduit des résultats de la sous-section précédente : les faisceaux lisses $\text{c-ind}^{G/X}_{K/U} e$ où $e$ est l'obet final du topos $(\CC/U)_{K-\lss}\top$ forment un système de générateurs de $\CC_{G-\lss}\top$ puisque pour tout $\F \in \CC_{G-\lss}\top$ 
$$
\Hom ( \text{c-ind}^{G/X}_{K/U} e, \F) = \F(U)^K
$$
et $\F(U)=\underset{K}{\limi} \F (U)^K$.
\qed

\begin{rema}
Le topos $\CC\top_{G-\lss}$ est annelé par $\underline{\La}$. On en déduit donc d'une autre façon que la catégorie $\La-\CC\top_{G-\lss}$ possède suffisamment d'injectifs.
\end{rema}

\begin{rema}
Les limites projectives quelconques, non-finies, dans le topos $\CC\top_{G-\lss}$ se calculent de la façon suivante : $\underset{i\in I}{\limp}\F_i$ dans $\CC\top_{G-\lss}$ est égal à
$(\underset{i\in I}{\limp} \F_i )^{\lss}$, la partie lisse de la limite projective des faisceaux usuels. 
\end{rema}

\subsection{Le complexe de cohomologie d'un $G$-faisceau lisse}

\subsubsection{Compléments sur la cohomologie de Cech}

Dans cette sous-section on oublie momentanément les notations précédentes.
On reprend ici le formalisme du complexe de Chech développé dans l'exposé V
de SGA4 au cas d'un recouvrement $(U_i\ldrt X)$ pour lequel les
produits fibrés $U_{i_1}\times_X\dots\times_X U_{i_p}$ n'existent
pas. 
\\

Soit donc $\CC$ une catégorie, $\La$ un anneau et $(U_i\ldrt X)_i$ une
famille finie de morphismes dans $\CC$. On note $R\in \CC^{\widehat{\;\;}}$
le crible de $X$ engendré par cette famille. Il y a un couple de foncteurs
adjoints
$$
(j_!,j^*) : \coprod_{i\in I} \CC^{\widehat{\;\;}} /U_i \ldrt \CC^{\widehat{\;\;}}/R
$$
On note $\La\in  \CC^{\widehat{\;\;}}/R$ le préfaisceau constant.

\`A un tel couple de foncteurs adjoints est associé un complexe
simplicial de $\La$-modules (cf. \cite{Illusie1} chapitre I section 1.5) $D^\bullet$
$$
\xymatrix{
D^{p} \ar@<1ex>[r] \ar@{..>}[r] \ar@<-1ex>[r]&  D^{p+1}&  \dots
\ar@<1ex>[r] \ar@<-1ex>[r] & D^{0}
}
$$
$$
D^p = \underbrace{(j_!j^*)\circ \dots\circ  (j_!j^*)}_{-p+1\text{ fois}} \La
$$
dont les flèches simpliciales sont données par les applications
d'adjonction.
De plus ce complexe est augmenté ver $\La$, $D^\bullet \ldrt \La$.

\begin{lemm}
$D^\bullet \ldrt \La$ est une résolution projective de $\La$ dans
$\CC^{\widehat{\;\;}}/R$. 
\end{lemm}
\dem
La famille $(U_i\ldrt R)_{i\in I}$ étant épimorphique vers l'objet
final de $\CC^{\widehat{\;\;}}/R$ il suffit de montrer que
$j^*D^\bullet \ldrt j^*\La$ est une résolution. Mais le complexe $[j^*
D^\bullet \ldrt j^*\La]$ est homotope à zéro (cf.  \cite{Illusie1}
chapitre I section 1.5). Les $D^p$ sont projectifs car $j_!$ et $j^*$,
possédant tous deux un adjoint à droite, envoient projectifs sur projectifs.
\qed
\\

Si $\kappa : \CC^{\widehat{\;\;}}/R\ldrt  \CC^{\widehat{\;\;}}$,
$\La_R = \kappa_! \La$ on obtient donc une résolution projective
$$
[C^\bullet \ldrt \La_R ] = \kappa_! [D^\bullet \ldrt \La]
$$
On a concrètement en notant $\forall \{i_1,\dots,i_p\}\subset I$ 
$$
\La_{{h_{U_{i_1}}\times_R \dots \times_R h_{U_{i_p}}}}
=\delta_! \La \; \text{ où } \; \delta : \CC^{\widehat{\;\;}}/h_{U_{i_1}}\times_R
\dots \times_R h_{U_{i_p}} \ldrt \CC^{\widehat{\;\;}}
$$
où comme d'habitude $h_Y$ désigne le préfaisceau représenté par $Y$,
ce qui peut encore se récrire
$$
\La_{{h_{U_{i_1}}\times_R \dots \times_R h_{U_{i_1}}}} =
\La_{U_{i_1}}\otimes_{\La_R} \dots \otimes_{\La_R} \La_{U_{i_p}}
$$ 
$$
C^{-p+1} = \bigoplus_{\{i_1,\dots, i_p\}\subset I}
\La_{h_{U_{i_1}}\times_R \dots \times_R h_{U_{i_p}}}
$$
Si $F\in \La-\CC^{\widehat{\;\;}}$ est un préfaisceau de $\La$-modules
on peut donc calculer
$$
\forall q\geq 0\;\; H^q (R,F) = \text{Ext}^q_\La ( \La_R, F) = H^q
(\Hom_{\La} ( C^\bullet ,F))
$$
où $H^q (R,-)$ désigne la cohomologie dans le topos
$\CC^{\widehat{\;\;}}$. 

Supposons maintenant que $\CC$ est un site et soit $\F\in \La-\CC\top$
un faisceau de $\La$-modules. On note comme dans l'exposé V de SGA4
$\; \mathcal{H}^0 : \CC\top \ldrt \CC^{\widehat{\;\;}}$ le plongement
canonique. On a alors
$$
H^q (R, \mathcal{H}^0 (\F)) = H^q (\Hom (C^\bullet, \mathcal{H}^0
(\F))) = H^q ( \Hom (a C^\bullet, \F))
$$
où $aC^\bullet$ est le complexe de faisceaux associé au complexe de
préfaisceaux $C^\bullet$. Ce complexe se calcule comme précédemment en
remplaçant $\CC^{\widehat{\;\;}}$ par $\CC\top$, les catégories
$\CC^{\widehat{\;\;}}/(-)$ par $\CC\top/(-)$ et les foncteurs
$(j_!,j^*)$ par leurs analogues faisceautiques. 

Rappelons que si $F\in C^{\widehat{\;\;}}$ alors $\CC\top/F =
\CC\top/aF$. \'Etant donné que $aR=X$ on en déduit que 
$$
aC^{-p+1}= \bigoplus_{\{i_1,\dots, i_p \}\subset I}
\underline{\La}_{{h_{U_{i_1}}\times_{h_X} \dots \times_{h_X}
    h_{U_{i_p}}}} = \bigoplus_{\{i_1,\dots, i_p \}\subset I}
\underline{\La}_{U_{i_1}}\underset{\underline{\La}_X}{\otimes} \dots
\underset{{\underline{\La}_X}}{\otimes} \underline{\La}_{U_{i_p}}
$$

\subsection{Le théorème d'acyclicité}\label{KBIGzg3675nt}

On reprend maintenant les hypothèses des sections antérieures.

{\it Hypothèses : On suppose que si $(U_i\ldrt Y)_{i\in I}$ est une famille couvrante finie
  alors $\forall p\;\forall \{i_1,\dots,i_p\}\subset I$ le préfaisceau
  $h_{U_{i_1}}\times_{h_Y}\dots \times_{h_Y} h_{U_{i_p}}$ est
  quasicompact dans $\CC^{\widehat{\;\;}}$ muni de la topologie définie
  dans la section 5 de l'exposé II de SGA4. Cette dernière condition,
  vérifiée par exemple si la topologie sur $\CC$ est induite par une prétopologie,
  signifie qu'il existe une famille finie $(V_\a\ldrt
  h_{U_{i_1}}\times_{h_Y}\dots \times_{h_Y} h_{U_{i_p}})_\a$ où
  $V_\a\in \CC$ telle que $\forall \F\in \CC\top$ $\; \Hom (
  h_{U_{i_1}}\times_{h_Y}\dots \times_{h_Y} h_{U_{i_p}},\F)
  \hookrightarrow \prod_\a \F (V_\a)$.
}

\begin{theo}\label{OOPfyzY136Uetet}
Sous les hypothèses précédentes si $\F\in \La-\CC\top_{G-\lss}$ est
un $G$-faisceau lisse de $\La$-modules injectif alors le faisceau sous-jacent
de $\La$-modules est flasque.
\end{theo}
\dem
Soit $(U_i\ldrt Y)_{i\in I}$ une famille couvrante finie. Soit $R$ le
crible de $Y$ engendré par cette famille. D'après la proposition 4.3
de l'exposé V de SGA4 on doit montrer que 
$$
\forall q>0\; H^q (R,\mathcal{H}^0 (\F)) =0
$$
Utilisons les résultats et notations de la sous-section précédente
pour écrire
$$
H^q (R,\mathcal{H}^0 (\F)) = H^q (\Hom_{\La-\CC\top} (aC^\bullet, \F))
$$
où $aC^{-p+1} =  \bigoplus_{\{i_1,\dots, i_p \}\subset I}
\underline{\La}_{{h_{U_{i_1}}\times_{h_X} \dots \times_{h_X}
    h_{U_{i_p}}}}$. Soit donc $\{i_1,\dots,i_p\}\subset I$ et $F=h_{U_{i_1}}\times_{h_X} \dots \times_{h_X}
    h_{U_{i_p}}\in \CC^{\widehat{\;\;}}$. Considérons des sous-groupes
    compacts ouverts $K$ suffisamment petits tels que l'action de $K$
    sur $X$ s'étende à chacun des $U_i$ et donc à $F$. Un tel
    sous-groupe $K$ agit sur 
$$
  \F (F) = \Hom_{\La-\CC\top} (\underline{\La}_F,\F)
$$
Il y a une application injective 
$$
\underset{K}{\limi} \F (F)\hookrightarrow \F (F)
$$
Cette application est bijective car si $(V_\a\ldrt F)_\a$ est une
famille finie comme dans les hypothèses de cette section alors 
$\F (F) \hookrightarrow \prod_{\a} \F (V_\a)$ et $\forall \a\; \F
(V_\a) =\bigcup_K \F (V_\a)^K$. . Donc
$$
\Hom (aC^\bullet, \F) =\underset{K}{\limi} \Hom (aC^\bullet, \F)^K
$$
La limite inductive étant filtrante il suffit de montrer que $\forall
K$ suffisamment petit
$\Hom (aC^\bullet, \F)^K$ est acyclique en degré $>0$.

Soit donc $K$ petit. On a un couple de foncteurs adjoints comme dans
la section \ref{kfjduytb35sfo}
$$
(j_{K!},j^*_K) : \coprod_{i\in I} ( \CC/U_i)\top_{K-\lss} \ldrt \CC\top_{G-\lss}
$$
d'où une résolution dans $\La-\CC\top_{G-\lss}$
$$
N^\bullet \ldrt \underline{\La}
$$
$$
N^{-p+1} =\underbrace{ ( j_{K!} j^*_K)\circ \dots \circ (j_{K!}
  j^*_K)}_{p-\text{ fois}} 
=\bigoplus_{\{i_1,\dots, i_p \}\subset I}
\text{c-ind}^{G/X}_{K/ {h_{U_{i_1}}\times_{h_X} \dots \times_{h_X}
    h_{U_{i_p}}}} \underline{\La}
$$
et $\Hom_{\La-\CC\top} (aC^\bullet, \F)^K = \Hom_{\La-\CC\top_{G-\lss}}
(N^\bullet, \F)$. D'où le résultat.
\qed

\subsection{Cohomologie des $G$-faisceaux lisses}

\begin{theo}
Supposons les hypothèses de la section précédente vérifiées. 
Soit $\F\in \La-\CC\top_{G-\lss}$ un $G$-faisceau de $\La$-modules
lisse. Alors $\forall U\in \CC\; \forall q\geq 0$ le $\La[G]$-module
$H^q (U,\F)$ est lisse. De plus le foncteur 
 \begin{eqnarray*}
 \La-\CC\top_{G-\lss} & \ldrt & G-\DD^+ ( \La) \\
\F & \longmapsto & R\GG ( U,\F)
\end{eqnarray*}
se factorise de façon canonique par $\DD^+ ( \La[G]_{\lss})$ en un
foncteur $\F\longmapsto R\GG_{G} (U,\F)$. 
\end{theo}
\dem
Il suffit d'appliquer le théorème précédent et d'écrire le foncteur 
dérivé comme foncteur composé. Plus précisément, le foncteur 
\begin{eqnarray*}
F:\La-\CC\top_{G-\lss} & \ldrt & \La[G]-\text{Mod}_{\lss} \\
 \F & \longmapsto & \GG (U,\F)
\end{eqnarray*}
se factorise en $F=G\circ k$ où $k$ est le foncteur exact de
l'inclusion de $\La-\CC\top_{G-\lss}$ dans $\La-\CC\top$ et $G=\GG
(U,-)$.  D'après le théorème précédent $RF = RG\circ k$. 
\qed

\subsection{Faisceaux lisses sur une tour formée par un
  pro-torseur}\label{ldddhzZRHREY457H}

Soit $\CC$ un site.
\\
Soit $G$ un groupe fini et $Y\xrig{\; p\; }X$ un $G$-torseur au dessus
de $X$ dans $\CC$. Rappelons qu'il y a une équivalence
\begin{eqnarray*}
\text{Faisceaux sur } X & \xrig{\; \sim \;} & \text{Faisceaux
}G-\text{équivariants sur }Y \\
\F & \longmapsto & p^* \F
\end{eqnarray*}
où $G$-équivariant signifie une action compatible à celle de $G$ sur $Y$
(ce que les adeptes des champs notent $[G\bc Y] =X$). L'équivalence
inverse est donnée par 
$$
\G\longmapsto (p_*\G)^G
$$
Le but de cette section est de généraliser cela aux pro-torseur sous
un groupe profini.

\begin{exem}
Avant de nous lancer expliquons de que l'on veut faire dans un cas
particulièrement simple. Soit $k$ un corps de clôture séparable
$\overline{k}$ et $X$ un $k$-schéma de type fini. Le topos étale de
$X_{\overline{k}}$ s'identifie à la limite projective des topos étales
$X_L$ où $L|k$ parcourt les extensions galoisiennes de degré fini. En
effet, le site étale de $X_{\overline{k}}$ s'identifie à la limite
inductive des sites $X_L$, $L|k$ galois finie, au sens où les
$X_{\overline{k}}$-schémas étales (quasicompacts) sont les germes de $X_L$-schémas
étales lorsque $L$ varie et un morphisme de germe est couvrant en un
niveau suffisamment grand ssi il le devient sur
$X_{\overline{k}}$. Les faisceaux étales sur $X$
s'identifient alors via $p^*$ où $p:X_{\overline{k}} \ldrt X$ aux
faisceaux étales sur $X_{\overline{k}}$ munis d'une action de $\Gal
(\overline{k}|k)$ compatible à celle de $\Gal (\overline{k}|k)$ sur
$X_{\overline{k}}$ qui est discrète au sens où les sections sur un schéma
étale quasicompact ont un stabilisateur ouvert. La condition est une
condition de continuité sur la donnée de descente relativement au
pro-$\Gal (\overline{k}|k)$-torseur $X_{\overline{k}}\ldrt X$.
\end{exem}

Soit donc $G$ un groupe profini et $(G_i)_{i\geq 0}$ une famille
décroissante de sous-groupes ouverts distingués dans $G$ telle que 
$\bigcap_{i\geq 0} G_i =\{e\}$ et $G_0 =G$. Supposons nous donné une
tour dans $\CC$
$$
\dots \ldrt X_{i+1} \ldrt X_i \ldrt \dots \ldrt X_0
$$
où l'on notera $\forall i\geq j\; \ph_{ij}: X_i\ldrt X_j$. Supposons
que cette tour est munie d'une action de $G$ au sens où $\forall i\;
G/G_i$ agit sur $X_i$ et $\forall i\geq j$ le morphisme $X_i\ldrt X_j$
est compatible à ces actions via $G/G_i\twoheadrightarrow G/G_j$.

Supposons de plus que $\forall i$ l'action de $G/G_i$ sur $X_i$ fasse
de $X_i$ un $G/G_i$-torseur au dessus de $X_0$ (et donc $\forall i\geq
j\; X_i\ldrt X_j$ est un $G_i/G_j$-torseur). 

Considérons le système projectif de topos 
$$
\dots \ldrt \CC\top/X_{i+1} \ldrt \CC\top/X_{i} \ldrt \dots \ldrt \CC\top/X_{0}
$$
où rappelons que les topos $\CC\top/U_i$ s'identifient aux $(\CC/U_i)\top$. 
La limite projective de ce système $\mathcal{T}_\infty=\underset{i\geq 0}{\limp}
(\CC\top/X_i)$ est la catégories des système $(\F_i)_{i\geq 0}$, $\F_i
\in \CC\top/X_i$, munis d'isomorphisme $c_{ij}:\forall i\geq j\;
\ph_{ij*}\F_i\iso \F_j$ vérifiant une condition de cocyle. 
\\

Il y a une ``action'' de $G$ sur $\mathcal{T}_\infty$ au sens où pour tout tout $g\in G$ il y a une équivalence de topos $\a_g :\mathcal{T}_\infty \ldrt \mathcal{T}_\infty$
définie par $\a_g = (g^*,g_*)$ où $g^* ((\F_{i})_i, (c_{ij})_{i\geq j}) = ((g^* \F_i)_i, (g^*(c_{ij}))_{i\geq j})$ vérifiant $\forall g,g'\in G$ il a une transformation naturelle
$\a_{gg'}\iso\a_g\circ \a_{g'}$ satisfaisant la condition une condition de cocyle.
On peut donc définir 
$$
G- \mathcal{T}_\infty
$$
la catégorie des $G$-objets dans $\mathcal{T}_\infty$ comme étant les $A\in \mathcal{T}_\infty$ munis d'isomorphismes $\forall g\in G\; g^* A\iso A$ satisfaisant une condition de cocyle.
On vérifie qu'en fait
$$
G-\mathcal{T}_\infty \iso \underset{i\geq 0}{\limp} (\CC/X_i)\top_G
$$ 
On remarquera qu'il faudra faire attention à ne pas prendre $(\CC/X_i)\top_{G/G_i}$ mais bien 
$(\CC/X_i)\top_G$ et que $(\CC/X_i)\top_G$ est par définition la
catégorie des faisceaux $G$-équivariants au dessus de $X_i$ au sens où
l'action est compatible l'action de $G$ sur
$X_i$ via $G\twoheadrightarrow G/G_i$; en particulier un $\F\in
(\CC/X_i)\top_G$ est muni d'une action ``abstraite'' de $G_i$.
\\

Pour tout $i\geq 0$ l'objet $X_i$ muni de son action de $G$ via
$G\twoheadrightarrow G/G_i$ vérifie les hypothèses de la section
\ref{lfizyr25si} puisque $G_i$ agit trivialement sur $X_i$. On peut donc
définir $(\CC/X_i)\top_{G-\lss}$ la catégorie des $G$-faisceaux
lisses. Un $\F\in (\CC/X_i)\top_{G-\lss}$ est lisse ssi le faisceau
muni d'une action de $G_i$ $\text{Res}^{G/X_i}_{G_i/X_i} \F$ l'est.

\begin{lemm}
Soit $(\F_i)_{i\geq 0}\in G-\mathcal{T}_\infty$ où $\forall i\; \F_i
\in (\CC/X_i)\top_G$. Sont équivalents
\begin{itemize}
\item $\forall i\; \F_i$ est lisse
\item $\F_0$ est lisse
\end{itemize}
\end{lemm}
\dem
Supposons $\F_0$ lisse. On a un isomorphisme $\forall i\; \F_0\iso
\ph_{i0*} \F_i$. Donc, si $K\subset G_i$, $K$ agissant trivialement
sur $U_i$
$$
\F_0^K  \iso \ph_{i0*} (\F_i^K)
$$
(puisque $\ph_{i0*}$ est le foncteur image directe d'un morphisme de
topos il commute aux limites projectives quelconques et donc aux
invariants sous $K$). Le faisceau $\F_i$ est lisse ssi
$$
\a: \underset{K\subset G_i}{\limi} \F_i^K \hookrightarrow \F_i
$$
est un isomorphisme. Il y a un diagramme
$$
\xymatrix{
\underset{K}{\limi} \F_0^K \ar[r]^\simeq \ar[rd]_\simeq & \underset{K}{\limi}
\ph_{i0*} (\F_i^K) \ar[r]^\simeq & \ph_{i0*} ( \underset{K}{\limi}
\F_i^K) \ar[d]^{\ph_{i0*}(\a)} \\
 & \F_0 \ar[r]^\simeq & \ph_{i0*} \F_i
}
$$
(où $\ph_{i0*}$ commute aux limites inductives car $X_i\ldrt X_0$
étant un torseur sous un groupe fini $\ph_{i0*} ``=\ph_{i0!} \, ``$ est
adjoint à gauche de $\ph_{i0}^* ``=\ph_{i0}^{!} ``$). Donc $\ph_{i0*}
(\a)$ est un isomorphisme ce qui implique que $\a$ en est un puisque
$\ph_{i0}$ est couvrant.
\qed

Par définition on notera $(G-\mathcal{T}_\infty)_{\lss}$ les systèmes
projectifs de faisceaux satisfaisant les hypothèses du lemme précédent.

\begin{theo}\label{kegjuegt2T35tp}
Supposons que tout objet de $\CC$ est quasicompact. Il y a alors une
équivalence de catégories entre $(G-\mathcal{T}_\infty)_{\lss}$, les
faisceaux en niveau infini munis d'une action lisse de $G$ compatible
à l'action sur la tour,
 et $(\CC/X_0)\top$, les faisceaux sur la base de la tour.
\end{theo}
\dem
Soit $(\F_i)_{i\geq 0}$ un objet de $(G-\mathcal{T}_\infty)_{\lss}$ où
$\forall i\; \F_i$ est muni d'une action de $G$ compatible à celle de
$G/G_i$ sur $X_i$. Associons lui $\F_0^G\in (\CC/X_0)\top$ le faisceau
des invariants de $\F_0$ sous $G$.

Réciproquement, étant donné $\G\in (\CC/X_0)\top$ posons 
$$
\forall i\geq 0\; \F_i= \underset{j\geq i}{\limi} \ph_{ji*} \ph_{j0}^* \G
$$
\'Etant donné que $\ph_{j0}^* \G$ est muni d'une action de $G/G_j$
compatible à celle de sur $X_j$ si $j\geq i$ le faisceau
$\ph_{ji*}\ph_{ji}^* \G$ est muni d'une action de $G/G_j$ compatible à
l'action de $G/G_j$ sur $X_i$ via $G/G_j\twoheadrightarrow G/G_i$. De
plus $\F_i$ est lisses comme limite inductive de faisceaux lisses
(lemme \ref{kfjsfiz24694zf}). 

Reste à voir que cela définit bien des équivalences inverses. 
\\
Si $(\F_i)_i \in (G-\mathcal{T}_\infty)\top_{\lss}$, $\F_i$ étant
lisse
$$
\F_i=\underset{j\geq i}{\limi} \F_i^{G_j}
$$
L'isomorphisme $\F_i\iso \ph_{ji*}\F_j$ induit $\F_i^{G_j} \iso
\ph_{ji*} (\F_j^{G_j})$. Le faisceau $\F_j^{G_j}$ sur $X_j$ est muni
d'une action de $G/G_j$ compatible à celle sur $X_j$. D'après les
rappels du début de cette section, puisque $X_j\ldrt X_0$ est un
$G/G_j$-torseur, cela implique
\begin{eqnarray*}
\F_j^{G_j} \simeq \ph_{j0}^*\left [ \ph_{j0*} (\F_j^{G_j}) \right
]^{G/G_j} &=& \ph_{j0}^* \left ( (\ph_{j0*} \F_j)^{G} \right ) \\
&\simeq & \ph_{j0}^* ( \F_0^G)
\end{eqnarray*}
Donc au final 
$$
\F_i \simeq \underset{j\geq i}{\limi} \ph_{ji*}\ph_{j0}^*  (\F_0^G)
$$
Dans l'autre sens si $\G\in (\CC/X_0)\top$, soit 
$$
\F_0  = \underset{i\geq 0}{\limi} \ph_{i0*} \ph_{i0}^* \G 
$$
Utilisant l'hypothèse de quasicompacité on vérifie que le préfaisceau
limite inductive précédent est séparé. On en déduit alors en utilisant
la définition du faisceau associé à un préfaisceau séparé couplée à
l'hypothèse de quasicompacité que 
\begin{eqnarray*}
\F_0^G = \underset{i\geq 0}{\limi}  \left [ (\ph_{i0*} \ph_{i0}^*
  \G)^G \right ] &=& \underset{i\geq 0}{\limi} \left [ (\ph_{i0*} \ph_{i0}^*
  \G)^{G/G_i} \right ] \\
&=& \underset{i\geq 0}{\limi} \G =\G
\end{eqnarray*}
\qed

\section{Cohomologie à support compacte équivariante-lisse des
  espaces analytiques de Berkovich}\label{kfuhDGT23576gffgf}

Soit $X$ un $K$-espace analytique de Berkovich paracompact.
L'espace $X$ vérifie cette hypothèse si par exemple $X=\X^{an}$ où
$\X$ est un schéma formel localement formellement de type fini sur
$\spf (\O_K)$. On fixe un anneau $\La$ qui servira comme coefficients.

\subsection{Les quatre suites spectrales de cohomologie de Cech permettant de
  calculer la cohomologie à support compact}

On note $X_\et$ le site étale de $X$ tel que défini dans
\cite{Berk1}. 

\begin{theo}\label{46sfjuf48yefg}
Soit $\F$ un faisceau de $\La$-modules sur $X_\et$. 
\begin{itemize}
\item
Soit $(U_i)_{i\in I}$ un recouvrement ouvert de $X$ tel que $\forall i\;
\{j\;|\; U_j\cap U_i\neq \emptyset\}$ soit fini. Il y a alors
une suite spectrale concentrée dans le cadrant $(p\leq 0, q\geq 0)$
$$
E^{pq}_1 = \bigoplus_{\a \subset I\atop |\a|=-p+1} H^q_c (
U(\a),\F)\limpl H^{p+q}_c (X,\F)
$$
où $U(\a) = \bigcap_{i\in \a} U_i$. 
\item Soit $(U_i)_{i\in I}$ un recouvrement ouvert localement fini de
  $X$ tel que $\forall i\; \overline{U_i}$ soit compact. Il y a alors
  une suite spectrale concentrée dans le cadran $(p\geq 0,q\geq 0)$
$$
E^{pq}_1 = \bigoplus_{\a \subset I\atop |\a|=p+1} H^q (U(\a),\F)
\limpl H^{p+q}_c ( X,\F)
$$
où $U(\a) = \bigcap_{i\in \a } U_i$.
\item Soit $(W_i)_{i\in I}$ un recouvrement localement fini de $X$
  formé de domaines analytiques fermés compacts. Il y a alors une
  suite spectrale concentrée dans le cadran $(p\geq0,q\geq 0)$ 
$$
E^{pq}_1 = \bigoplus_{\a \subset I\atop |\a|=p+1} H^q (W(\a),\F)
\limpl H^{p+q}_c ( X,\F)
$$
où $W(\a) = \bigcap_{i\in \a} W_i$.
\item Sous les mêmes hypothèses que précédemment, supposant de plus
  que les intérieurs $\overset{\circ}{W_i}$ recouvrent $X$, il y a une suite
  spectrale concentrée dans le cadran $(p\leq 0,q\geq 0)$
$$
E^{pq}_1 = \bigoplus_{\a \subset I\atop |\a|=-p+1} H^q_{W(\a)} (X,\F)
\limpl  H^{p+q}_c ( X,\F)
$$
\end{itemize}
\end{theo} 
\dem
On considère les deux morphismes surjectifs
$$
j:\coprod_{i\in I} U_i \ldrt X\;\;\; \; i: \coprod_{i\in I} W_i \ldrt X
$$
auxquels sont associés quatre couples de foncteurs adjoints $(j_!,j^*)$
pour la première suite spectrale, $(j^*,j_*)$ pour la seconde,
$(i^*,i_*)$ pour la troisième et $(i_*,i^!)$ pour la dernière.  Soit
$\G$ un faisceau de $\La$-modules sur $X_\et$. 
\`A chacun de ces couples de foncteurs adjoints est associé une
résolution simpliciale, $C^\bullet \ldrt \F$, ou co-simpliciale,
$\G\ldrt C^\bullet$ où si $j_\a : U(\a)\ldrt X$, $i_\a: W(\a)\ldrt X$,
$$
C^p = \bigoplus_{\a\subset I\atop |\a|=-p+1} j_{\a!}j_\a^* \G\;\;\;\;\;\;
\GG_c (X,C^p) = \bigoplus_{\a\subset I\atop |\a|=-p+1} \GG_c ( U(\a),\G)
$$
pour la première suite spectrale 
$$
C^p=\prod_{\a\subset I\atop |\a| = p+1} j_{\a*}j_\a^* \G\;\;\;\;\;\; \GG_c (
X,C^p) =  \bigoplus_{\a\subset I\atop |\a|=p+1} \GG ( U(\a),\G)
$$
pour la seconde (utiliser l'hypothèse que les $\overline{U}_i$ sont
compacts et les $U_i$ localement finis pour voir que $\GG_c (X,C^p)
\subset \GG (X,C^p) =\prod_{|\a|=p+1} \GG (U(\a),\G)$ est le
sous-module formé par la somme directe),
$$
C^p =\prod_{\a\subset I\atop |\a|=p+1} i_{\a*}i_\a^* \G\;\;\;\;\;\; \GG_c
(X,\G) = \bigoplus_{\a\subset I\atop |\a|=p+1} \GG ( W(\a),\G)
$$
pour la troisième,
$$
C^p = \prod_{\a\subset I\atop |\a|= -p+1} i_{\a*} i_\a^! \G\;\;\;\;\;\;
\GG_c ( X,C^p) = \bigoplus_{\a\subset I\atop |\a|=-p+1}  \GG_{W(\a)} ( X,\G)
$$
pour la dernière (pour vérifier que cette dernière est une résolution
il faut utiliser que lorsqu'on lui applique $i^!$ elle devient une
résolution et que puisque les intérieurs des $W_i$ recouvrent $X$
cela implique que c'est une résolution).

Soit $\pi : X_\et\top \ldrt |X|\top$ la projection du topos étale sur
celui de l'espace paracompact $|X|$.  
Rappelons qu'un faisceau $\mathcal{H}$ sur $X_\et$ est dit mou si sa projection
$\pi_* \mathcal{H}$ est un faisceau mou au sens de Godement et si $\forall x\in
X$ le module galoisien lisse $\mathcal{H}_x = i_x^*\mathcal{H}$, $i_x:\mathcal{M}
(\mathcal{K}(x))\ldrt X$, est flasque. Si $\mathcal{H}$ est mou alors $\forall
i>0\; R^i\pi_* \mathcal{H}=0$ grâce à l'hypothèse sur les fibres, et $\pi_* \mathcal{H}$
est $\GG_c(|X|,-)$-acyclique. Les faisceaux mous sur $X_\et$
sont donc $\GG_c(X,-)$-acycliques. Tout faisceau injectif est mou. 
On notera également que la notion d'être mou est une notion locale sur
$|X|$. 

Soit $\G$  un faisceau de $\La$-modules injectif.
\begin{itemize}
\item Pour la première suite spectrale $\dpt{\bigoplus_{\a\subset I\atop
    |\a|=-p+1} j_{\a!}j_\a^* \G}$ est un faisceau mou car une somme
  directe finie de faisceaux mous est molle, la notion de mollesse est
  locale sur $|X|$, les $(U_i)_i$ sont localement finis donc la somme
  directe précédente est localement une somme finie 
et enfin pour
  une immersion ouverte $\ph: V\hookrightarrow X$ $\;\ph^*$ et $\ph_!$
  transforment faisceaux mous en faisceaux mous. On déduit également
  de ces arguments que $\forall \a\; j_\a^*\G$ étant mou est $\GG_c
  (U(\a),-)$-acyclique. 
\item Pour la seconde suite spectrale $\dpt{\prod_{\a\subset I\atop
    |\a|=p+1} j_{\a*}j_\a^* \G}$ est un faisceau injectif comme
produit de faisceaux injectifs (si $\ph :V\hookrightarrow X$ est une
immersion ouverte alors $\ph^*$, resp. $\ph_*$, a comme adjoint à
gauche
$\ph_!$, resp. $\ph^*$, donc $\ph^*$ et $\ph_*$ envoient les injectifs
sur des injectifs). De même $\forall \a\; j_a^*\G$ est injectif. 
\item Pour la troisième $\dpt{C^p =\prod_{\a\subset I\atop |\a|=p+1}
    i_{\a*}i_\a^* \G}$ est mou car la notion de mollesse est locale sur
  $|X|$, 
localement sur $|X|$ le produit
  précédent est fini, si $\ph : Z\hookrightarrow X$ est l'immersion
  définie par un domaine analytique compact alors $\ph^*$ et
  $\ph_*=\ph_!$ conservent la mollesse. De même $\forall \a\; i_\a^*\G$
  étant mou et $W(\a)$ compact il est $\GG (W(\a),-)$-acyclique.
\item Enfin pour la dernière suite spectrale $\dpt{C^p =
    \prod_{\a\subset I\atop |\a|= -p+1} i_{\a*} i_\a^! \G}$ est
  injectif comme produit d'injectifs puisque si $\ph :
  Z\hookrightarrow X$ alors $\ph_*$ et $\ph^!$, qui possèdent tous deux
  des adjoints à gauche, envoient injectifs sur injectifs. De même
  $\forall \a \; i_\a^!\G$ est injectif.
\end{itemize}

Les résolutions $C^\bullet (\G)$ sont fonctorielles en $\G$. Si
$\F\ldrt \mathcal{I}^\bullet$ est une résolution injective le complexe
total du complexe double $\GG_c(X,C^\bullet (\mathcal{I}^\bullet))$
fournit les suites spectrales voulues.
\qed

\subsection{Le site quasi-étale de $X$ et cohomologie à support compact}

On note $X_\et$ le site étale de $X$, 
$X_\qet$ son site quasi-étale  (\cite{Berk2}) et enfin $X_\qetc$
le site quasi-étale compact qui est
le ``sous-site'' de $X_\qet$ formé des morphismes quasi-étales $U\ldrt X$ 
avec $U$ compact.

Il y a un morphisme de sites $X_\qet\ldrt X_\et$ qui induit un
morphisme de topos (\cite{Berk2} section 3) 
$$
(\mu^*,\mu_*) : X_\qet\top \ldrt X_{\et}\top
$$
Il y a également un morphisme de sites $X_\qet\ldrt X_\qetc$ qui
induit une équivalence de topos 
$$
X_\qet\top\iso X_\qetc\top
$$
On en déduit un morphisme de topos
$$
(\nu^*,\nu_{*}) : X_\qetc\top \ldrt X_\et\top
$$
auquel s'appliquent tous les résultats de la section 3 de \cite{Berk2} 
concernant $(\mu^*,\mu_*)$. 

Concrètement si $\F$ est un faisceau sur
$X_\et$ et $f:W\ldrt X$  quasi-étale avec $W$ compact on a $\GG(W,\nu^*\F)
 = \GG (W,f^*\F)$. Si $\G$ est un faisceau sur $X_\qetc$ et
$f: U\ldrt X$ est étale alors $\GG (U,\nu_* \G) = \underset{W\subset
  U}{\limp} \GG (W,\F)$ où $W$ parcourt les domaines analytiques
compacts dans $U$. 

\begin{rema}
Si $X^{rig}$ est l'espace rigide associé à $X$ (\cite{Berk1} section 1.6)
alors $X_\qet\top$ s'identifie au topos rig-étale de $X^{rig}$. Le
foncteur $\mu^*$ est pleinement fidèle d'image essentielle les
faisceaux rig-étale surconvergents. De même si $X^{ad}$ désigne
l'espace adique sur $\spa (\O_K,K)$ associé à $X$ (\cite{Hu1} section
8.3) alors $X_\qet\top$
s'identifie au topos étale de $X^{ad}$ et $X_\et\top$ au topos étale
partiellement propre de $X^{ad}$. Dans le cadre des espaces adiques le
morphisme de topos $(\mu^*,\mu_*)$ est alors celui noté $\theta_X$
dans  
le chapitre 8 de \cite{Hu1}. 
\end{rema}

\begin{prop}\label{kfjsuf546EFfz}
Soit $\F$ un faisceau de $\La$-modules sur $X_\et$ tel que
$\nu^*\F$ soit flasque sur $X_\qetc$. Alors $\F$ est $R\GG_c
(X,-)$-acyclique. 
\end{prop}
\dem
Appliquons la troisième suite spectrale du théorème
\ref{46sfjuf48yefg} dont on reprend les notations.
L'espace $X$ étant paracompact 
il existe un recouvrement $(W_i)_{i\in I}$ tel que dans les hypothèses
de ce théorème. De plus $\forall \a\subset I$ de cardinal fini 
$$\forall q>0\; H^q (W(\a),\F) = H^q ( W(\a),\nu^*\F)$$
 (\cite{Berk2}
théorème 3.3) qui est nul puisque $W(\a)$ est dans $X_\qetc$. 
\qed

\begin{prop}
Soit $\F$ un faisceau flasque sur $X_{\qetc}$. Alors $\F$ est
$\nu_{*}$-acyclique. 
\end{prop}
\dem
Pour $q>0$ le faisceau $R^q\nu_{*} \F$ est le faisceau associé au
préfaisceau qui  à $f:U\ldrt X$ étale associe $H^q
(U_{\qetc},f^*\F)$. Mais tout point d'un tel $U$ possède un voisinage
$W\subset U$ où $W\ldrt X$ est dans $X_{\qetc}$. Alors $H^q (
W_{\qetc},(f^*\F)_{|W})= H^q (W,\F)=0$.
\qed 

\begin{prop}\label{jfysfyztTGG4}
Soit $\F$ un faisceau flasque sur $X_\qetc$. Alors $\nu_{*}\F$ est
$R\GG_c (X,-)$-acyclique et $\nu^*\nu_* \F$ est flasque sur
$X_{\qetc}$.  
\end{prop}
\dem
D'après la proposition \ref{kfhztr26Ie} il suffit de montrer que
$\nu^*\nu_{*}\F$ est flasque. Mais si $f:W\ldrt X$ est un morphisme
dans $X_{\qetc}$ alors $H^q ( W,\nu^* \nu_{*} \F) = H^q ( W, \nu_{*} ( f^* \F)) = H^q (W,\F)$ d'après la
proposition précédente.
\qed 

\begin{coro}\label{kfyvyzrput24zssf9}
Soit $\F$ un faisceau de groupes abéliens sur $X_\et$ tel que
$\nu^*\F$ soit flasque sur $X_\qetc$. Alors $\F$ est mou.
\end{coro}
\dem
Soit $x\in X$ et $i_x:\mathcal{M} ( \mathcal{K} (x))\ldrt X$. Alors on
a une égalité de 
modules galoisiens $i_x^* \F =i_x^*\nu^* \F$ qui est donc flasque. Soit
maintenant $U$ un ouvert de $X$ et $K\subset |U|$ un compact. On doit
montrer que l'application $\GG (U,\F)\ldrt \GG (K,\F)$ est surjective
où $K$ est vu comme germe d'espace analytique dans $X$ au sens de la
section 3.4 de \cite{Berk2}. Mais il y a une suite exacte 
$$
0\ldrt \GG_c (U\setminus K,F)\ldrt \GG_c (U,\F) \ldrt \GG (K,\F) \ldrt
H^1_c ( U\setminus K,\F)
$$
Le terme de droite est nul d'après la proposition précédente.
\qed

\subsection{Faisceaux équivariants lisses}

On applique maintenant le formalisme de la section
\ref{lsfkjsut25tggpoi}.
Soit $G$ un groupe topologique possédant un sous-groupe ouvert profini
et agissant continûment sur $X$ au sens de la section 6 de
\cite{Berk2}. D'après le théorème 7.1 de \cite{Berk2} le site
$X_\qetc$ satisfait aux hypothèses de la section
\ref{lsfkjsut25tggpoi}. On note $X_{\qetc-G-\lss}\top$ le topos
des faisceaux $G$-équivariants lisses sur $X_\qetc$.

\begin{defi}
On note $X\top_{\et-G-\lss}$ la catégorie des $G$-faisceaux
 sur $X_\et$ tels que $\nu^*\F$ soit lisse. 
\end{defi}

En d'autres termes un $G$-faisceau étale $\F$ est lisse ssi $\forall
f:W\ldrt X$ quasi-étale avec $W$ compact et $K$ un sous-groupe 
  compact ouvert de $G$ dont l'action se relève à $W$ le $K$-ensemble $\GG (W,f^* \F)$ est 
lisse.

\begin{lemm}\label{jkfytghputzr35egv}
La catégorie $\La-X_{\et-G-\lss}\top$ est abélienne et possède
suffisamment d'injectifs. Elle s'identifie à l'image essentielle par
$\nu_{*}$ de $\La-X\top_{\qetc-G-\lss}$.
\end{lemm}
\dem
Le fait que la catégorie soit abélienne résulte de l'exactitude de
$\nu^*$ et de ce que $\La-X\top_{\qetc-G-\lss}$ est abélienne
(corollaire \ref{kfuyzr35fzfr46}). L'assertion concernant l'image
essentielle résulte de ce que l'adjonction $Id \ldrt \nu_{*}\nu^*$
est un isomorphisme (\cite{Berk2} corollaire 3.5). L'assertion
concernant les injectifs résulte de ce que $\La-X_{\qetc-G-\lss}$
possède suffisamment d'injectifs (lemme \ref{kfusfufuzr45r}),  
 si $\F\in
\La-X_{\qetc-G-\lss}$  et
$u:\nu^*\F\hookrightarrow \mathcal{I}$ est un plongement de
$\nu^*\F$ dans un injectif alors $\F\iso \nu_{*}\nu^* \F 
\ldrt \nu_{*} \mathcal{I}$ est un morphisme vers un objet
injectif. Ce morphisme est un monomorphisme car si $\mathcal{H} = \ker
\nu_{*}u$ alors $\nu^*$ étant exact il y une suite exacte
$$
0\ldrt \nu^*\mathcal{H}\ldrt \nu^*\F \ldrt \nu^*\nu_{*}\mathcal{I} \ldrt
\mathcal{I} 
$$
et donc $\nu^* \mathcal{H}=0$ ce qui implique $\mathcal{H}=0$.
\qed

Comme corollaire de la démonstration précédente on a. 

\begin{coro} \label{kfhztr26Ie}
Si $\F\in \La-X\top_{\et-G-\lss}$ et
$\nu^* \F \ldrt \mathcal{I}^\bullet$ est une résolution injective
dans $\La-X\top_{\qetc-G-\lss}$ alors 
$$
\F\ldrt \nu_{*}\mathcal{I}^\bullet
$$
en est une dans $\La-X\top_{\et-G-\lss}$ et de plus les faisceaux $
\nu_{*}\mathcal{I}^\bullet$ sont mous.
\end{coro}

\begin{prop}
La catégorie
$X_{\et-G-\lss}$ est un topos. 
\end{prop}
\dem L'assertion concernant l'existence de limites projectives finies, de limites inductives quelconques et de ``bons'' quotients par des relations d'équivalence résulte du cas de $X_{\qetc-G-\lss}$ puisque $\nu^*$ commute à ces opérations.
Reste à exhiber un petit système de générateurs. Pour cela remarquons que les ouverts étales distingués engendrent la topologie de $X_\et$. Maintenant si $f: U\ldrt X$ est un morphisme étale avec $U$ distingué, c'est à dire il existe une factorisation $U\hookrightarrow W \ldrt X$ avec $W\ldrt X$ quasi-étale, $W$ compact et $W\setminus U$ un domaine analytique compact alors pour un sous-groupe compact-ouvert $K$ suffisamment petit dans $G$ l'action de $K$ se prolonge à $U$ et $\F(U) = \underset{K}{\limi} \F (U)^K$. 
 Si $h_U$ désigne le faisceau représenté par $U$ sur $X_\et$ on peut alors former comme dans la section \ref{kfjduytb35sfo} le faisceau induit
$$
\text{c-ind}^{G/X}_{K/U} h_U = \bigoplus_{K\bc G} f_! h_U
$$
où $f_!$ désigne les images à support propre au sens de \cite{Berk1}. Il est alors aisé de voir que cela forme un système de générateurs.
\qed

\begin{rema}\label{dsgjdsguuet25YETG}
Les limites projectives quelconques, non-finies, dans
$X\top_{\et-H-\lss}$ se calculent de la façon suivante : si
$(\F_i)_{i\in I}$ est un système projectif de $G$-faisceau étales
lisses alors sa limite projective est
$$
\nu_* \left [ (\underset{i\in I}{\limp} \nu^*\F_i)^{\lss} \right ]
$$
où la limite projective dans l'expression précédente est prise dans la
catégorie des faisceaux sur le site quasi-étale. 
\end{rema}

\subsection{Le complexe de cohomologie lisse}

Reprenons les hypothèses de la section précédente.

\begin{theo}\label{PPiuzgtgf26Usge}
Soit $\F$ un $G$-faisceau de $\La$-modules lisse sur $X_\et$. Alors
$\forall q\geq 0\; H^q_c ( X,\F)$ est un $G$-module lisse. De plus 
le foncteur 
\begin{eqnarray*}
 \La-X\top_{\et-G-\lss} & \ldrt & G-\DD^+ ( \La) \\
\F & \ldrt & R\GG_c (X,\F)
\end{eqnarray*}
où $G-\DD^+ (\La)$ est la catégorie formée des objets de $\DD^+ (\La)$
munis d'une action de $G$, 
se factorise de façon canonique en un foncteur 
$$
R\GG_{c-G} ( X, -) :
\La-X\top_{\et-G-\lss} \ldrt \DD^+ ( \La[G]_{\lss})
$$
où $\DD^+ (\La[G]_{\lss})$ désigne la catégorie dérivée des
$\La$-modules munis d'une action lisse de $G$.  
Ce foncteur est le foncteur dérivé de $\GG_{c} (X,-) :
\La-X\top_{\et-G-\lss} \ldrt \La[G]_{\lss}$. 
\end{theo}
\dem
Comme dans la troisième suite spectrale du théorème
\ref{46sfjuf48yefg} si $(W_i)_i$ est comme dans l'énoncé de ce même
théorème il y a un plongement
$$
\GG_c ( X,\F) \hookrightarrow \bigoplus_{i\in I}  \GG (W_i,\F) =
\bigoplus_{i\in I} \GG(W_i,\nu^* \F) 
$$
et donc si $\F$ est lisse $\GG_c (X,\F)$ est un $G$-module lisse. 
D'après le corollaire \ref{kfhztr26Ie} il suffit alors de montrer que
si $\mathcal{I}$ est un injectif de $\La-X\top_{\qetc-G-\lss}$ alors
$\nu_{*}\mathcal{I}$  est $\GG_c (X,-)$-acyclique. Mais d'après le
théorème \ref{OOPfyzY136Uetet} le faisceau $\mathcal{I}$ sur $X_\qetc$
est flasque. Il suffit alors d'appliquer la proposition
\ref{jfysfyztTGG4}.
\qed

\subsection{Les quatre résolutions/suites spectrales permettant de
  calculer la cohomologie à support compact équivariantes lisse}

Comme précédemment $X$ est muni d'une action continue de
$G$. Rappelons qu'un ouvert $U$ de $X$ est dit distingué s'il s'écrit
$U=W_1\setminus W_2$ où $W_1,W_2$ sont deux domaines analytiques
compacts dans $X$. Un tel ouvert est stabilisé par un sous-groupe
ouvert de $G$.
\\

Le théorème qui suit fait suite au théorème \ref{46sfjuf48yefg}. 

\begin{theo}
Soit $\F$ un $G$-faisceau lisse de $\La$-modules sur $X_\et$. 
Soit $\nu^*\F\ldrt \mathcal{I}^\bullet$ une résolution injective de
$\F$ dans $\La-X_{\qetc-G-\lss}$ et $\F\ldrt \nu_* \mathcal{I}^\bullet$
la résolution injective associée dans $\La-X_{\et-G-\lss}$ (proposition
\ref{kfhztr26Ie}). Notons $\F\ldrt \mathcal{J}^\bullet$ cette
résolution.
\begin{itemize}
\item Soit $(U_i)_{i\in I}$ un recouvrement ouvert de $X$ tel que $\forall i\;
\{j\;|\; U_j\cap U_i\neq \emptyset\}$ soit fini. 
Supposons ce recouvrement invariant sous $G$ au sens où $\forall g\in
G\; \forall i\in I \; \exists j\in I\; g.U_i=U_j$. 
Il y a alors un
isomorphisme dans $\DD^+ (\La[G]_{\lss})$ 
$$
R\GG_{c-G} ( X,\F) \simeq \text{Tot } D^\bullet (\mathcal{J}^\bullet)
$$
où $\G\longmapsto D^\bullet (\G)$ est le complexe simplicial
fonctoriel en $\G$ 
$$
D^p (\G) = \bigoplus_{\a\subset I\atop |\a|=-p+1} \GG_c ( U(\a), \G)
$$
\`A cette résolution est associé la suite spectrale de
$\La[G]$-modules lisses
$$
E^{pq}_1 = H^q_c ( U(\a),\F)\limpl H^{p+q}_c ( X,\F)
$$
\item Soit $(U_i)_{i\in I}$ un recouvrement invariant sous $G$ localement fini de $X$ par
  des ouverts distingués. 
Il y a alors un
isomorphisme dans $\DD^+ (\La[G]_{\lss})$ 
$$
R\GG_{c-G} ( X,\F) \simeq \text{Tot } D^\bullet (\mathcal{J}^\bullet)
$$
où $\G\longmapsto D^\bullet (\G)$ est le complexe co-simplicial
fonctoriel en $\G$ 
$$
D^p (\G) = \bigoplus_{\a\subset I\atop |\a|=p+1} \GG  ( U(\a), \G)
$$
\`A cette résolution est associé la suite spectrale de
$\La[G]$-modules lisses
$$
E^{pq}_1 = H^q ( U(\a),\F)\limpl H^{p+q}_c ( X,\F)
$$
\item  Soit $(W_i)_{i\in I}$ un recouvrement invariant sous $G$ localement fini de $X$ par
  des domaines analytiques compacts.  
Il y a alors un
isomorphisme dans $\DD^+ (\La[G]_{\lss})$ 
$$
R\GG_{c-G} ( X,\F) \simeq \text{Tot } D^\bullet (\mathcal{J}^\bullet)
$$
où $\G\longmapsto D^\bullet (\G)$ est le complexe co-simplicial
fonctoriel en $\G$ 
$$
D^p (\G) = \bigoplus_{\a\subset I\atop |\a|=p+1} \GG  ( W(\a), \G)
$$
\`A cette résolution est associé la suite spectrale de
$\La[G]$-modules lisses
$$
E^{pq}_1 = H^q ( W(\a),\F)\limpl H^{p+q}_c ( X,\F)
$$
\item Soit $(W_i)_{i\in I}$ un recouvrement invariant sous $G$ localement fini de $X$ par
  des domaines analytiques compacts tel que les
  $(\overset{\circ}{W}_i)_i$ recouvrent $X$.   
Il y a alors un
isomorphisme dans $\DD^+ (\La[G]_{\lss})$ 
$$
R\GG_{c-G} ( X,\F) \simeq \text{Tot } D^\bullet (\mathcal{J}^\bullet)
$$
où $\G\longmapsto D^\bullet (\G)$ est le complexe simplicial
fonctoriel en $\G$ 
$$
D^p (\G) = \bigoplus_{\a\subset I\atop |\a|=p+1} \GG_{W(\a)}  ( X, \G)
$$
\`A cette résolution est associé la suite spectrale de
$\La[G]$-modules lisses
$$
E^{pq}_1 = H^q_{W(\a)} ( X,\F)\limpl H^{p+q}_c ( X,\F)
$$
\end{itemize}
\end{theo}
\dem
Reprenons la démonstration du théorème 
\ref{46sfjuf48yefg} ainsi que ses notations. 
On utilise le corollaire  $\ref{kfhztr26Ie}$ qui nous dit que nos
faisceaux $\mathcal{J}^\bullet$ sont mous.

Pour la première suite spectrale l'argument donné dans la
démonstration du théorème \ref{46sfjuf48yefg} fonctionne donc encore
puisque comme expliqué si $\G$ est mou  $\dpt{\bigoplus_{\a\subset I\atop |\a=-p+1}
  j_{\a!}j_\a^* \G}$ l'est encore. De plus  si $\G$ est
$G$-équivariant lisse 
ce faisceau $G$-équivariant
est lisse puisque si $f:W\ldrt X$ est quasi-étale avec $W$ compact
alors
 $$
\GG (W, f^*\bigoplus_{\a\subset I\atop |\a=-p+1}
  j_{\a!}j_\a^* \G) =  \bigoplus_{\a\subset I\atop |\a=-p+1} \GG_c (
  f^{-1} (U(\a)), f^*\F))
$$
et d'après le théorème \ref{PPiuzgtgf26Usge} les sections à support
compact sont lisses. 

Pour la seconde suite spectrale par contre en général si $\G$ est mou
$\dpt{\prod_{\a\subset I \atop |\a|=p+1} j_{\a*} j_\a^* \G}$ n'est pas
mou. Cependant il est $\GG_c (X,-)$-acyclique. En effet,  si
$\ph: U\hookrightarrow X$ est une immersion ouverte et $\mathcal{H}$
un faisceau mou sur $U_\et$ alors $\mathcal{H}$ est $\ph_*$-acyclique
puisque $R^q\ph_* \mathcal{H}$ est le faisceau associé au préfaisceau 
qui à $f:V\ldrt X$ étale associe $H^q (f^{-1} (U),f^*\mathcal{H})$ et
$f^*\mathcal{H}$ est mou. Donc si $V$ est un ouvert relativement
compact de $X$, donc ne rencontrant qu'un nombre fini de $U_i$, 
$$
H^q (V, \prod_{\a\subset I \atop |\a|=p+1} j_{\a*} j_\a^* \G) =
\bigoplus_{\a\subset I \atop |\a|=p+1} H^q ( V\bigcap U(\a), \G) = 0
$$
On conclu alors en appliquant la deuxième suite spectrale du théorème
\ref{46sfjuf48yefg} que l'on a un faisceau $\GG_c (X,-)$-acyclique. 
De plus, utilisant l'hypothèse que les ouverts $U_i$ sont distingués
on vérifie aisément que si $\G$ est équivariant lisse alors 
 $\dpt{\prod_{\a\subset I \atop |\a|=p+1} j_{\a*} j_\a^* \G}$ l'est
 également.

La troisième suite spectrale ne pose pas de problème puisque tous les
faisceaux sont mous et lisses.

La quatrième suite spectrale est laissée en exercice au lecteur.
\qed

\section{Cohomologie à support compact équivariante-lisse des espaces
  rigides généralisés}

  Soit $\X$ un schéma formel $\pi$-adique sans
  $\pi$-torsion. Supposons le muni d'une action continue d'un groupe
  topologique $G$.
 La conditions de continuité signifie que pour tout
  ouvert quasicompact $\mathfrak{U}$ de $\X$ tout entier $n$
et tout
$f\in \GG(\mathfrak{U},\O_\X\otimes \O_K/\pi^n \O_K)$ 
 il
  existe un sous-groupe ouvert de $G$ stabilisant $\mathfrak{U}$ et
  laissant la section $f$ invariante : $\forall k\in K\;\;
  k^* f = f$.

\begin{rema}
Une définition naïve de la notion d'action continue consisterait à
dire que pour tout $\mathfrak{U}$ quasicompact et tout $n$ il existe
un sous-groupe compact ouvert $K$ agissant trivialement sur
$\mathfrak{U}\otimes \O_K/\pi^n \O_K$, c'est à dire $G$ agit
continuement pour ``la topologie forte'' alors que la définition
qu'on a donnée est au sens ``faible''.
 Cependant avec cette définition
le groupe topologique $\Gal (\overline{\mathbb{F}}_p |\Fp)$ n'agirait
pas continûment sur $\spec (\overline{\mathbb{F}}_p)$ et plus
généralement si $X$ est un schéma sur $\spec (\Fp)$ le groupe $\Gal
(\overline{\mathbb{F}}_p |\Fp)$  n'agirait pas continûment sur
$X_{\overline{\mathbb{F}}_p}$. Par exemple si $X$ est un espace rigide
sur $\Qp$ $\;\Gal (\Qpb |\Qp)$ n'agirait pas continûment sur
$X\hat{\otimes} \C_p$. Plus généralement si $G$ est un groupe profini
et  $(\X_K)_{K\subset G}$ est une ``tour'' de schémas formels
$\pi$-adiques munie d'une action de $G$, $G$ n'agirait pas continûment
sur $\underset{K}{\limp} \X_K$ en général. 
\end{rema}

Commençons par vérifier que le site $\X_{\E-rig-\et}$ satisfait aux différentes
hypothèses de la section \ref{lsfkjsut25tggpoi}. Tout d'abord d'après le
corollaire \ref{DSSBeztez46I} tout objet de ce site est quasicompact.

\begin{prop}
 le site
$\X_{\E-rig-\et}$ satisfait à l'hypothèse de continuité de la section
\ref{lfizyr25si}.
\end{prop}
\dem
Soit $\YY\ldrt \X$ un morphisme de type $(\E)$. Ce morphisme se
factorisant par un ouvert quasicompact, ouvert stabilisé par un
sous-groupe ouvert de $G$, 
on peut supposer $\X$ quasicompact.

Supposons d'abord $\X$ affine. \'Ecrivons $\X=\underset{i\in I}{\limp}
\X_i$ où les $\X_i$ sont admissibles affines et $\O_{\X_i}\subset
\O_\X$.  
D'après le théorème \ref{DGufbrzgzr3277Jne} il existe $i_0\in I$ ainsi que $\YY'\ldrt
\X_{i_0}$ de type $(\E)$ et un isomorphisme
$$
\YY\simeq \YY'\times_{\X_{i_0}} \X
$$
On applique le théorème de rigidité, la proposition
\ref{DGTeztJTUtuu2258}. Soit $N$ un entier associé à $\YY\ldrt \X$ comme dans cette
proposition. Par hypothèse, $\X_{i_0}$ étant topologiquement de type
fini sur $\spf (\O_K)$, il existe un sous-groupe ouvert $K$ dans $G$
tel que $\forall k\in K$ le diagramme 
$$
\xymatrix{
\X\otimes \O_K/\pi^N \O_K \ar[r]^k \ar[d]^q & \X\otimes\O_K/\pi^N\O_K \ar[d]^q \\
\X_{i_0} \otimes \O_K/\pi^N \O_K \ar[r]^{Id} & \X_{i_0} \otimes
\O_K/\pi^N \O_K
}
$$
commute. Alors $k^{-1} (\YY) = k^{-1} ( q^{-1} (\YY'))$ et 
donc $k^{-1} ( \YY)\otimes \O_K/\pi^N\O_K = q^{-1}( \YY' \otimes
\O_K/\pi^N\O_K )= \YY\otimes\O_K/\pi^N\O_K$, comme $\X$-schémas. 
D'après la proposition \ref{DGTeztJTUtuu2258} on en déduit qu'il
existe un  prolongement $(\beta_k)_{k\in K}$ de l'action de $K$ à $\YY$ 
$$
\xymatrix@R=4mm@C=6mm{
\YY \ar[r]^{\beta_k} \ar[d] & \YY \ar[d] \\
\X \ar[r]^k & \X
}
$$
Nous avons besoin de caractériser ce prolongement de manière unique afin de le rendre ``canonique''.
Pour cela soit $(\mathcal{V}'_\a)_\a$ un recouvrement affine fini de $\YY'$ et $\forall\a\; (f_{\a\b})_{\b}$,
$f_{\a\b}\in \GG(\mathcal{V}'_\a,\O_{\YY'})$, un système fini de générateurs 
topologiques de la $\GG(\X_{i_0},\O_{\X_{i_0}})$-algèbre topologiquement de type fini $\GG(\mathcal{V}'_\a,\O_{\YY'})$. Soit $(\mathcal{V}_\a)_\a= (\mathcal{V}'_\a)_\a\times_{\YY'} \YY$ et
$\forall \a,\b\; f_{\a\b}\otimes 1 \in \GG (\mathcal{V}_\a,\O_\YY)$ les sections obtenues par changement de base, sections qui engendrent topologiquement la $\GG (\X,\O_\X)$-algèbre $\GG ( \mathcal{V}_\a,\O_\YY)$.
\\
Alors d'après l'assertion d'injectivité dans la proposition \ref{DGTeztJTUtuu2258}
$(\b_k)_{k\in K}$ est l'unique relèvement de l'action de $K$ à $\YY$ tel que
$$
\forall \a,\b\;\forall k\in K\;\; \b_k^*(f_{\a\b}\otimes 1) \equiv f_{\a\b}\otimes 1
\text{ mod } \pi^N
$$
Maintenant soit $(\mathcal{W}_\gamma)_\gamma$  un autre recouvrement affine fini de $\YY$ et $\forall \gamma\; (g_{\gamma \delta})_\delta$, $g_{\gamma \delta}\in \GG (\mathcal{W}_\gamma, \O_\YY)$, un système de générateurs topologiques de $\GG (\mathcal{W}_\gamma,\O_\YY)$ comme $\GG (\X,\O_\X)$-algèbre. On vérifie aisément qu'il existe un sous-groupe ouvert $K'\subset K$ tel que 
$$
\forall \gamma,\delta\;\forall k\in K'\;\; \beta_k^* g_{\gamma\delta}\equiv g_{\gamma\delta}\text{ mod }\pi^N
$$

A partir de cette constatation la cas $\X$ général se déduit du cas affine par recollement en
utilisant cette dernière assertion (et toujours quitte à se restreindre à des sous-groupes ouverts plus petits). 
Plus précisément si l'entier $N$ est associé à $\YY\ldrt \X$ comme dans la proposition  \ref{DGTeztJTUtuu2258} soient $\X=\bigcup_\a \mathfrak{U}_\a$,resp. $\YY= \bigcup_{\a,\b} \mathcal{V}_{\a\b}$, deux recouvrements affines finis de $\X$, resp. $\YY$, tels que $\mathcal{V}_{\a\b}\ldrt \mathcal{U}_\a$. Soit $\forall \a,\b\; (f_{\a\b\gamma})_\gamma$, $f_{\a\b\gamma}\in \GG ( \mathcal{V}_{\a\b},\O_\YY)$ un système de générateurs topologiques de la $\GG (\mathfrak{U}_\a,\O_\X)$-algèbre $\GG ( \mathcal{V}_{\a\b},\O_\YY)$. On montre alors qu'il existe un sous-groupe ouvert $K\subset G$ tel qu'il existe un unique relèvement ($\beta_k)_{k\in K}$ de l'action de $K$ sur $\X$ à $\YY$ vérifiant
$$
\forall \a,\b,\gamma\;\forall k\in K\;\; \beta_k^* f_{\a\b\gamma} \equiv f_{\a\b\gamma}\text{ mod }\pi^N
$$

 Le fait que cette action relevée est
``canonique'' au sens de la section \ref{lfizyr25si} se déduit par une
nouvelle application de l'assertion d'injectivité dans la proposition 
 \ref{DGTeztJTUtuu2258}.
\qed

\begin{prop}
Le site $\X_{\E-rig-\et}$ satisfait à l'hypothèse de la section
\ref{KBIGzg3675nt}. \'Etant donnée une famille finie
$(\mathcal{V}_i\ldrt \YY)_{i\in I}$ de morphismes dans la catégorie
$\E_\X$ des morphismes de type $(\E)$ au dessus de $\X$ si
$$
\mathfrak{Z} = (\underset{i\in I,\YY}{\times} \mathcal{V}_i )^{adh}
$$
désigne le produit fibré de cette famille au dessus de $\YY$ il existe
une famille finie $(\mathfrak{U}_\a\ldrt \ZZ)_\a$  de morphismes au
dessus de $\X$ où $\forall \a\; \mathfrak{U}_\a\ldrt \X$ est de type
$(\E)$ et que $\forall \F\in \X\top_{\E-rig-\et}$ l'application
$$
\Hom_{\X_{\E-rig-\et}^{\widehat{\;\;}}} ( \underset{i\in I,
  h_{\YY^{rig}}}{\times} h_{\mathcal{V}_i^{rig}},\F) \ldrt \prod_{\a} \F ( \mathfrak{U}_\a)
$$ 
soit injective.
\end{prop}
\dem
Le morphisme $\YY\ldrt \X$ se factorisant par un ouvert quasicompact
on peut supposer $\X$ quasicompact, ce que nous ferons.

Supposons le résultat démontré lorsque $\X$ est affine. Soit
$\X=\bigcup_k \mathcal{W}_k$ un recouvrement affine. Appliquant le
résultat dans le cas affine on en déduit l'existence de familles 
$\forall k\; (\mathfrak{U}_\a\ldrt \ZZ\times_\X \mathcal{W}_k)_{\a\in
  A_k}$ vérifiant les hypothèses voulues lorsqu'on tire la situation
en arrière à l'ouvert $\mathcal{W}_k$. On vérifie qu'alors
$(\mathfrak{U}_\a\ldrt \ZZ)_{k,\a\in A_k}$ vérifie bien les
conclusions de l'énoncé.

Soit donc maintenant $\X$ affine. \'Ecrivons $\X =\underset{j\in
  J}{\limp} \X_j$ où $\forall j\; \X_j$ est admissible
affine. D'après le théorème \ref{DGufbrzgzr3277Jne} il existe $j_0\in
J$ ainsi qu'une famille $(\mathcal{V}'_i\ldrt \YY')_{i\in I}$ dans
$\E_{\X_{j_0}}$ induisant la famille $(\mathcal{V}_i\ldrt \YY)_{i\in
  I}$ en niveau infini
$$
\left [ (\mathcal{V}'_i\ldrt \YY')_{i\in I} \times_{\X_{j_0}} \X
\right ]^{adh} \simeq (\mathcal{V}_i\ldrt \YY)_{i\in I}
$$
Le produit fibré $\ZZ' = ( \underset{i\in I,\YY'}{\times}
\mathcal{V}'_i)^{adh}$ est un $\X_{j_0}$-schéma formel
rig-étale. D'après le théorème \ref{rhfueo5pjgf} il existe 
un $\X_{j_0}$-morphisme $\mathcal{W}\ldrt \ZZ'$ tel que
$\mathcal{W}\ldrt \X_{j_0}$ soit de type $(\E)$ et $\mathcal{W}^{rig}\ldrt \ZZ'^{rig}$ soit un
morphisme couvrant. 
\\
Alors $(\mathcal{W}\times_{\X_{j_0}} \X)^{adh} \ldrt \ZZ$ satisfait
aux hypothèses puisque $\forall j\geq j_0$ le morphisme 
$(\mathcal{W}\times_{\X_{j_0}}\X_j )^{rig} \ldrt (
\ZZ'\times_{\X_{j^0}} \X_j)^{rig}$ est couvrant et
$\X_{\E-rig-\et}\top \simeq \underset{j\geq j_0}{\limp} (\X_{j})\top_{rig-\et}$.
\qed
 
Soit maintenant le foncteur 
\begin{eqnarray*}
\GG_! ( \X^{rig},-): \La-(\X_{\E-rig-\et})_{G-\lss}\top &\ldrt & \La[G]_\lss \\
\F & \longmapsto & \GG_! (\X^{rig},\F)
\end{eqnarray*}
qui est bien à valeurs dans les $\La[G]$-modules lisses puisque si
$(\mathcal{W}_i)_{i\in I}$ est un recouvrement localement fini par des
ouverts quasicompacts alors 
$$
\GG_! ( \X^{rig},\F)\hookrightarrow \bigoplus_{i\in I} \GG ( \mathcal{W}_i,\F)
$$
Considérons le foncteur dérivé $R\GG_! ( \X^{rig},-)$ à valeurs dans
$\DD^+ ( \La [G]_\lss)$.
D'après le corollaire \ref{DGkygygez380EZDP} et le théorème
\ref{OOPfyzY136Uetet} sa cohomologie calcule la cohomologie à support
compact munie de son action lisse de $G$. 

\section{Cohomologie à support compact équivariante-lisse des tours
  d'espaces analytiques}\label{kjsfyqdfgEZTdytzrz4658}

\subsection{Hypothèses et notations}

Soient $G$ et $H$ deux groupes topologiques possédant un sous-groupe
ouvert profini. On suppose fixé un sous-groupe compact ouvert $G^0$ de
$G$. 

Soit $\mathcal{S}$ la catégorie dont les objets sont les sous-groupes
ouverts de $G^0$,
$$
\forall K,K' \in \text{Ob} (\mathcal{S})\;\; \Hom_{\mathcal{S}} (
K,K') = \{ \bar{g}\in K\bc G\;|\; g^{-1} K g\subset K'\}
$$
et 
$$\forall \bar{g}_1\in \Hom_{\mathcal{S}} ( K,K')\; \forall
\bar{g}_2 \in \Hom_{\mathcal{S}} ( K',K'')\;\; \bar{g}_2\circ
\bar{g_1} = \overline{g_1 g_2}$$
(on vérifie aisément que cette opération de composition est bien
définie). 

Supposons nous donné un foncteur $K\longmapsto X_K$ de $\mathcal{S}$
dans la catégorie des espaces analytiques de Berkovich paracompacts munis d'une
action continue de $H$.
En particulier $\forall K\subset G^0$ $\;\forall g\in G$ tel que
$g^{-1} K g \subset G^0$ il y a un isomorphisme $H$-équivariant
$$
g:X_K\iso X_{g^{-1} K g}
$$
trivial si $g\in K$  
et si $K'\subset K\subset G^0$ il y a un morphisme $H$-équivariant 
$$
\pi_{K',K} : X_{K'}\ldrt X_K
$$
Ainsi si $K'\lhd K$ $\; X_{K'}$ est muni d'une action de $K/K'$ au dessus de $X_K$. 

Nous ferons l'hypothèse supplémentaire suivante : si $K'\lhd K$ alors
$\pi_{K',K} : X_{K'}\ldrt X_K$ est un $K/K'$-torseur étale. On en
déduit que $\forall K'\subset K\; \pi_{K',K}$ est étale fini et que
l'on a un ``pro-$K$-torseur'' $(X_{K'})_{K'\subset K} \ldrt X_K$ et en particulier un ``pro-$G^0$-torseur'' au dessus de l'objet final de notre tour $(X_K)_{K\subset G^0}\ldrt X_{G^0}$.

Il y a également des correspondances de Hecke : si $K\subset G^0$ et
$g\in G$ 
 il y a un diagramme
$$
\xymatrix{
 & X_{K\cap gKg^{-1}}  \ar[r]^g \ar[ld]_{\pi_{K\cap g K g^{-1} ,K}} & X_{g^{-1} K g\cap
   K}\ar[rd]^{\pi_{g^{-1}Kg\cap K,K}}
 \\
X_K & & & X_K
}
$$

\begin{exem}
On peut prendre $X_K =\M_K$ un espace de Rapoport-Zink (\cite{RZ}) où $G=G(\Qp)$,
$G^0$ est un sous-groupe compact maximal et $H=J_b$, par exemple la tour d'espaces de Lubin-Tate considérée dans \cite{Cellulaire} avec $G=\GL_n (F)$ et $H=D^\times$ ou bien celle de Drinfeld avec $G=D^\times$ et $H=\GL_n (F)$. Plus généralement si $\breve{E}$ désigne le complété de l'extension maximal non-ramifiée  du corps reflex associé à l'espace de Rapoport-Zink on peut prendre $X_K = \M_K\hat{\otimes}_{\breve{E}} \widehat{\overline{\breve{E}}}$, $G=G(\Qp)$ et $H=J_b\times \Gal (\overline{\breve{E}}|\breve{E})$. Un autre point de vue pour incorporer l'action de Galois consiste à prendre $G=G(\Qp)\times \Gal (\overline{\breve{E}}|\breve{E})$, $H=J_b$, si $K=K_1\times K_2$ où $K_1\subset G(\Zp)$ et $K_2 = \Gal (\overline{\breve{E}} | M)$, $X_{K_1\times K_2} = \M_{K_1}\otimes_{\breve{E}}  M$ (avec les conditions imposées sur notre tour cela est suffisant pour la définir puisque les groupes de la forme $K_1\times K_2$ forment une base de voisinage de l'identité dans $G(\Qp)\times \Gal (\overline{\breve{E}}|\breve{E})$). 
\end{exem}

\subsection{Faisceaux de Hecke sur la tour}
\subsubsection{Définitions}

Les systèmes $K\longmapsto (X_K)\top_{\et-H-\lss}$, resp. $K\longmapsto
(X_K)\top_{\qetc-H-\lss}$, des $H$-faisceaux  étales lisses, resp. sur le site quasi-étale compact,   forment deux systèmes de topos fibrés au dessus de la catégorie $\mathcal{S}$.

\begin{defi}
Nous noterons $\mathcal{H}_\et$, resp. $\mathcal{H}_\qetc$, les topos
totaux associés aux topos fibrés précédents (\cite{Illusie3} chapitre
VI). On les appellera  faisceaux de Hecke étales, resp. quasi-étales, sur la tour.
\end{defi}

Par exemple $\mathcal{H}_\et$ est la catégorie des $(\F_K)_{K\subset
  G^0}$, où $\F_K$ est un $H$-faisceau étale lisse sur $X_K$, munis de morphismes $H$-équivariants
$$
\forall g\in G\text{ tel que } g^{-1} K g\subset G^0\;\; g^*\F_{g^{-1} K g} \ldrt \F_K
$$
égaux à l'identité si $g\in K$ ,
satisfaisant à certaines conditions de cocyle lorsqu'on les compose et de morphismes 
$H$-équivariants 
compatibles aux morphismes précédents 
$$
\forall K'\subset K\; \; \pi_{K',K}^* \F_K\ldrt \F_{K'}
$$
eux-mêmes soumis à des conditions de compatibilité lorsque l'on a trois groupes $K''\subset K'\subset K$.
On remarquera en particulier que les conditions de cocyle imposent que les morphismes $g^*\F_{g^{-1}Kg}\ldrt \F_K$ soient des isomorphismes.

Dans la suite on notera en abrégé $(\F_K)_K$ un tel objet du topos total, sans faire référence aux morphismes de compatibilité précédents qui seront sous-entendus.

La terminologie faisceaux de Hecke provient de ce que si $(\F_K)_K$
est comme précédemment,  $K\subset G^0$ et $g\in G$ il y a  une correspondance cohomologique de $\F_K$
dans lui-même supportée par la correspondance de Hecke associée à $K$
et $g$ 
$$
g^*\pi_{g^{-1}K g\cap K,K}^* \F_K \ldrt \pi_{K\cap gKg^{-1},K}^*\F_K
$$
et de plus cette correspondance cohomologique ne dépend que de la
double classe $KgK\in K\bc G/K$ au sens qui suit. Si $k_1,k_2\in K$ et
$g'=k_1 gk_2$ il
y a un isomorphisme entre correspondances de $X_K$ dans lui-même
$$
\xymatrix{
&   X_{K\cap g'Kg'^{-1}} \ar[d]_{k_1^{-1}}^\simeq \ar[r]^{g'} \ar@(lu,u)[ldd]_(.6){\pi_{K\cap g' K g'^{-1} ,K}} &  X_{g'^{-1} K g'\cap
   K}\ar[d]^{k_2}_\simeq \ar@(ru,u)[rdd]^(.6){\pi_{g'^{-1}Kg'\cap K,K}} \\
 & X_{K\cap gKg^{-1}}  \ar[r]^g \ar[ld]|{\pi_{K\cap g K g^{-1} ,K}} & X_{g^{-1} K g\cap
   K}\ar[rd]|{\pi_{g^{-1}Kg\cap K,K}}
 \\
X_K & & & X_K
}
$$
Via cet isomorphisme les deux correspondances cohomologiques associées
à $g$ et $g'$ se correspondent.

\begin{defi}
Pour $\bullet\in \{\et,\qetc\}$ on notera
$\mathcal{H}_{\bullet}^{cart}$ les objets cartésiens du topos total
$\mathcal{H}_\bullet$.
\end{defi}

En d'autres termes il s'agit de la sous-catégorie de
$\mathcal{H}_\bullet$ formée des $(\F_K)_K$ tels que 
$$
\forall K'\subset K\subset G^0\;\; \pi_{K',K}^*\F_K\iso \F_{K'}
$$
soit un isomorphisme. En particulier si $K'\lhd K$ on a $\F_K =
(\pi_{K',K*} \F_{K'})^{K/K'}$.
\\
 On notera que le foncteur  $(\F_K)_K \longmapsto
\F_{G^0}$ de
 $\mathcal{H}_\bullet^{cart}$ dans le topos $(X_{G^0})\top_{H-\lss}$ des faisceaux sur la base
de la tour est fidèle.

\begin{exem}
Considérons l'espace de Rapoport-Zink des déformations d'un groupe
$p$-divisible étale muni d'une $G$-structure où $G$ est un groupe
réductif non-ramifié sur $\Qp$. On a alors $G=G(\Qp) = H$, $G^0=
G(\Zp)$. 
Alors si $L=\widehat{\Qp^{nr}}$ on a
$$
\M_K= \coprod_{K\bc G} \mathcal{M} (L)
$$
sur lequel $H$ agit à gauche et $G$ à droite par correspondances de
Hecke. Les faisceaux de Hecke cartésiens sur la tour s'identifient alors aux
ensembles munis d'une action lisse de $G(\Qp)\times \Gal
(\overline{L}|L)$ (il s'agit d'un cas très particulier de la
proposition \ref{DVKET479dsgf2} puisque dans ce cas là l'espace des
périodes est réduit à un point).
\end{exem}

\subsubsection{Propriétés de base des catégories de faisceaux de Hecke}

Il y a des triplets de foncteurs adjoints pour tout $K\subset G^0$ et $\bullet\in \{\et,\qetc \}$
$$
\xymatrix@C=14mm{
(X_K)_{\bullet-H-\lss} \ar@<1.2ex>[r]^(.6){i_{K*}} \ar@<-1.2ex>[r]_(.6){i_{K!}} & \mathcal{H}_\bullet \ar[l]|(.4){i_K^*}
}
$$
où 
\begin{itemize}
\item $i_K^* ((\F_{K'})_{K'}) =\F_K $ qui est le foncteur ``restriction à un étage''
\item $\dpt{(i_{K*} \F)_{K'} = \prod_{\a: K\underset{\mathcal{S}}{\drt}
    K'} \a_*\F_K}$ (produit dans la catégorie des
  $H$-faisceaux lisses, cf. remarque \ref{dsgjdsguuet25YETG})
\item $\dpt{(i_{K!} \F)_{K'} = \bigoplus_{\a: K'\underset{\mathcal{S}}{\drt} K} \a^*\F_K}$
\end{itemize}
(cf. \cite{Illusie3} chapitre VI section 5). On notera en particulier que l'existence de l'adjoint à droite $i_{K!}$ de $i_K^*$ implique que si $(\F_K)_K$ est un objet injectif de $\La-\mathcal{H}_\bullet$ alors $\forall K\; \F_K$ en est un dans $\La-(X_K)_{\bullet-H-\lss}$. 

Si $K'\ldrt K$ est un morphisme dans $\mathcal{S}$ il y a un diagramme 1-commutatif de morphisme de topos
$$
\xymatrix{
(X_{K'})_{\qetc-H}\top \ar[d]\ar[r]^{(\nu^*,\nu_*)} & (X_{K'})_{\et-H}\top \ar[d] \\
(X_{K})_{\qetc-H}\top \ar[r]^{(\nu^*,\nu_*)} & (X_{K})_{\et-H}\top 
}
$$
(utiliser le fait que le morphisme $X_{K'}\ldrt X_K$ est étale fini) qui induit un morphisme de topos fibrés d'où un morphisme de topos
$$
(\nu^*,\nu_*) : \mathcal{H}_\qetc \ldrt \mathcal{H}_\et
$$
qui commute à la restriction aux étages : $\nu^* ((\F_K)_K) = (\nu^* \F_K)_K$ et $\nu_* ((\F_K)_K) = (\nu_*\F_K)_K$. Ce morphisme satisfait donc à toutes les propriétés déjà énoncées dans le cas d'un espace de Berkovich seul (i.e. pas dans le cas d'une tour). En particulier $\nu^* : \mathcal{H}_\et \ldrt \mathcal{H}_\qetc$ est pleinement fidèle et il y a un isomorphisme $Id\iso \nu_*\nu^*$. 

\begin{lemm}\label{FSKDGKZERTY36458HJB}
Les catégories $\La-\mathcal{H}_{\qetc}$ et $\La-\mathcal{H}_{\et}$ sont abéliennes et possèdent suffisamment d'injectifs. De plus si $(\F_K)_K$ est un objet injectif de l'une de ces catégories alors $\forall K\; \F_K$ est un objet injectif de la catégorie des faisceaux $H$-lisses sur $X_K$ correspondante. Si $(\F_K)_K \in \mathcal{H}_{\et}$ et $\nu^* (\F_K)_K \ldrt (\mathcal{I}^\bullet_K )_K$ est une résolution injective  dans $\mathcal{H}_{\qetc}$ alors 
$(\F_K)_K \ldrt \nu_* (\mathcal{I}^\bullet_K)_K$ en est une dans $\mathcal{H}_{\et}$ telle que les $\nu_* \mathcal{I}^\bullet_K$ soient des faisceau mous.
\end{lemm}

\subsubsection{Propriétés de base des catégories de faisceaux de Hecke cartésiens}
\label{DGKDG23565UGJETOP8}

Soit $\bullet\in \{\et,\qetc\}$. L'inclusion $\mathcal{H}_\bullet^{cart}\hookrightarrow \mathcal{H}_\bullet$ possède un adjoint à gauche $\kappa$ défini par $\forall (\F_K)_K \;\; \kappa ((\F_K)_K) =(\G_K)_K$ où
$$
\G_K = \underset{K'\lhd K}{\limi} (\pi_{K',K*}\F_{K'})^{K/K'}
$$
dont on vérifie facilement que c'est naturellement un faisceau de Hecke cartésien. De plus $\kappa$ est un foncteur exact. 

\begin{prop}
La catégorie $\mathcal{H}_\bullet^{cart}$ est un topos. 
\end{prop}
\dem
L'existence de limites projectives finies résulte de ce que si $K'\subset K\subset G^0$ $\;\pi_{K',K}^*$ commute à ces limites et donc si $((\F_{i,K})_K)_{i\in I}$ est un système projectif fini de faisceaux cartésiens sa limite projective dans $\mathcal{H}_\bullet$,
$(\underset{i\in I}{\limp} \F_{i,K})_K$, est cartésienne. Il en est de même pour les limites inductives quelconques et les ``bons'' quotients par des relations d'équivalence. Reste à voir l'existence d'une famille génératrice. Mais si $(C_\a)_{\a\in A}$ en est une de $\mathcal{H}_\bullet$ alors $(\kappa (C_\a))_{\a\in A}$ en est une de $\mathcal{H}_\bullet^{cart}$. 
\qed
\\

Il résulte des discussions précédentes que si $\F\ldrt \mathcal{I}^\bullet$ est une résolution injective dans $\La-\mathcal{H}_\bullet^{cart}$ alors c'en est une dans $\La-\mathcal{H}_\bullet$. Comme dans la section précédente il y a un morphisme de topos $(\nu^*,\nu_*):\mathcal{H}^{cart}_{\qetc}\ldrt \mathcal{H}^{cart}_{\et}$. Et le lemme \ref{FSKDGKZERTY36458HJB} reste valable dans ce contexte en remplaçant $\mathcal{H}$ par $\mathcal{H}^{cart}$. 

\subsection{Le complexe de cohomologie à support compact de la tour}\label{DGKSRY246EHENrf}

Si $(\F_K)_K\in \mathcal{H}_{\et-\lss}$ et $K'\subset K$ étant donné que $X_{K'}\ldrt X_K$ est étale fini il y a un morphisme
$$
\GG_c (X_K,\F_K)\xrig{\; \pi_{K',K}^*\;} \GG_c ( X_{K'},\pi_{K',K}^* \F_K ) \ldrt \GG_c ( X_{K'},\F_{K'})
$$
où le deuxième morphisme est induit par les morphismes de transition dans la définition du topos total. De plus on a vu dans la section \ref{kfuhDGT23576gffgf} que ce sont des $H$-modules lisses. On obtient ainsi un système inductif de $H$-modules lisses muni d'une action de $G$.

\begin{defi}
On note 
$\GG_c (\mathcal{H},-) : \La-\mathcal{H}_{\et-\lss} \ldrt \La[G\times H]_{\lss}$ le foncteur qui à $(\F_K)_K$ associe
$$
\underset{K\subset G^0}{\limi} \GG_c ( X_K,\F_K)
$$
On note $R\GG_c (\mathcal{H}_\et,-)$ le foncteur dérivé correspondant à valeurs dans $\DD^+ ( \La[G\times H]_{\lss})$. 
\end{defi}

\begin{theo}
Le foncteur $(\F_K)_K\longmapsto R\GG_c ( \mathcal{H}, (\F_K)_K)$ est tel que 
$$
\forall i\; H^i ( R\GG_c (\mathcal{H}_\et,(\F_K)_K)) = \underset{K}{\limi} H^i_c (X_K,\F_K)
$$
\end{theo}
\dem
C'est une conséquence du lemme \ref{FSKDGKZERTY36458HJB} qui affirme
que $(\F_K)_K$ possède une résolution injective qui après oubli de l'action de $H$ est molle étage par étage.
\qed

\begin{rema}
D'après les résultats de la section \ref{DGKDG23565UGJETOP8} le foncteur dérivé précédent évalué sur un faisceau de Hecke cartésien peut se calculer avec une résolution injective de faisceaux cartésiens. C'est par exemple le cas du faisceau constant $\underline{\La}$. 
\end{rema}

\begin{rema}
Si l'ordre du groupe profini $K_0$, $K_0\subset G^0$, est inversible dans $\La$ alors le foncteur $M\longmapsto M^{K_0}$ de $\La[G\times H]_{\lss}\ldrt \La [H]_{\text{\lss}}$ est exact et $R\GG_c (\mathcal{H}_\et, (\F_K)_K)^{K_0} = R\GG_c (X_{K_0},\F_{K_0})$. Ce foncteur d'invariants sous $K_0$ se factorise en fait en un foncteur $\La[G\times H]_{\lss} \ldrt \mathcal{H} (K_0\bc G/K_0)\otimes_{\La} \La[H]_{\lss}$ où $\mathcal{H} (K_0\bc G/K_0)$ désigne l'algèbre de Hecke des fonctions bi-invariantes sous $K_0$ à valeurs dans $\La$ et $R\GG_c (\mathcal{H}_\et, (\F_K)_K)^{K_0}$ est le complexe de modules sur l'algèbre de Hecke donné par l'action des correspondances de Hecke cohomologiques. 
\end{rema} 
 
\subsection{Objets cartésiens sur la tour et domaine fondamental pour
  l'action des correspondances de Hecke}\label{SFSjygez26Uetejej3}
  
\subsubsection{Introduction : recollement de faisceaux équivariants}\label{SDGKIR54Y66Ttf}

Soit $X$ un espace topologique muni d'un recouvrement ouvert
$(U_i)_{i\in I}$. Notons $U_{ij}= U_i\cap U_j$ et $U_{ijk}= U_i\cap
U_j\cap U_k$. 
Le topos $X\top$ des faisceaux sur $X$ est
équivalent aux objets cartésiens du topos total du diagramme suivant 
$$
\xymatrix{
\coprod_{i,j,k} U_{ijk}\top \ar@<1.2ex>[r] \ar[r] \ar@<-1.2ex>[r] &
\coprod_{i,j} U_{ij}\top \ar@<0.8ex>[r]\ar@<-0.8ex>[r] & \coprod_{i} U_i\top
}
$$
Traduit en termes usuels cela signifie que se donner un faisceau sur
$X$ est équivalent à se donner une collection $(\F_i)_{i\in I}$ de
faisceaux sur les $(U_i)_i$ munis d'isomorphismes $\b_{ij}:
\F_{i|U_{ij}}\iso \F_{j|U_{ij}}$ vérifiant $\b_{ij|U_{ijk}}\circ
\b_{jk|U_{ijk}} =\b_{ik|U_{ijk}}$. 

Soit maintenant $G$ un groupe agissant sur $X$. Se donner un
$G$-faisceau sur $X$ est équivalent à se donner un objet cartésien du
topos total associé au diagramme usuel
$$
\xymatrix{
G\times G\times X\top  \ar@<1.2ex>[r] \ar[r] \ar@<-1.2ex>[r] & G\times X\top
\ar@<0.8ex>[r]\ar@<-0.8ex>[r] & X\top
}
$$
c'est à dire des isomorphismes $\forall g\in G\; \a_g:g^*\F\iso \F$
vérifiant $\forall g_1,g_2\; \a_{g_2}\circ g_2^*\a_{g_1} =\a_{g_1
  g_2}$. 

Supposons maintenant que $G$ agisse sur $I$ et que $\forall i\; g.U_i
=U_{g.i}$. Se donner un faisceau $G$-équivariant sur $X$ est alors
équivalent à se donner un objet cartésien du topos total associé au diagramme
$$
\xymatrix{
\coprod_{i,j,k} U_{ijk}\top \ar@<1.2ex>[r] \ar[r] \ar@<-1.2ex>[r] &
\coprod_{i,j} U_{ij}\top \ar@<0.8ex>[r]\ar@<-0.8ex>[r] & \coprod_{i}
U_i\top \\
 & G\times \coprod_{i,j} U_{ij}\top  \ar@<0.8ex>[r]\ar@<-0.8ex>[r]
 \ar@<0.8ex>[u]\ar@<-0.8ex>[u]   &
 G\times \coprod_{i} U_i\top   \ar@<0.8ex>[u]\ar@<-0.8ex>[u] \\
& & G\times G\times  \coprod_{i} U_i\top  \ar@<1.2ex>[u] \ar[u] \ar@<-1.2ex>[u]
}
$$
où on laisse le lecteur deviner les différents flèches de ce diagramme.
Traduit en termes usuels cela signifie que se donner un $G$-faisceau
sur $X$ est équivalent à se donner des faisceaux $(\F_i)_{i\in I}$ sur
les $(U_i)_i$ munis 
\begin{itemize}
\item
d'isomorphismes $\b_{ij} :\F_{i|U_{ij}}\iso \F_{j|U_{ij}}$ vérifiant $\b_{ij|U_{ijk}}\circ
\b_{jk|U_{ijk}} =\b_{ik|U_{ijk}}$
\item
 d'isomorphismes $\forall g,i\;\;
\a_{g,i} : g^* \F_{g.i}\iso \F_i$ vérifiant $\a_{g_1,i}\circ
g_1^*\a_{g_2,g_2 g_1.i} =\a_{g_2 g_1,i}$
\item enfin on demande que le diagramme suivant commute 
$$
\xymatrix{
g^* ( \F_{g.i|U_{g.i,g.j}}) \ar[d]_{g^* \b_{g.i,g.j}} \ar[r]^(.6){\a_{g,i}} & \F_{i|U_{ij}}
\ar[d]^{\b_{ij}} \\
g^* (\F_{g.j|U_{g.i,g.j}}) \ar[r]^(.6){\a_{g,j}} & \F_{j|U_{ij}} 
}
$$
\end{itemize}
\vspace{4mm}

Soit maintenant $\forall i\in I \; G_i=\text{Stab}_G (U_i)$. Posons $\forall A\subset I$ $\; G_A=\bigcap_{i\in A} G_i$ qui stabilise donc $\bigcap_{i\in A} U_i$. Pour un $i\in I$ notons $G_i-U_i\top$ le topos des $G_i$-faisceaux sur $U_i$.

Pour $g\in G$ on a $G_{g.i}=g G_i g^{-1}$ et il y a une équivalence
$$
g^* : G_{g.i}-U_{g.i}\top \iso G_i-U_i\top
$$
induite par le morphisme d'espaces équivariants $(g,\text{int}_g): (U_i,G_i)\ldrt (U_{g.i},G_{g.i})$. Si $g\in G_i$ cette équivalence est le foncteur 
\begin{eqnarray*}
G_i-U_i\top &\iso & G_i-U_i\top \\
\F & \mapsto & \F^g
\end{eqnarray*}
où $\F^g$ est le $G_i$-faisceau ``tordu par l'automorphisme intérieur'' $\text{int}_g$ de $G_i$.

Pour un tel $g\in G_i$ il y a un isomorphisme de foncteurs 
$$
\gamma_{i,g} : (-)^g\iso Id
$$
donné pour un $\F\in G_i-U_i\top $ par l'action de $g$ : $g^*\F\ldrt \F$. 

Soit maintenant $(\F_i)_{i\in I}\in \prod_{i\in I} G_i-U_i\top$ une collection de faisceaux équivariants. Se donner un $G$-faisceau $\F$ sur $X$ tel que $\forall i\in I\; \F_{|U_i}\simeq \F_i$ comme $G_i$-faisceaux est équivalent à la donnée 
\begin{itemize}
\item d'isomorphismes de $G_{ij}$-faisceaux $\b_{ij} :\F_{i|U_{ij}}\iso \F_{j|U_{ij}}$ vérifiant $\b_{ij|U_{ijk}}\circ
\b_{jk|U_{ijk}} =\b_{ik|U_{ijk}}$
\item
 d'isomorphismes de $G_i$-faisceaux   $\forall g,i\;\;
\a_{g,i} : g^* \F_{g.i}\iso \F_i$ vérifiant $\a_{g_1,i}\circ
g_1^*\a_{g_2,g_2 g_1.i} =\a_{g_2 g_1,i}$ et tels que si $g\in G_i$ on ait
$\a_{g,i} = \gamma_{i,g}$ l'isomorphisme donné par la structure de $G_i$-faisceau sur $\F_i$
\item enfin on demande que le diagramme suivant commute 
$$
\xymatrix{
g^* ( \F_{g.i|U_{g.i,g.j}}) \ar[d]_{g^* \b_{g.i,g.j}} \ar[r]^(.6){\a_{g,i}} & \F_{i|U_{ij}}
\ar[d]^{\b_{ij}} \\
g^* (\F_{g.j|U_{g.i,g.j}}) \ar[r]^(.6){\a_{g,j}} & \F_{j|U_{ij}} 
}
$$
\end{itemize}

\subsubsection{Hypothèses et notations}

Nous ferons désormais les hypothèses suivantes.
Il existe un domaine analytique compact $W$ dans la base de notre tour
$X_{G^0}$ tel que 
\begin{itemize}
\item Lorsque $(g,h)\in G\times H$ parcourt $G\times H$ l'image par
  les correspondances de Hecke 
$$
\xymatrix{
 & X_{G^0\cap g G^0 g^{-1}} \ar[r]^{(g,h)}\ar[ld]_{\pi_{G^0\cap g G^0
     g^{-1},G^0}} 
 & X_{g^{-1} G^0 g\cap G^0}
 \ar[rd]^{\pi_{g^{-1} G^0 g\cap G^0,G^0}} \\
X_{G^0} & && X_{G^0}
}
$$
 de $W$, les $\pi_{g^{-1} G^0 g\cap G^0,G^0} \left ( (g,h). \pi_{G^0\cap g G^0
     g^{-1},G^0}^{-1} (W)\right )$, recouvrent $X_{G^0}$.
\item Soit $H^0$ le  stabilisateur dans $H$ de $W$ (qui est nécessairement un
  sous-groupe ouvert). On suppose l'existence d'un sous-groupe central
  $Z$ dans $G\times H$ tel que $Z$ agisse trivialement sur la
  tour. Notons alors  $$\mathcal{I} = (G^0\bc G\times H/H^0)/Z$$
 et  $\forall z\in\mathcal{I} \;\; z.W\subset X_{G^0}$
  l'image de $W$ donnée par la correspondance de Hecke $z$. On suppose
  alors que 
$$
\forall z\in \mathcal{I}\;\; \{z'\in \mathcal{I}\;|\; z'.W\cap z.W
\neq \emptyset \} \text{ est fini}
$$
\end{itemize}

\begin{exem}
Prenons le cas des espaces de Lubin-Tate étudié dans
\cite{Cellulaire}. Avec les notations de l'article \cite{Cellulaire}
l'espace de Rapoport-Zink sans niveau $\M$ s'écrit $\M= \coprod_\Z {\X}^{an}$ où $\X$
désigne l'espace de Lubin-Tate. Alors le domaine fondamental de
Gross-Hopkins dans la composante indexée par $0\in \Z$ satisfait aux
hypothèses demandée à $W$. L'ensemble $\mathcal{I}$ consiste en les sommets de
l'immeuble introduit dans \cite{Cellulaire}.
\end{exem}

\begin{exem}
Dans le cas de l'espace de Rapoport-Zink associé à l'espace de
Drinfeld $\Omega$, $\M=\coprod_\Z \Omega$ on peut prendre pour $W$
l'image réciproque d'un simplexe maximal dans l'immeuble via la
rétraction de $\Omega$ sur l'immeuble.
\end{exem}

\subsubsection{Les topos équivariants fibres en un sommet de ``l'immeuble''}

\begin{defi}
Soit $z=[g,h]\in \mathcal{I}$. On note $G_z=g^{-1} G^0 g$,
$H_z=h H^0 h^{-1}$. Pour $K\subset G_z\cap G^0$ un sous-groupe ouvert on note 
$$
W_{z,K} = (g,h). \pi_{gKg^{-1},G^0}^{-1} (W) \subset X_K
$$
$$
\xymatrix{
 & \ar[ld]^{\pi_{gKg^{-1},G^0}}  X_{gKg^{-1}} \ar[r]^{(g,h)} & X_K \\
X_{G^0}
}
$$
De même si $I\subset \mathcal{I}$ est un sous-ensemble fini on note
$G_I=\bigcap_{z\in I} G_z$, $H_I =\bigcap_{z\in I} H_z$ et pour
$K\subset G_I\cap G^0$ $\; W_{I,K} = \bigcap_{z\in i} W_{z,K}\subset X_K$.  
\end{defi}

Remarquons que si $K\subset G_I\cap G^0$, $g\in G_I$ et $h\in H_I$
sont tels que $g^{-1} K g \subset G^0$ alors $g^{-1} K g \subset
G_I\cap G^0$ et $(g,h)$ induit un isomorphisme
$$
(g,h) : W_{I,K}\iso W_{I,g^{-1} K g}
$$
\\
On obtient ainsi une tour $(W_{I,K})_{K\subset G_I\cap G^0}$ satisfaisant aux
même conditions que la tour $K\longmapsto X_K$ en remplaçant $G$ par
$G_I$, $G^0$ par $G_I\cap G^0$ et $H$ par $H_I$. Cette tour est plus simple
car pour celle-ci $G_I$ est compact. 

\begin{defi}
Soit $I\subset \mathcal{I}$ un ensemble fini. On note
$\mathcal{H}_{\qetc,I}^{cart}$ la catégorie des objets cartésiens du topos total du
topos fibré $G_I\cap G^0 \supset K \longmapsto (W_{I,K})_{\qetc-H_I-\lss}\top$ au
dessus de la catégorie des sous-groupes ouverts de $G_I\cap G^0 $ munis des
morphismes $\Hom (K,K') = \{g\in K\bc G_I \;|\; g^{-1} K g \subset
K'\}$. On note de même $\mathcal{H}_{\et,I}^{cart}$ en remplaçant
quasi-étale par étale. 
\end{defi}

\begin{lemm}\label{GRYRT2345ffdz}
Soit $\bullet\in \{\et,\qetc\}$. 
Soit $I$ comme précédemment et $K_0\subset G_I\cap G^0$ tel que 
$K_0\lhd G_I$. Alors l'application qui à $(\F_K)_{K\subset G_I\cap G^0} \in
\mathcal{H}_{\bullet, I}^{cart}$ associe $\F_{K_0}$ muni de son action de
$G_I/K_0\times H_I$ compatible à celle sur $W_{I,K_0}$
est une équivalence de
catégories entre $\mathcal{H}_{\bullet, I}^{cart}$ et
$(W_{K_0,I})_{\bullet,G_I/K_0\times H_I-\lss}\top$. 
\end{lemm}
\dem
C'est une conséquence du rappel du début de la section \ref{ldddhzZRHREY457H} sur les
faisceaux équivariants sur un torseur. En effet, si $K\subset K_0$ il
suffit de poser $\F_K= \pi_{K,K_0}^*\F_{K_0}$. Pour $K\subset G_I\cap
G^0$ général il
suffit de choisir $K'\lhd G_I$ tel que $K'\subset K_0$ et $K'\subset
K$. Alors
$$
\F_K= \left (\pi_{K',K*} ( \pi_{K',K_0}^*\F_{K_0})\right )^{K/K'}
$$
On vérifie alors que l'action de $G_I/K_0$ sur $\F_{K_0}$ définit bien
un objet cartésien de notre topos fibré. 
\qed 

On a même le lemme suivant.

\begin{lemm}
Soient $\bullet\in \{\et,\qetc\}$
et $I$ comme précédemment. Soit $K\subset G_I\cap G^0$ tel que $K\lhd G_I$. Notons $W_I :=W_{I,K}/(G_I/K)$. Alors $W_{I,K}\ldrt W_I$ est un $(G_I/K)$-torseur étale et il y a une équivalence
$$
\mathcal{H}^{cart}_{\bullet,H_I-\lss} \iso (W_I)\top_{\bullet,H_I-\lss}
$$
\end{lemm}
\dem
Elle ne pose pas de problème.
\qed

\begin{exem}
Dans le contexte de l'article \cite{Cellulaire} on trouve donc que pour un sommet $z$ de l'immeuble les objets de $\mathcal{H}^{cart}_{\bullet,z}$ s'identifient aux $\O_D^\times$-faisceaux lisses sur la cellule sans niveau associée au sommet $z$. 
\end{exem}

Il y a une action de $G\times H$ sur $\mathcal{I}$ (l'action de $G$ se
faisant à droite) telle que $\forall I \; G_{(g,h).I} = g^{-1} G_I g$ et
$H_{(g,h).I} = h H_I h^{-1}$. Cette action induit un isomorphisme pour
$I$ et $K\subset G_I\cap G^0$ suffisamment petits
$$
(g,h) : W_{I,K}\iso W_{(g,h).I,g^{-1} K g}
$$
Cela n'induit pas un isomorphisme entre les tours $(W_{I,K})_{K\subset
  G_I\cap G^0}$ et $(W_{(g,h).I,K})_{K\subset G_{(g,h).I}\cap G^0}$
car il faut se restreindre à des $K$ suffisamment petits.
Néanmoins d'après le lemme \ref{GRYRT2345ffdz} cela induit une équivalence entre
les systèmes de faisceaux cartésiens sur les tours 
$$
((g,h )_*,(g,h)^*): \mathcal{H}^{cart}_{\bullet,I} \iso
\mathcal{H}^{cart}_{\bullet, (g,h).I}
$$
De plus si $J\subset I$ il y a une application de restriction 
\begin{eqnarray*}
\mathcal{H}_{\bullet, J}^{cart} &\ldrt & \mathcal{H}_{\bullet, I}^{cart} \\
(\F_K)_{K\subset G_J\cap G^0} & \longmapsto & (\F_{K|W_{I,K}})_{K\subset G_I\cap G^0}
\end{eqnarray*}

\subsubsection{Reconstruction des systèmes de faisceaux cartésiens sur la
  tour à partir de ceux sur les tours en chaque sommet de l'immeuble}

\begin{theo}\label{jytarmsvioizr35zsfjA}
Il y a une équivalence de catégories entre les systèmes cartésiens de faisceaux quasi-étales sur la tour, $\mathcal{H}^{cart}_{\qetc}$, est les  $(\F_i)_{i\in\mathcal{I}}$, où $\F_i\in \mathcal{H}^{cart}_{\qetc,i}$
\begin{itemize}
\item munis d'isomorphismes $\forall i,j\; \b_{ij}: \F_{i|\mathcal{H}^{cart}_{ij}} \iso \F_{j|\mathcal{H}^{cart}_{ij}}$ satisfaisant 
$$
\beta_{ij|\mathcal{H}^{cart}_{ijk}} \circ \beta_{jk|\mathcal{H}^{cart}_{ijk}} = \beta_{ik|\mathcal{H}^{cart}_{ijk}}
$$
\item munis d'isomorphismes $\forall (g,h)\in G\times H\;\forall i\;\; \a_{(g,h),i} :(g,h)^*\F_{(g,h).i}\iso \F_i$ satisfaisant
$$
\forall (g,h),(g',h')\;\forall i\;\; \a_{(g',h'),(g,h).i} \circ (g,h)^*\a_{(g,h),i} = \a_{(g'g,h'h),i}
$$
et tels que si $(g,h)\in G_i\times H_i$ alors $\a_{(g,h),i}$ soit  l'action de $G_i\times H_i$ sur le
système de faisceau $\F_i$ 
\item et enfin tels que le diagramme suivant commute
$$
\xymatrix@C=20mm@R=12mm{
((g,h)^*\F_{(g,h).i})_{|\mathcal{H}^{cart}_{ij}} \ar[r]^{(\a_{(g,h),i})_{|\mathcal{H}_{ij}^{cart}}} \ar[d]_{(g,h)^*\beta_{(g,h).i,(g,h).j}} & \F_{i|\mathcal{H}_{ij}^{cart}} \ar[d]^{\beta_{ij}} \\
((g,h)^*\F_{(g,h).j})_{|\mathcal{H}^{cart}_{ij}} \ar[r]^{(\a_{(g,h),j})_{|\mathcal{H}_{ij}^{cart}}} &  \F_{j|\mathcal{H}_{ij}^{cart}}
}
$$
\end{itemize}
Lorsque les images par les correspondances de Hecke de $\overset{\circ}{W}$ recouvrent $X_{G^0}$ on peut remplacer quasi-étale par étale.
\end{theo}
\dem
Il y a un foncteur qui à $(\F_K)_{K\subset G^0}\in\mathcal{H}^{cart}_{\qetc}$ associe la collection des $(\F_{K|W_{i,K}})_{K\subset G_i\cap G^0} \in \mathcal{H}_{\qetc,i}^{cart}$ lorsque $i$ varie dans $\mathcal{I}$ et dont on vérifie que c'est bien un objet cartésien comme dans l'énoncé du théorème.

Montrons que ce foncteur induit une équivalence. 
Soit donc $(\F_{i,K})_{i\in \mathcal{I},K\subset G_i\cap G^0}$ un système de faisceaux comme dans l'énoncé. \'Ecrivons 
$$
\mathcal{I} = \bigcup_{k\geq 1} A_k
$$
où $\forall k\; A_k$ est un ensemble fini, $A_k\subset A_{k+1}$ et $(G^0\times H^0).A_k=A_k$.
Alors
$$
X_{G^0} = \bigcup_{k\geq 1\atop K\subset G_{A_k}\cap G^0} \pi_{K,G^0}\left (  \bigcup_{i\in A_k} W_{i,K} \right )
$$
où dans l'union précédente pour un $k$ donné $K\subset G_{A_k}\cap G^0$ est quelconque, on peut prendre en particulier $K_k := G_{A_k}\cap G^0\lhd G^0$. Il y a un revêtement galoisien de groupe $G^0/K_k$
$$
\xymatrix@C=2cm{
\bigcup_{i\in A_k} W_{i,K_k} \ar[d] \ar@{^(->}[r]^{\text{domaine analytique}\atop\text{compact}} & X_{K_k} \ar[d] \\
\pi_{K_k,G^0} \left (  \bigcup_{i\in A_k} W_{i,K} \right ) \ar@{^(->}[r] & X_{G^0}
}
$$
Les $H^0$-faisceaux quasi-étales lisses $(\F_{i,K_k})_{i\in A_k}$ se recollent grâce aux données de recollement de l'énoncé en un $G^0/K_k\times H^0$-faisceau lisse quasi-étale $\G_k$ sur 
$\bigcup_{i\in A_k} W_{i,K_k}$. Ce faisceau descend en un $H^0$-faisceau lisse $\mathcal{H}_k =(\pi_{K_k,G^0*} \G_k)^{G^0/K_k}$ sur $\pi_{K_k,G^0} \left (  \bigcup_{i\in A_k} W_{i,K} \right )$. Lorsque $k'\geq k$ il y a des identifications $\mathcal{H}_{k'|\pi_{K_k,G^0}(\cup_{i\in A_k} W_{i,K_k})} = \mathcal{H}_k$. Cela définit par recollement  un $H^0$-faisceau lisse $\F_{K^0}$ sur $X_{G^0}$. On pose alors $\forall K\subset G^0\;\; \F_K=\pi_{K,G^0}^* \F_{G^0}$.

On doit maintenant montrer que ce système de faisceaux est muni d'une action de $G\times H$, c'est à dire définir des isomorphismes naturels $\forall K\subset G^0\;\forall (g,h)$  tels que $g^{-1} K g\subset G^0$$\;\; (g,h)^*\F_{g^{-1} K g}\iso \F_K$.  Cela est fastidieux mais ne pose pas de problèmes particuliers. On en laisse donc la vérification au lecteur.
\qed

\begin{rema}
Quitte à remplacer $W$ par un nombre fini de ses itérés sous des
correspondances de Hecke on peut toujours supposer que les itérés sous
les correspondances de Hecke de $\overset{\circ}{W}$ recouvrent
$X_{G^0}$. 
\end{rema}

\subsection{Faisceaux cartésiens sur la tour et espaces de périodes}\label{dgskbds46jtjez}

Supposons maintenant l'existence d'un morphisme étale d'espaces
analytiques
$$
\Pi : X_{G^0}\ldrt P
$$
tel que $P$ soit muni d'une action lisse de $H$ et $\Pi$ commute à
cette action. Supposons de plus que $\Pi$ est $G$-invariant au sens où
$\forall g\in G$ le diagramme suivant commute
$$
\xymatrix@C=8mm{
 & X_{G^0\cap g G^0 g^{-1}} \ar[ld]_{\pi_{G^0\cap g G^0 g^{-1},G^0}} \ar[rr]^g && X_{g^{-1} G^0 g\cap
   G^0} \ar[rd]^{\pi_{g^{-1} G^0 g\cap
   G^0,G^0}} \\
X_{G^0} \ar[rrd]_\Pi &&& & X_{G^0} \ar[lld]^\Pi \\
 & & P
}
$$
On peut alors définir deux foncteurs pour $\bullet \in \{\et,\qetc \}$

\begin{eqnarray*}
\a: P\top_{\bullet-H-\lss} & \ldrt & \mathcal{H}_\bullet^{cart} \\
\G & \longmapsto & \left ( (\Pi\circ \pi_{K,G^0})^*\G\right
)_{K\subset G^0}
\end{eqnarray*}

\begin{eqnarray*}
\beta : \mathcal{H}_\bullet^{cart} &\ldrt & P\top_{\bullet-H-\lss} \\
(\F_K)_{K\subset G^0} & \longmapsto & \left [ \underset{K\subset
    G^0}{\limi} (\Pi\circ \pi_{K,G^0})_* \F_K \right ]^{G}
\end{eqnarray*}
où l'action de $G$ sur $ \underset{K\subset
    G^0}{\limi} (\Pi\circ \pi_{K,G^0})_* \F_K $ se fait de la façon
  suivante. Pour $g\in G$ et $K$ suffisamment petit tel que $K\subset
  G^0\cap g^{-1} G^0 g$ il y a un isomorphisme $\F_K\simeq
  g_*\F_{gKg^{-1}}$ qui induit 
\begin{eqnarray*}
(\Pi\circ \pi_{K,G^0})_*\F_K &\simeq & (\Pi\circ \pi_{K,G^0})_* (g_*
\F_{gKg^{-1}}) \\
 &=& (\Pi\circ \pi_{K,G^0}\circ g)_* \F_{gKg^{-1}} \\
&=& (\Pi\circ \pi_{gKg^{-1},G^0})_* \F_{gKg^{-1}} \ldrt \underset{K'\subset
    G^0}{\limi} (\Pi\circ \pi_{K',G^0})_* \F_{K'}
\end{eqnarray*}
On vérifie que lorsque $K$ varie ces morphismes sont compatibles et
que donc $g$ induit un automorphisme du faisceau limite inductive. 

Une autre description du foncteur $\beta$ en termes de correspondances
de Hecke sphériques plutôt que d'invariants sous $G$ est la suivante.
Si $K\subset G^0$ étant donné que $\pi_{K,G^0}^*\F_{G^0} \iso \F_K$ il
y a un monomorphisme $\F_{G^0}\hookrightarrow (\pi_{K,G^0})_*\F_K$ tel
que si $K\lhd G^0$ on ait $\left [ (\pi_{K,G^0})_*\F_K\right ]^{G^0}
=\F_{G^0}$. Donc
\begin{eqnarray*}
\beta ((\F_K)_K) = \left [ \underset{K\subset
    G^0}{\limi} (\Pi\circ \pi_{K,G^0})_* \F_K \right ]^{G} \subset  \left [ \underset{K\subset
    G^0}{\limi} (\Pi\circ \pi_{K,G^0})_* \F_K \right ]^{G^0} &= &  \underset{K\lhd
    G^0}{\limi} [(\Pi\circ \pi_{K,G^0})_* \F_K]^{G^0} \\ &=&\Pi_* \F_{G^0}
\end{eqnarray*}
Maintenant si $g\in G$ il y a deux flèches 
$$
\xymatrix{
 \Pi_*\F_{G^0} \ar@{=}[r] \ar[rrd]_{\delta_g} & (\Pi \circ \pi_{G^0\cap
  g^{-1} G^0 g,G^0})_*\F_{G^0\cap
  g^{-1} G^0 g} \ar[r]^\sim &   (\Pi \circ \pi_{G^0\cap
  g^{-1} G^0 g,G^0})_* g_*\F_{gG^0 g^{-1}} \ar@{=}[d] \\
& & (\Pi\circ \pi_{gG^0
  g^{-1} \cap G^0,G^0})_* \F_{gG^0
  g^{-1} \cap G^0}
}
$$
ainsi que la flèche tautologique 
$$
\eta_g: \Pi_*\F_{G^0}\hookrightarrow  (\Pi\circ \pi_{gG^0
  g^{-1} \cap G^0,G^0})_* \F_{gG^0
  g^{-1} \cap G^0} 
$$
Alors 
$$
\b ((\F_K)_K) = \bigcap_{\bar{g}\in G/G^0} 
\ker \left ( \xymatrix{ \Pi_*\F_{G^0}\ar@<.6ex>[r]^(.25){\eta_g}\ar@<-.6ex>[r]_(.25){\delta_g} &   (\Pi\circ \pi_{gG^0
  g^{-1} \cap G^0,G^0})_* \F_{gG^0
  g^{-1} \cap G^0}} \right )
$$

\begin{lemm}
Les deux foncteur précédents $\xymatrix{P\top_{\bullet-H-\lss} \ar@<.6ex>[r]^\a 
  & \mathcal{H}^{cart}_\bullet \ar@<.6ex>[l]^{\b}}$ 
entre faisceaux de Hecke cartésiens sur la tour et faisceaux sur l'espace de
périodes 
sont adjoints l'un
de l'autre : $\a$ est l'adjoint à gauche de $\beta$. 
\end{lemm}
\dem
Cela résulte facilement de la seconde description du foncteur
$\beta$. 
\qed

\begin{rema}
Si $\Pi$ est surjectif alors $\a$ est  fidèle puisque
$\Pi:X_{G^0}\ldrt P$
est un recouvrement étale.
\end{rema}

\begin{prop}\label{DVKET479dsgf2}
Supposons l'existence d'un domaine analytique compact $W\subset
X_{G^0}$ satisfaisant aux hypothèses de la section \ref{SFSjygez26Uetejej3} et tel que de
plus $\Pi_{|W}$ soit un isomorphisme entre $W$ et $\Pi (W)$. Alors
$\beta$ est pleinement fidèle et donc $\a$ et $\beta$ induisent des
équivalences de catégories entre faisceaux sur l'espace de périodes et
faisceaux cartésiens sur la tour.
\end{prop}
\dem 
Grâce à la description donnée dans la section \ref{SFSjygez26Uetejej3} des faisceaux
cartésiens comme recollés de faisceaux en chaque sommet de
``l'immeuble'' $\mathcal{I}$ on vérifie que
$$
\beta (( \F_K)_K)_{|\Pi (W)} = \Pi_* ((\F_{G^0})_{|W})
$$
\qed

\begin{rema}
Dans la proposition précédente pour traiter le cas du site étale il
faut remarquer que si $\Pi_{|W}$ est un isomorphisme alors puisque
$\Pi$ est étale il existe un domaine analytique compact $W'$ contenant
$W$ dans son intérieur et tel que $\Pi_{|W'}$ soit un isomorphisme. 
\end{rema}

\subsection{Rajout d'une donnée de descente}

Soit $F|\Qp$ une extension de degré fini et
$\breve{F}=\widehat{F^{nr}}$. On note $\s\in \text{Aut} (\breve{F})$
le Frobenius de $F^{nr}|F$. Plaçons nous dans le cadre précédent et
supposons que notre tour $K\longmapsto X_K$ est une tour de
$\breve{F}$-espaces analytiques.

Supposons de plus que la tour est munie d'une donnée de descente à $F$
c'est à dire un morphisme entre foncteur de $\mathcal{S}$ vers la
catégorie des $\breve{F}$-espaces analytiques munis d'une action de
$H$
$$
\a: (X_K)_{K\subset G^0} \ldrt (X_K^{(\s)})_{K \subset G^0}
$$
Cette donnée de descente ne sera pas nécessairement supposée
effective. Notons $pr:  (X_K^{(\s)})_{K}\ldrt (X_K)_K$ la projection
de $(X_K)_K\hat{\otimes}_{\breve{F},\s} \breve{F}$ vers $(X_K)_K$.

On peut alors considérer les faisceaux de Hecke sur la tour munis
d'une donnée de descente c'est à dire pour $\bullet \in \{\et,\qetc
\}$ le topos total du diagramme de topos 
$$
\xymatrix{
\mathcal{H}_\bullet ( (X_K)_K)\ar@<.6ex>[r]^{Id}
\ar@<-.6ex>[r]_{pr\circ \a} & \mathcal{H}_\bullet ((X_K)_K)
}
$$
dont les objets  sont les $(\F_K)_K$ munis de morphismes
$H$-équivariants compatibles aux morphismes de la tour $\a^*
\F_K^{(\s )}\ldrt \F_K$.

Notons alors $\,_\s\mathcal{H}_\bullet$ le topos des faisceaux de
Hecke munis d'une donnée de descente et
$\,_\s\mathcal{H}_\bullet^{cart}$ le même en remplaçant faisceaux par
faisceaux cartésiens.

Si $(\F_K)_K\in \,_\s\mathcal{H}_\et$ notons $(\overline{\F}_K)_K$ le
faisceau étendu à la tour $(X_K\otimes_{\breve{F}}
\widehat{\overline{\breve{F}}})_K$ où l'on voit cette dernière tour
comme associée non-plus aux groupes $G$ et $H$ mais à $G$ et $H\times
I_F$ où $I_F=\Gal (\overline{\breve{F}} |\breve{F})$. On a bien
$(\overline{\F}_K)_K \in \mathcal{H}_\et ((X_K\otimes_{\breve{F}}
\widehat{\overline{\breve{F}}})_K)$ (l'action de l'inertie $I_F$ est
lisse). 
Couplé à la donnée de descente cela définit un foncteur
\begin{eqnarray*}
\GG_c (\mathcal{H}_\et\hat{\otimes}_{\breve{F}}
\widehat{\overline{\breve{F}}}, -) : \,_\s \mathcal{H}_\et & \ldrt
& \La [G\times H\times W_F]_{\lss} \\
(\F_K)_K & \longmapsto & \underset{K\subset G^0}{\limi} \GG_c (
X_K\hat{\otimes}\widehat{\overline{\breve{F}}}, \overline{\F}_K)
\end{eqnarray*}
et on note $R\GG_c ( \mathcal{H}_\et\hat{\otimes}_{\breve{F}}
\widehat{\overline{\breve{F}}}, -) :  \,_\s \mathcal{H}_\et \ldrt
\DD^+ ( \La [G\times H\times W_F]_{\lss} )$ le foncteur dérivé
associé.
Il jouit de toutes les propriétés précédentes énoncées dans la section
\ref{DGKSRY246EHENrf} (utiliser toujours le fait qu'un objet injectif
d'un topos total associé à un topos fibré est injectif étage par étage).
\\

Plaçons nous maintenant dans le cadre de la section
\ref{dgskbds46jtjez}. Supposons de plus que l'espace de périodes $P$ soit
muni d'une donnée de descente et que le morphisme des périodes $\Pi$
soit compatible à cette donnée. On obtient alors une équivalence
entre faisceaux de Hecke cartésiens sur la tour et $H$-faisceaux
lisses sur $P$ munis d'une donnée de descente. 
\\

Dans les considérations précédentes on peut également considérer la
variante suivante : remplacer ce qu'on a appelé donnée de descente par
la notion plus forte qui consiste à demander que le morphisme
$\a^*\F_K^{(\s)}\ldrt \F_K$ soit un isomorphisme. Tout ce que l'on
vient de dire s'adapte à cette situation.

\section{Faisceaux de Hecke cartésiens et faisceaux rigides équivariants en
  niveau infini} 
\subsection{Faisceaux lisses sur une tour} 
 
Soient $G$ et $H$ deux groupes profinis et $(\X_K)_{K\subset G}$ une tour équivariante de $H$-schémas formels admissibles quasicompacts. Plus précisément si $\mathcal{S}$ est la catégorie dont les objets sont les sous-groupes ouverts de $G$ et
$$
\forall K,K'\in \mathcal{S}\;\; \Hom_{\mathcal{S}} ( K,K') = \{ g\in K\bc g\;|\; g^{-1} K g\subset K' \}
$$
on suppose que $K\longmapsto \X_K$ définit un foncteur de $\mathcal{S}$ dans la catégorie 
des schémas formels admissibles quasicompacts munis d'une action continue de $H$. On suppose de plus que pour $K'\subset K$ le morphisme $\X_{K'}\ldrt \X_K$ est affine et si $K'\lhd K$ alors
$\X_{K'}^{rig}\ldrt \X_K^{rig}$ est un $K/K'$-torseur étale.

Soit $\X_\infty = \underset{K}{\limp} \X_K$ qui est donc muni d'une action de $G\times H$.

\begin{lemm}
L'action de $G\times H$ sur $\X_\infty$ est continue.
\end{lemm}

Rappelons (théorème \ref{jkfjdgyuzr345Yfg}) qu'il y a une équivalence
de topos
$$
(\X_\infty )\top_{\E-rig-\et} \iso \underset{K}{\limp} (\X_K)\top_{rig-\et}
$$

\begin{prop}\label{sfytarzte647zezv}
Il y a des équivalences de topos 
$$
(\X_\infty)\top_{\E-rig-\et-G\times H-\lss} \simeq (\X_G)\top_{rig-
  \et-H-\lss} 
$$
qui sont équivalents au topos des faisceaux de Hecke quasi-étales cartésiens $H$-lisses sur la tour
$K\longmapsto \X_K^{an}$. 
\end{prop}
\dem
Il s'agit d'une application du théorème \ref{kegjuegt2T35tp}. On
laisse en effet le lecteur vérifier que via l'équivalence $(\X_\infty
)\top_{\E-rig-\et} \iso \underset{K\lhd G}{\limp} (\X_K)\top_{rig-\et}$ les
faisceaux $G\times H$-lisses correspondent aux systèmes projectifs de
$G\times H$-faisceaux lisses. 
\qed

\subsection{Principaux résultats}

Reprenons toutes les notations de la section
\ref{kjsfyqdfgEZTdytzrz4658}. On supposera que $\forall I\subset
\mathcal{I}$ la tour d'espaces analytiques $(W_{I,K})_{K\subset
  G_I\cap G^0}$ est de la forme $(\mathcal{W}_{I,K}^{an})_{K\subset
  G_I\cap G^0}$ où 
$$
(\mathcal{W}_{I,K})_{K\subset G_I\cap G^0}
$$
est une tour de schémas formels admissibles quasicompacts
munis d'une action lisse de $H$ sur $\spf
(\O_{\breve{F}})$
dont les morphismes de transition sont affines.

On supposera de plus que si $I\subset J$ l'inclusion de tours de
domaines analytiques $(W_{J,K})_{K\subset G_J\cap G^0} \hookrightarrow
(W_{I,K})_{K\subset G_I\cap G^0}$ est donnée par une immersion ouverte 
$$
(\mathcal{W}_{J,K})_{K\subset G_J\cap G^0} \hookrightarrow
(\mathcal{W}_{I,K})_{K\subset G_I\cap G^0}
$$
Cela permet de définir un schéma formel $\pi$-adique sans
$\pi$-torsion sur $\spf
(\O_{\breve{F}})$ 
$$
\X_\infty = \bigcup_{z\in \mathcal{I}} \underset{K\subset G_z\cap G^0}{\limp} \mathcal{W}_{z,K}
$$
obtenu par recollement des limites projectives des tours au dessus de
chaque sommet de l'immeuble comme dans l'article \cite{Cellulaire}. 
Ce schéma formel est muni d'une action cellulaire continue de $G\times
H$.

\begin{theo}\label{kdiSFPZ24riMf}
Il y a des équivalences de topos
$$
(\X_\infty)\top_{\E-rig-\et-G\times H-\lss} \simeq
\mathcal{H}^{cart}_{\qetc-H-\lss} \simeq P\top_{\qetc-H-\lss}
$$
entre faisceaux $\E$-rig-étales $G\times H$-équivariants lisses sur
$\X_\infty$,
systèmes de faisceaux de Hecke  cartésiens quasi-étales $H$-lisses sur
la tour et faisceaux quasi-étales $H$-lisses sur l'espace des périodes.
\end{theo}
\dem
Il suffit d'empiler la proposition \ref{sfytarzte647zezv} avec  le
théorème \ref{jytarmsvioizr35zsfjA} et la proposition \ref{DVKET479dsgf2}.
\qed
\\

Soit $\X$ un schéma formel $\pi$-adique sans $\pi$-torsion. 
Si $x:\spf (V)\ldrt \X$ est un point rigide de $\X$ 
et $\F \in\X_{\E-rig-\et}\top$ on peut définir la fibre de $\F$ en $x$
$x^*\F$ qui est un ensemble muni d'une action lisse de $\Gal
(\overline{V}\unp | V\unp)$. 

On a défini dans la section les points analytique d'un schéma formel
$\pi$-adique sans $\pi$-torsion et la notion de spécialisation dans $|\X^{rig}|$.

\begin{defi}
Soit $\X$ un schéma formel $\pi$-adique sans $\pi$-torsion. Un faisceau $\F\in
\X_{\E-rig-\et}\top$ sera dit surconvergent si $\forall x:\spf
(V)\ldrt \X$ , $\forall y:\spf (V')\ldrt \X$ deux points rigides tels
que $x\succ y$ le morphisme de spécialisation 
$$
y^* \F\ldrt x^*\F
$$
est un isomorphisme.
\end{defi}

On peut alors démontrer le théorème suivant.

\begin{theo}
Dans les équivalences de topos précédentes le sous-topos des faisceaux
étales sur l'espace de périodes correspondent au faisceaux surconvergents de $(\X_\infty)\top_{\E-rig-\et-G\times H-\lss}$.
\end{theo}

\begin{theo}\label{kudgtZROOR24Ttg}
Soit $(\F_K)_K$ un système de faisceaux de Hecke étales cartésien sur la
tour et $R\GG_c (\mathcal{H}_\et\hat{\otimes}
\widehat{\overline{\breve{F}}}, (\F_K)_K) \in \DD^+ ( \La[G\times
H\times I_F]_\lss)$ son complexe de cohomologie lisse. Soit $\G \in
\La-(\X_\infty)\top_{\E-rig-\et-G\times H-\lss}$ le faisceau
équivariant lisse associé sur $\X_\infty$.  Il y a alors un
isomorphisme dans la catégorie dérivée équivariante lisse
$$
R\GG_c (\mathcal{H}_\et\hat{\otimes}
\widehat{\overline{\breve{F}}}, (\F_K)_K)
\simeq R\GG_! ( (\X_\infty\hat{\otimes}
\O_{\widehat{\overline{\breve{F}}}})^{rig}, \G)
$$
\end{theo}
\dem
Utiliser le théorème \ref{kdiSFPZ24riMf} couplé aux résolutions de Cech
permettant de calculer les différents complexes de cohomologie à
support compact pour des résolutions flasques, comme  le lemme \ref{FSKDGKZERTY36458HJB}.
\qed

\section{Application aux tours de Lubin-Tate et de Drinfeld}

 \subsection{La correspondance de Jacquet-Langlands locale géométrique}

\begin{theo}\label{QFKSDuqza236gjyza} 
Soit $(\P^{n-1})^{tordu}$ l'espace des périodes de Gross-Hopkins
descendu à
$F$ via la donnée de descente de Rapoport-Zink. Il s'agit de l'espace
analytique de Berkovich associé à une 
variété de Severi-Brauer définie par l'algèbre à division $D$
d'invariant $\frac{1}{n}$ sur $F$. Soit $\Omega\subset \P^{n-1}_{/F}$
l'espace de Drinfeld sur $F$. Il y a une équivalence de topos entre
$D^\times$-faisceaux étales  sur $ (\P^{n-1})^{tordu}$
dont l'action de $D^\times$ est lisse
 et $\GL_n
(F)$-faisceaux étales sur $\Omega$
pour lesquels l'action de $\GL_n (F)$ est lisse
$$
\text{JL} : \left ((\P^{n-1})^{tordu}\right )\top_{\et-D^\times-\lss}
\iso \Omega\top_{\et-GL_n (F)-\lss}
$$
Il en est de même en
remplaçant étale par quasi-étale (i.e. faisceaux sur le site étale de
l'espace rigide plutôt que sur l'espace analytique de Berkovich). 
\end{theo}
\dem
Dans \cite{iso4} on a construit un isomorphisme
$\GL_n (F)\times D^\times$-équivariant entre
schémas formels $\pi$-adiques  sur $\spf (\O_{\widehat{F^{nr}}})$
$$
\X'_\infty \iso \mathcal{Y}'_\infty
$$
où $\X'_\infty$ est un éclaté équivariant du schéma formel $\X_\infty$
construit
dans \cite{Cellulaire} et associé à la tour de Lubin-Tate et
$\mathcal{Y}'_\infty$ est un éclaté d'un schéma formel
$\mathcal{Y}_\infty$ associé à la tour de Drinfeld. 
Par invariance du topos $\E$-rig-étale par éclatements formels
admissibles on a 
$$
(\X'_\infty)\top_{\E-rig-\et-GL_n (F)\times D^\times-\lss} \simeq
(\X_\infty)\top_{\E-rig-\et-GL_n (F)\times D^\times-\lss}$$
et
$$
(\mathcal{Y}'_\infty)\top_{\E-rig-\et-GL_n (F)\times D^\times-\lss} \simeq
(\mathcal{Y}_\infty)\top_{\E-rig-\et-GL_n (F)\times D^\times-\lss} 
$$
D'après le théorème \ref{kdiSFPZ24riMf} 
$$
(\X_\infty)\top_{\E-rig-\et-GL_n (F)\times D^\times-\lss} \simeq (\mathbb{P}^{n-1})\top_{\qetc-D^\times-\lss}
$$
et 
$$
(\mathcal{Y}_\infty)\top_{\E-rig-\et-GL_n (F)\times D^\times-\lss}
\simeq (\Omega\hat{\otimes} \widehat{F^{nr}})\top_{\qetc-GL_n (F)-\lss}
$$
D'où le résultat sur $\widehat{F^{nr}}$
et le site quasi-étale. Le cas du site étale se déduit en remplaçant
le topos $\E$-rig-étale par le topos surconvergent. 
 La descente de
$\widehat{F^{nr}}$ à $F$ ne pose pas de problème une fois que l'on a
vérifié que l'isomorphisme $\X'_\infty \iso \mathcal{Y}'_\infty$ est
compatible à la donnée de descente de Rapoport-Zink.
\qed

\subsection{Comparaison des complexes de cohomologie des deux tours}

\begin{theo}\label{ksdsrpytg79324Fgp}
Soit $(\M_K^{\LT})_{K\subset GL_n (\O_F)}$ la tour de Lubin-Tate
au sens de Rapoport-Zink (une union disjointe indexée par $\Z$ de la
tour de Lubin-Tate usuelle),
formée d'espaces analytiques sur $\breve{F}$, munie d'une action
``horizontale'' de $D^\times$
 d'une action ``verticale'' de $\GL_n (F)$ par correspondances de
 Hecke et d'une donnée de descente à $F$. Soit 
$$
\breve{\pi}^{\LT}_K : \M_K^{\LT} \ldrt \P^{n-1}_{/\breve{F}}
$$
l'application des périodes de Gross-Hopkins. \`A
un faisceau $\F$ de $ \left ((\P^{n-1})^{tordu}\right
)\top_{\et-D^\times-\lss}$ 
 est associé un système de faisceaux ``de Hecke'' $\left
   ((\breve{\pi}_K^{\LT})^* \F\right )_{K\subset GL_n (\O_F)}$ sur la tour de
   Lubin-Tate muni 
de correspondances de Hecke, d'une action de $D^\times$  et
d'une donnée de descente à $F$ compatible à celle
   sur la tour (ces trois actions commutant). Soit le foncteur 
\begin{eqnarray*}
\GG_c (\LT,-) : \La- \left ((\P^{n-1})^{tordu}\right )\top_{\et-D^\times-\lss} &\ldrt
& \La [\GL_n (F)\times D^\times \times W_F]-\text{Mod-lisses} \\
\F & \longmapsto & \underset{K\subset GL_n (\O_F)}{\limi} \GG_c (
\M_K^{\LT}\hat{\otimes} \widehat{\overline{\breve{F}}}, \breve{\pi}_K^{\LT*} \F)
\end{eqnarray*}

Soit $(\M_K^{\D})_{K\subset \O_D^\times}$ la tour de Drinfeld au sens de
Rapoport-Zink (une union disjointe indexée par $\Z$ de la tour
usuelle de revêtements de l'espace $\Omega$ de Drinfeld), formée
d'espaces analytiques sur $\breve{F}$, munie d'une action
``horizontale'' de $\GL_n (F)$, ``verticale'' de $D^\times$ et d'une
donnée de descente à $\breve{F}$. Soit 
$$
\breve{\pi}^\D_K : \M^{\D}_K \ldrt \Omega_{/\breve{F}} 
$$
l'application des périodes sur la tour de Drinfeld. \`A un faisceau
$\G$ de $\Omega\top_{\et-GL_n (F)-\lss}$ est associé un système de faisceaux de
Hecke muni d'une donnée de descente à $F$ sur la tour de Drinfeld,
$(\breve{\pi}^{\D*}_K \G)_{K\subset \O_D^\times}$. Soit le foncteur
\begin{eqnarray*}
\GG_c (\D,-) : \La-\Omega\top_{\et-GL_n (F)-\lss} &\ldrt & \La [\GL_n
(F)\times D^\times \times W_F]-\text{Mod-lisses} \\
\G  & \longmapsto &  \underset{K\subset GL_n (\O_F)}{\limi} \GG_c 
(\M_K^{\D}\hat{\otimes} \widehat{\overline{\breve{F}}}, \breve{\pi}_K^{\D*} \G)
\end{eqnarray*}
Il y a alors un isomorphisme de foncteurs à valeurs dans $\DD^+ (\La[\GL_n (F)\times D^\times
\times W_F]_{\lss})$ 
$$
R\GG_c (\D,-)\circ \text{JL} \iso R\GG_c (\LT,-)
$$
En particulier si $\underline{\La}$ désigne le faisceau constant,
  $\text{JL} (\underline{\La})=\underline{\La}$ et donc il y a un
  isomorphisme entre la cohomologie à support compact à coefficients
  constants des tours de Lubin-Tate et de Drinfeld.
\end{theo}
\dem
Il s'agit d'une application du théorème \ref{kudgtZROOR24Ttg} 
et du théorème principal de \cite{iso4} 
comme
dans la démonstration précédente.
\qed

\bibliographystyle{plain}
\bibliography{biblio}
\end{document}